\def\DateTime{April 6, 2021}
\def\Version{Version $1.0$}

\def\no{\if01}
\def\iftwelvept{\no}

\iftwelvept
\documentclass[leqno,12pt]{amsart}
\usepackage{fullpage}

\else
\documentclass[leqno]{amsart}

\fi

\usepackage{amscd}
\usepackage{amsmath}
\usepackage{amssymb}
\usepackage{amsfonts}
\usepackage{latexsym}
\usepackage{color}

\usepackage{colonequals}

\usepackage{lmodern}
\usepackage[T1]{fontenc}

\theoremstyle{plain}
\newtheorem{theorem}{Theorem}[section]
\newtheorem{proposition}[theorem]{Proposition}
\newtheorem{lemma}[theorem]{Lemma}

\theoremstyle{definition}
\newtheorem{definition}[theorem]{Definition}

\theoremstyle{remark}
\newtheorem{remark}[theorem]{Remark}

\newenvironment{pf}{\begin{proof}}{\end{proof}}
\makeatletter
\renewcommand{\theequation}{%
\thesection.\arabic{equation}}
\@addtoreset{equation}{section}
\makeatother
\renewcommand{\theclaim}{\arabic{claim}}

\newcommand{\ZZ}{{\bf{Z}}}

\newcommand{\RR}{{\bf{R}}}
\newcommand{\CC}{{\bf{C}}}
\newcommand{\PP}{{\bf{P}}}

\newcommand{\LL}{{\mathbb{L}_{K3}}}
\newcommand{\KK}{{\mathbb{K}}}
\newcommand{\LAM}{{\bf\Lambda}}
\newcommand{\HH}{{\mathfrak
{H}}}
\newcommand{\EE}{{\mathbb{E}}}

\newcommand{\UU}{{\mathbb{U}}}

\newcommand{\SL}{{\operatorname{SL}}}

\renewcommand{\setminus}{\smallsetminus}
\newcommand{\Ima}{{\operatorname{Im}}}

\newcommand{\dbar}{{\bar{d}}}
\newcommand{\Km}{{\operatorname{Km}}}

\newcommand{\rank}{{\operatorname{rank}}}

\newcommand{\Hom}{{\operatorname{Hom}}}

\newcommand{\rest}[2]{\left.{#1}\right\vert_{{#2}}}

\renewcommand{\d}{\delta}

\renewcommand{\l}{\lambda}

\renewcommand{\r}{\rho}

\renewcommand{\t}{\tau}

\newcommand{\D}{\Delta}

\renewcommand{\O}{\Omega}

\newcommand{\Ccal}{{\mathcal C}}
\newcommand{\Dcal}{{\mathcal D}}

\newcommand{\Hcal}{{\mathcal H}}

\newcommand{\Mcal}{{\mathcal M}}

\newcommand{\Dfrak}{{\mathfrak D}}
\newcommand{\Hfrak}{{\mathfrak H}}

\newcommand{\ebf}{{\bf e}}
\newcommand{\fbf}{{\bf f}}

\begin{document}

\title[]
{$j$-invariant and Borcherds $\Phi$-function}
\author{Shu Kawaguchi}
\address{
Department of Mathematical Sciences, 
Doshisha University, 
Kyotanabe Kyoto 610-0394, 
Japan}
\email{kawaguch@mail.doshisha.ac.jp}
\author{Shigeru Mukai}
\address{
Research Institute for Mathematical Sciences, 
Kyoto University, 
Kyoto 606-8502, 
Japan}
\email{mukai@kurims.kyoto-u.ac.jp}
\author{Ken-Ichi Yoshikawa}
\address{
Department of Mathematics, 
Faculty of Science,
Kyoto University, 
Kyoto 606-8502, 
Japan}
\email{yosikawa@math.kyoto-u.ac.jp}
\date{\DateTime, (\Version)}

\begin{abstract}
We give a formula that relates the difference of the $j$-invariants  
with the Borcherds $\Phi$-function, an automorphic form on the period domain 
for Enriques surfaces characterizing the discriminant divisor. 
\end{abstract}

\maketitle


\section{Introduction}
\label{sect:introduction}
The $j$-invariant $j(\tau)$ is the ${\rm SL}_{2}({\bf Z})$-invariant holomorphic function on the complex upper half-plane ${\mathfrak H}$
with Fourier series expansion at the cusp $j(\tau)=e^{-2\pi i\tau}+744+196884\,e^{2\pi i\tau}+\cdots$,
which induces an isomorphism from the moduli space of elliptic curves to ${\bf C}$. 
The $j$-invariant is fundamental in many branches of mathematics such as the elliptic function theory, number theory, 
the theory of automorphic forms, and monstrous moonshine... 
Besides the $j$-invariant itself, it has been discovered that the difference $j(\tau)-j(\tau')$ enjoys beautiful properties as well: 
Gross--Zagier's result \cite{GZ} on singular moduli and the denominator formula for the monster Lie algebra 
(see Borcherds \cite{Borcherds92} and \cite{Borcherds95}), for example. 
In this paper, we show that the difference $j(\tau)-j(\tau')$ is closely related to the {\em  Borcherds $\Phi$-function} \cite{Borcherds96}, 
a remarkable automorphic form on the period domain for Enriques surfaces.  

For $\tau\in{\mathfrak H}$, we set $E_{\tau}={\bf C}/{\bf Z}+\tau{\bf Z}$. 
Let ${\mathfrak D}$ be the $\Gamma(2)\times \Gamma(2)$-orbit of 
the diagonal locus of ${\mathfrak H}\times{\mathfrak H}$, where 
$\Gamma(2) \subset \SL_2(\ZZ)$ is the principal  congruence subgroup of level $2$. Then, for any 
$(\tau, \tau') \in {\mathfrak H}\times{\mathfrak H}\setminus{\mathfrak D}$, 
one can construct $15$ fixed-point-free involutions on ${\rm Km}(E_{\tau}\times E_{\tau'})$, the Kummer surface of product type associated to $E_{\tau}\times E_{\tau'}$, 
which induce $15$ distinct Enriques surfaces: Indeed, if $(\tau, \tau')\in{\mathfrak H}\times{\mathfrak H}$ is very general, then, up to conjugacy, 
these $15$ involutions are the only 
fixed-point-free involutions on ${\rm Km}(E_{\tau}\times E_{\tau'})$  (see \cite{Kondo86}, \cite{Mukai10}, \cite{MukaiNamikawa84}, \cite{Ohashi07}). 
By \cite{Mukai10}, \cite{Ohashi07}, there is a one-to-one correspondence between the set of these $15$ conjugacy classes of involutions and the set of non-zero elements of the discriminant group  $A_{\mathbb K}$ 
of the lattice ${\mathbb K} = \UU(2) \oplus \UU(2)$
(see Sect.~\ref{sec:lattice} for the notation about lattices). 
Recall that the discriminant form of ${\mathbb K}$ takes its values in ${\mathbf Z}/2{\mathbf Z}$. 
According to the corresponding value of the discriminant form of $A_{\mathbb K}$, 
the $15$ conjugacy classes of involutions are divided into the $6$ {\em odd} involutions and the $9$ {\em even} involutions. 

Let $\LAM = {\mathbb U}(2) \oplus {\mathbb U} \oplus {\mathbb E}_{8}(2)$ be the Enriques lattice. 
Let  $\Omega_{\LAM}^+$ denote the period domain for Enriques surfaces (cf.~\eqref{eqn:def:omega:LAM}). 
In \cite{Borcherds96}, Borcherds discovered a remarkable automorphic form $\Phi$, called the Borcherds $\Phi$-function, of weight $4$ on $\Omega_{\LAM}^+$ vanishing exactly 
on the discriminant locus. 
For $\ell \in \{1, 2\}$, 
once a primitive isotropic vector ${\bf e}_{\ell} \in {\mathbb U}(\ell)$ and an isomorphism of $\Omega_{\LAM}^+$ 
with the tube domain $\mathbb{M}_\ell \otimes \RR + i C^+_{\mathbb{M}_\ell}$ (see \eqref{eqn:tube:domain:level2}) are fixed, 
where ${\mathbb M}_{\ell} = {\bf e}_{\ell}^{\perp}/{\bf e}_{\ell} = {\mathbb U}(2/\ell) \oplus {\mathbb E}_{8}(2)$,  
then $\Phi$ is identified with the holomorphic function given by an explicit infinite product,
denoted by $\Phi_\ell$, on $\mathbb{M}_\ell \otimes \RR + i C^+_{\mathbb{M}_\ell}$ 
(see Sect.~\ref{subsec:Borcherds:Phi} for details). 

For each $\gamma \in A_{\KK}\setminus\{0\}$, let $\iota_\gamma$ be the attached fixed-point free involution on ${\rm Km}(E_{\tau}\times E_{\tau'})$. One can naturally associate an integer $\ell(\gamma) \in \{1, 2\}$, and 
a period map $\varphi_\gamma\colon {\mathfrak H}\times{\mathfrak H} \to \mathbb{M}_{\ell(\gamma)} \otimes \RR + i C^+_{\mathbb{M}_{\ell(\gamma)}}$ for the family of Kummer surfaces ${\rm Km}(E_{\tau}\times E_{\tau'})$ over ${\mathfrak H}\times{\mathfrak H}$ by fixing a normalized marking (see Definition~\ref{def:normalized:marking}). 
We set $\Phi_\gamma = \Phi_{\ell(\gamma)} \circ \varphi_\gamma$. 
Then $\Phi_{\gamma}^{2}$ is independent of the choice of a normalized marking, and is 
an automorphic form on ${\mathfrak H}\times{\mathfrak H}$ of weight 
$8$ for $\Gamma(2) \times \Gamma(2)$. 
(See Sect.~\ref{sect:parity:invol} and Sect.~\ref{subsec:actions:period:map}.)

The main result of this paper states that the difference $j(\tau)-j(\tau')$ 
is related to the Borcherds $\Phi$-function in the following simple way.

\begin{theorem}
\label{thm:MainTheorem}
For any $(\tau, \tau') \in {\mathfrak H}\times{\mathfrak H} \setminus{\mathfrak D}$,  
\begin{equation}
\label{eqn:Main:Theorem}
2^{-96}
\left( j(\tau) - j(\tau') \right)^{12} 
=
\frac
{\prod_{\gamma\,{\rm odd}}\Phi_{\gamma}(\tau,\tau')^{6}}
{\prod_{\gamma\,{\rm even}}\Phi_{\gamma}(\tau,\tau')^{4}}. 
\end{equation}
\end{theorem}

We will prove Theorem~\ref{thm:MainTheorem} by first showing \eqref{eqn:Main:Theorem} up to a constant and then showing that the constant is equal to $2^{-96}$.

To prove \eqref{thm:MainTheorem} up to a constant, 
since the logarithm of the absolute value of each side of \eqref{thm:MainTheorem} is shown to be an ${\rm SL}_{2}({\bf Z})\times{\rm SL}_{2}({\bf Z})$-invariant 
pluriharmonic function on ${\mathfrak H}\times{\mathfrak H}$, it suffices to compare their singularities. 
This will be done by analyzing the period map $\varphi_{\gamma}$ and its extension to the Baily-Borel compactification, in particular its intersection with the discriminant locus and the boundary locus of the moduli space of Enriques surfaces.
Using lattice theory, we will show that such a geometric property of $\varphi_{\gamma}$ is determined by the parity of $\gamma$.

Determining the constant in \eqref{thm:MainTheorem} is a different and somewhat a more delicate problem. 
We will do this by computing the leading terms of  the denominator and the numerator of the right-hand side of \eqref{eqn:Main:Theorem} near the cusp $(+i\infty,+i\infty)$.
For the computation of ${\prod_{\gamma\,{\rm even}}\Phi_{\gamma}(\tau,\tau')}^4$, we use the formula in \cite[Cor.~7.6]{KawaguchiMukaiYoshikawa18}, 
which is a consequence of an algebraic expression of $\Phi$. For the computation of ${\prod_{\gamma\,{\rm odd}}\Phi_{\gamma}(\tau,\tau')}^6$, we carefully study the Borcherds products 
with respect to the $0$-dimensional cusps of levels one and two.  For the left-hand side of \eqref{eqn:Main:Theorem}, we use the denominator formula for the monster Lie algebra
\cite{Borcherds92} (see also \eqref{eqn:denominator:formula:monster:Lie:algebra}). All of these enable us to obtain the constant $2^{-96}$ on the left-hand side of \eqref{eqn:Main:Theorem}. 

For the results on the the discriminant locus and the boundary locus of the moduli space of Enriques surfaces, 
which we obtain in the course of the proof of Theorem~\ref{thm:MainTheorem} and may be of independent interest, see e.g. Lemma~\ref{lemma:boundary:component}, Propositions~\ref{prop:discriminant:fixed:point} and \ref{prop:characterization:zero:Kondo-Mukai}.

The organization of this paper is as follows. 
In Sect.~\ref{sec:Enriques:Phi}, we recall the moduli space of Enriques surfaces and the Borcherds $\Phi$-function. 
In Sect.~\ref{sec:involution}, we recall the period mapping $\varpi_{\gamma}$ for the Enriques surfaces ${\rm Km}(E_{\tau}\times E_{\tau'})/\iota_{\gamma}$,
prove the automorphy of $\Phi_{\gamma}$ and
study the intersection of the period mapping $\varpi_{\gamma}$ with the discriminant locus of $\Omega_{\LAM}^{+}$. 
In Sect.~\ref{sec:L}, we give an explicit formula for the denominator of the right hand side of \eqref{eqn:Main:Theorem}.
In Sect.~\ref{sec:KM}, we recall involutions of odd type and compute $dd^{c}\log(\cdot)$ of the Petersson norm of the numerator of the right hand side of 
\eqref{eqn:Main:Theorem} as a current on the second symmetric product of the compactified modular curve.
In Sect.~\ref{sec:comparison}, we compute the leading term of $\Phi_{\gamma}$ near $(+i\infty, +i\infty)$ 
for all odd $\gamma\in A_{\mathbb K}\setminus\{0\}$. 
In Sect.~\ref{sec:proof}, we prove Theorem~\ref{thm:MainTheorem}. 
In Sect.~\ref{sec:comparison:2}, we consider automorphic forms arising from non-Enriques involutions. 
We also list some open problems. 
An appendix is accompanied to give some technical results on lattices. 

\smallskip
{\bf Acknowledgements }
The first named author is partially supported by JSPS KAKENHI 18H01114 and 20H00111. 
The second named author is partially supported by JSPS KAKENHI 15H05738 and 16H06335. 
The third named author is partially supported by JSPS KAKENHI 16H03935 and 16H06335.  
The first and the third named authors thank Prof. Vincent Maillot for helpful discussions, 
and the first named author thanks Prof. Damian R\"{o}ssler for helpful discussions. 
The authors thank Prof. Noriko Yui for bringing the references \cite{LiYang20}, \cite{YangYinYu18}, \cite{YangYin19} 
to their attention.

\section{Enriques surfaces and the Borcherds $\Phi$-function}
\label{sec:Enriques:Phi}
 
\subsection{Lattices}
\label{sec:lattice}

A lattice $L$ is a pair of a free ${\bf Z}$-module of rank $r$ and a non-degenerate, integral, symmetric bilinear form 
$\langle\cdot,\cdot\rangle_{L}$ on it. 
When there is no possibility of confusion, we often write $\langle\cdot,\cdot\rangle$ for $\langle\cdot,\cdot\rangle_{L}$. 
A lattice $L$ is said to be even if $x^{2} \colonequals \langle x,x\rangle_{L} \equiv 0$ mod $2$ for any $x \in L$. 
For any integer $m$, we denote by $L(m)$ the lattice $L$ equipped with the rescaled bilinear form $m\langle\cdot,\cdot\rangle_{L}$.
The group of isometries of $L$ is denoted by $O(L)$.
The set of roots of $L$ is defined by $\Delta_{L} \colonequals \{d\in L \mid \langle d,d\rangle_{L}=-2\}$.
Let $L^{\lor} = \Hom(L, {\bf Z})\subset L\otimes{\bf Q}$ be the dual lattice of $L$.
The discriminant group of $L$ is the finite abelian group $A_{L} \colonequals L^{\lor}/L$. 
For a sublattice $S \subset L$, we write 
$S^{\perp_L}$ (and $S^\perp$ if no confusion is likely) for the orthogonal complement of $S$ in $L$. 

Let $L$ be an even lattice. For $\lambda\in L^{\lor}$, we write $\bar{\lambda} \colonequals \lambda+L\in A_{L}$. 
The discriminant form and the discriminant bilinear from on $A_{L}$ are denoted by $q_{L}$ and $b_{L}$, respectively.
We set ${\bf F}_{2} \colonequals {\bf Z}/2{\bf Z}$. If $A_{L}\cong{\bf F}_{2}^{\oplus l}$ for some $l\in{\bf Z}_{\geq0}$, then $L$ is said to be $2$-elementary.
For an even $2$-elementary lattice $L$, we define its parity 
as $\delta(L)\colonequals0$ if $q_{L}$ is ${\bf Z}/2{\bf Z}$-valued and $\delta(L)\colonequals1$ otherwise. 
We refer to \cite{Nikulin80} for more about lattices and discriminant forms.

In this paper, the following lattices are important.
Let ${\mathbb U}$ be the hyperbolic plane, i.e., the even unimodular lattice of signature $(1,1)$ 
and let ${\mathbb E}_{8}$ be the {\em negative-definite} even unimodular lattice of rank $8$. 
The $K3$ lattice is the even unimodular lattice
$$
{\mathbb L}_{K3}
 \colonequals 
{\mathbb U}\oplus{\mathbb U}\oplus{\mathbb U}\oplus
{\mathbb E}_{8}\oplus{\mathbb E}_{8}. 
$$
Then ${\mathbb L}_{K3}$ has rank $22$ with signature $(3, 19)$. 
The Enriques lattice is the even $2$-elementary lattice of rank $12$ defined as 
\[
  \LAM  \colonequals  \UU(2) \oplus \UU \oplus \EE_8(2).
\]
Then $\LAM$ has 
signature $(2, 10)$ and discriminant group $A_{\LAM} \cong {\bf F}_{2}^{\oplus 10}$.  
We fix a primitive embedding $\LAM \subset \LL$. 
Then $\LAM^{\perp_{\LL}} \cong \UU(2)\oplus\EE_8(2)$ (see \cite[VIII, \S19]{BPV84}).

\subsection{The moduli space of Enriques surfaces and the Baily--Borel compactification}
\label{sect:moduli:Enriques:Baily:Borel}

Recall that a compact connected complex surface is a $K3$ surface if it is simply connected and has trivial canonical bundle.
For a $K3$ surface $X$, the second cohomology group $H^2(X, \ZZ)$ endowed with the cup-product is isometric to $\LL$ (see \cite[p.241]{BPV84}). 
The N\'eron--Severi and transcendental lattices of $X$ are defined as  
${\rm NS}_{X} \colonequals H^{2}(X,{\bf Z})\cap H^{1,1}(X,{\bf R})$ and 
$T_{X} \colonequals {\rm NS}_{X}^{\perp_{H^{2}(X,{\bf Z})}}$, respectively.

A compact connected complex surface is an Enriques surface if it is not simply connected and its universal covering is a $K3$ surface. 
Let $Y$ be an Enriques surface with universal covering $K3$ surface $X$. Then $\pi_{1}(Y)\cong{\mathbf Z}/2{\mathbf Z}$ and 
the non-trivial covering transformation $\iota_{Y}\colon X\to X$ is then an anti-symplectic fixed-point-free involution. 
We set 
$$
H^{2}(X,{\bf Z})_{\pm}
 \colonequals 
\{l\in H^{2}(X,{\bf Z})\mid\iota_{Y}^{*}(l)=\pm l\}.
$$
By e.g.  \cite[VIII, \S19]{BPV84}, 
there is an isometry of lattices $\alpha\colon H^{2}(X,{\bf Z})\cong{\mathbb L}_{K3}$ such that 
\begin{equation}
\label{eqn:alpha:+:-}
\alpha(H^{2}(X,{\bf Z})_{+})
=
\LAM^{\perp_{\LL}},
\qquad
\alpha(H^{2}(X,{\bf Z})_{-})
=
\LAM.
\end{equation}
Such an isometry $\alpha$ is called a marking, and 
the pair $(Y,\alpha)$ satisfying \eqref{eqn:alpha:+:-} is called a marked Enriques surface. 

We set
\begin{equation}
\label{eqn:def:omega:LAM}
\Omega_{\LAM}
 \colonequals 
\{[\omega]\in{\bf P}(\LAM\otimes{\bf C})\mid 
\langle\omega,\omega\rangle=0,\,\langle\omega,\bar{\omega}\rangle>0\}.
\end{equation}
Then $\Omega_{\LAM}$ consists of two disjoint 
connected complex manifolds of dimension $10$, both of which is a bounded symmetric domain of type ${\rm IV}$. 
We fix one connected component $\Omega_{\LAM}^{+}$ of $\Omega_{\LAM}$. 
On $\Omega_{\LAM}$ acts $O(\LAM)$ projectively. Let $O^+(\LAM)\subset O(\LAM)$ be the subgroup of index $2$ preserving $\O^+_\LAM$. 
Then $O^+(\LAM)$ acts on $\O_\LAM^+$ properly discontinuously.
The period domain for Enriques surfaces is defined as  
$$
{\mathcal M} \colonequals \Omega_{\LAM}^{+}/O^{+}(\LAM)=\Omega_{\LAM}/O(\LAM). 
$$
The period of a marked Enriques surface $(Y,\alpha)$ is defined as 
\begin{equation}
\label{eqn:definition:period}
\varpi(Y,\alpha)
 \colonequals 
[\alpha(H^{0}(X,\Omega^{2}_{X}))]
\in
\Omega_{\LAM}^{+}
\end{equation}
and the period of an Enriques surface $Y$ is defined 
as  the $O^{+}(\LAM)$-orbit of $\varpi(Y,\alpha)$:
$$
\overline{\varpi}(Y) \colonequals [\varpi(Y,\alpha)]\in{\mathcal M}.
$$

For $d \in \LAM \otimes \RR$,  
we set 
\begin{equation}
\label{eqn:H:d}
H_d\colonequals \{\omega \in\Omega_{\LAM}^{+}\mid \langle\omega, d\rangle = 0\}.
\end{equation}
The discriminant locus of $\Omega^{+}_\LAM$ is the $O^{+}(\LAM)$-invariant divisor
$
{\mathcal H} 
 \colonequals 
\bigcup_{d\in\Delta_{\LAM}}H_d.
$
Set 
\[
\mathcal D \colonequals {\mathcal H}/O^{+}(\LAM). 
\]
Then $\overline{\varpi}(Y)\not\in{\mathcal D}$ for any Enriques surface $Y$.
Via the period map, the coarse moduli space of Enriques surfaces is isomorphic to the 
analytic space (\cite[VIII, \S19]{BPV84})
$$
(\Omega_{\LAM}^{+}\setminus{\mathcal H})/O^{+}(\LAM)
={\mathcal M}\setminus{\mathcal D}.
$$ 

Let $\Mcal^*$ be the Baily--Borel compactification of $\Mcal$. 
By Sterk~\cite[Props.~4.5, 4.6, 4.7]{Sterk}, the boundary $\Mcal^* \setminus \Mcal$ consists of  two $1$-dimensional components,
one of which 
is isomorphic to the modular curve $X(1)  \colonequals  (\SL_{2}({\bf Z})\backslash \Hfrak)^*$
and the other 
to the curve $X^1(2)  \colonequals  (\Gamma^1(2)\backslash \Hfrak)^*$.
Here $\Gamma^{1}(2) =\left\{\left. \binom{a\,b}{c\,d} \in \SL_2(\ZZ)\;\right\vert\; b \equiv 0\; \text{mod $2$}\right\}$ 
and the asterisk $*$ denotes the Baily--Borel  compactification. 
By slight abuse of notation, the components of ${\mathcal M}^{*}\setminus{\mathcal M}$ corresponding to $X(1)$ and $X^{1}(2)$ are denoted by the same symbols.  
Further, $\Mcal^* \setminus \Mcal$ has two $0$-dimensional cusps: 
The curves $X(1)$ and $X^1(2)$ intersect at one point, and this point gives one $0$-dimensional cusp;  
The other $0$-dimensional cusp lies on $X^1(2) \setminus X(1)$. 
 
By \cite[Sect.~2]{Scattone87}, 
from lattice-theoretical terms,  
the $1$-dimensional components of $\Mcal^* \setminus \Mcal$ correspond to the $O(\LAM)$-orbits of the
primitive totally isotropic sublattices of rank $2$ of $\LAM$, 
and the $0$-dimensional cusps of $\Mcal^* \setminus \Mcal$ correspond to 
the $O(\LAM)$-orbits of the primitive isotropic vectors of $\LAM$. 
The above results of Sterk in particular say 
that, up to the $O(\LAM)$-action, there are exactly two distinct such sublattices of rank $2$ of $\LAM$ and 
two distinct primitive isotropic vectors of $\LAM$. 
The following lemma gives their explicit representatives. 

\begin{lemma}
\label{lemma:isotoropic:sublattices}
Let ${\mathbb I}_{2,9}$ be an odd unimodular lattice of signature $(2,9)$. 
\begin{enumerate}
\item
For any root $d \in \Delta_{\LAM}$, we have $d^{\perp_{\LAM}}\cong{\mathbb I}_{2,9}(2)$. 
Let $\{ {\bf e}_1, {\bf e}_2, \ldots, {\bf e}_{11} \}$ be a 
free basis of ${\mathbb I}_{2,9}(2)$ with Gram matrix $2\, {\rm diag}({\bf 1}_2, -{\bf 1}_9)$. 
Let $F_1$ be the sublattice of $\LAM$ that corresponds to $\ZZ ({\bf e}_1 + {\bf e}_3) + \ZZ ({\bf e}_2 + {\bf e}_4)$ via $d^{\perp_{\LAM}}\cong{\mathbb I}_{2,9}(2)$. Then $F_1$ is a primitive totally isotropic sublattice of rank $2$ of $\LAM$, and corresponds to~$X(1)$. 
\item
Let ${\bf e}, {\bf e}'$ be the standard free basis of the left lattice $\UU(2)$ of $\LAM$ with Gram matrix $\binom{0\,2}{2\,0}$, 
and let ${\bf f}, {\bf f}'$  be the standard free basis of the middle lattice $\UU$ of $\LAM$ with Gram matrix $\binom{0\,1}{1\,0}$.  
We set $F_2 \colonequals \ZZ {\bf e} + \ZZ {\bf f}$. Then $F_2$ is a primitive totally isotropic sublattices of rank $2$ of $\LAM$, and corresponds to $X^1(2)$. 
\end{enumerate}
\end{lemma}

\begin{pf}
By \cite[Th.~3.6.2]{Nikulin80}, up to isometries of lattices, an even indefinite $2$-elementary lattice is determined by 
its signature, the rank of its  discriminant group, and its parity. It follows that 
$d^{\perp_{\LAM}}\cong{\mathbb I}_{2,9}(2)$ for any $d\in\Delta_{\LAM}$. 
Recall that $X(1)$ has one $0$-dimensional cusp and $X^1(2)$ has two $0$-dimensional cusps.  
Since ${\bf e}$ and ${\bf f}$ are primitive vectors of $F_2$ with 
$\langle {\bf e}, \LAM \rangle = 2 \ZZ$ and $\langle {\bf f}, \LAM \rangle = \ZZ$, 
the above results of Sterk imply that $F_2$ corresponds to $X^1(2)$, so $F_1$  corresponds to $X(1)$. 
\end{pf}

For a lattice $L$ of signature $(2, r-2)$, we define $\Omega_{L}$, $\Omega_{L}^{+}$, $O(L)$, $O^{+}(L)$, and $H_d$ in the same way as those for $\LAM$, 
and we set ${\mathcal M}_{L} \colonequals \Omega_{L}^{+}/O^{+}(L)={\mathcal M}_{L}/O(L)$.
The Baily--Borel  compactification of ${\mathcal M}_{L}$ is denoted by ${\mathcal M}_{L}^{*}$.
Recall that for any linearly independent vectors $d_{1},\ldots,d_{k} \in L \otimes \RR$, we have 
\begin{equation}
\label{eqn:Hd:d2}
H_{d_{1}}\cap\cdots\cap H_{d_{k}} \neq \emptyset \quad\Longleftrightarrow\quad (d_{i},d_{j})_{1\leq i,j\leq k}\,\, \text{is a negative-definite matrix}. 
\end{equation}

For a primitive isotropic vector ${\bf v}\in L$, the {\em level} of ${\bf v}$ is defined as 
the positive integer $\ell$ with $\langle{\bf v}, L\rangle  = \ell {\bf Z}$. In this case, we say that 
the $0$-dimensional cusp of ${\mathcal M}_{L}^{*}$ corresponding ${\bf v}$ has level $\ell$. 
For example, in the case of $\LAM$, 
by Sterk~\cite[Props.~4.5,~4.6,~4.7]{Sterk} (see also the proof of Lemma~\ref{lemma:isotoropic:sublattices}), 
one of the two $0$-dimensional cusps of ${\mathcal M}^{*}$, which lies on $X^1(2) \setminus X(1)$, has level $1$, 
and the other $0$-dimensional cusp, which is the intersection of $X^1(2)$ and $X(1)$, has level $2$. 

\begin{lemma}
\label{lemma:boundary:component}
Let ${\mathcal D}^{*}$ be the closure of ${\mathcal D}$ in ${\mathcal M}^{*}$.
Then ${\mathcal D}^{*}\setminus{\mathcal D}=X(1)$.
\end{lemma}

\begin{pf}
We fix a root $d\in\Delta_{\LAM}$. 
To simplify the notation, we write $d^\perp$ for $d^{\perp_{\LAM}} \subset \LAM$. 
We choose the connected component 
$\Omega_{d^{\perp}}^+$ of $\Omega_{d^{\perp}}$ 
such that $\Omega_{d^{\perp}}^+ = H_d \subset \Omega_{\LAM}^+$. 

\medskip
{\sl Step 1.}\;
We first relate ${\mathcal M}_{d^{\perp}}^{*}\setminus{\mathcal M}_{d^{\perp}}$, 
${\mathcal D}^{*}\setminus{\mathcal D}$,  
and ${\mathcal M}^{*}\setminus {\mathcal M}$. 
Let $O^{+}(\LAM)_{d} \colonequals \{g \in O^{+}(\LAM) \mid 
g(d) = d\}$ be the stabilizer of $d$ in $O^{+}(\LAM)$. 
Then both $O^{+}(\LAM)_{d}$ and $O^{+}(d^{\perp})$ act on $\Omega_{d^{\perp}}^+$. 
We will show in the appendix (see Lemma~\ref{lemma:appendix:surjectivity}) that 
the restriction map 
$O^{+}(\LAM)_{d}\ni g\mapsto \rest{g}{d^{\perp}}\in O^{+}(d^{\perp})$ 
is surjective. It follows that  
${\mathcal M}_{d^{\perp}} \colonequals \Omega_{d^{\perp}}^+/O^{+}(d^{\perp}) = 
\Omega_{d^{\perp}}^+/O^{+}(\LAM)_{d}$. 

We have the natural surjective map $\Omega_{d^{\perp}}^+/O^{+}(\LAM)_{d} \to (O^{+}(\LAM)\cdot\Omega_{d^{\perp}}^{+})/O^{+}(\LAM)$. 
Since $O^{+}(\LAM)$ acts on $\Delta_\LAM$ transitively (see \cite[Remark~3.6]{Sterk}), we have $O^{+}(\LAM)\cdot\Omega_{d^{\perp}}^{+} = O^{+}(\LAM)\cdot H_d = \mathcal{H}$. 
Thus $ (O^{+}(\LAM)\cdot\Omega_{d^{\perp}}^{+})/O^{+}(\LAM) = \mathcal{D}$. We have obtained 
the surjective morphism ${\mathcal M}_{d^{\perp}} \to \mathcal{D}$. 
By the description of the Baily--Borel compactification, we have the extended surjective morphisms  
${\mathcal M}_{d^{\perp}}^* \to \mathcal{D}^*$ and 
${\mathcal M}_{d^{\perp}}^{*}\setminus{\mathcal M}_{d^{\perp}}
\to\mathcal D^{*}\setminus \mathcal D$.

The inclusion $\Omega_{d^{\perp}}^+ = H_d \subset \Omega_{\LAM}^+$ 
induces the injective morphisms $\mathcal{D} \hookrightarrow {\mathcal M}$, 
$\mathcal{D}^* \hookrightarrow {\mathcal M}^*$, and 
$\mathcal D^{*}\setminus \mathcal D
\hookrightarrow
{\mathcal M}^{*}\setminus{\mathcal M}$.
To conclude, we have the morphisms 
\begin{equation}
\label{eqn:inclusions:boundary:components}
{\mathcal M}_{d^{\perp}}^{*}\setminus{\mathcal M}_{d^{\perp}}
\to
\mathcal D^{*}\setminus \mathcal D
\hookrightarrow
{\mathcal M}^{*}\setminus{\mathcal M},
\end{equation}
where the left morphism is surjective, and the right morphism is injective. 

\medskip
{\sl Step 2.}\;
We show that $D^{*}\setminus \mathcal D$ is irreducible. Indeed, 
by Lemma~\ref{lemma:isotoropic:sublattices}~(1), we have 
an isometry $d^{\perp}\cong{\mathbb I}_{2,9}(2)$, with which 
we identify $
{\mathcal M}_{{\mathbb I}_{2,9}(2)}^{*}\setminus{\mathcal M}_{{\mathbb I}_{2,9}(2)}
= 
{\mathcal M}_{d^{\perp}}^{*}\setminus{\mathcal M}_{d^{\perp}}$. 

By \cite[Prop.~1.17.1]{Nikulin80}, ${\mathbb I}_{2,9}(2)$ has a unique primitive totally isotropic sublattice of rank $2$ up to $O({\mathbb I}_{2,9}(2))$. 
It follows from \cite[Sect.~2.1]{Scattone87} that ${\mathcal M}_{{\mathbb I}_{2,9}(2)}^{*}\setminus{\mathcal M}_{{\mathbb I}_{2,9}(2)}$ is irreducible. 
Since the morphism ${\mathcal M}_{d^{\perp}}^{*}\setminus{\mathcal M}_{d^{\perp}}
\to
\mathcal D^{*}\setminus \mathcal D$ in \eqref{eqn:inclusions:boundary:components} is surjective, 
$D^{*}\setminus \mathcal D$ is irreducible.

\medskip
{\sl Step 3.}\;
We show that $D^{*}\setminus \mathcal D = X(1)$ in ${\mathcal M}^{*}$. 
Let $F\subset d^{\perp}$ be a primitive totally isotropic sublattice of rank $2$. 
We set $N(F) = \{g \in O^+(\LAM) \mid g(F) = F\}$. Let $\overline{\Omega}_{d^{\perp}}^{+}$ 
be the closure of ${\Omega}_{d^{\perp}}^{+}$ in $\{\omega \in \PP(d^{\perp} \otimes \CC) \mid (\omega, \omega) = 0\}$. 
Then the description of the boundary component of the Baily--Borel compactification (see \cite[Sect.~4.1]{Sterk}) implies that the closure of the boundary component 
corresponding to $F$ is ${\bf P}(F\otimes{\bf C})\cap\overline{\Omega}_{d^{\perp}}^{+}$, on which 
$N(F)$ acts naturally. 

We will show in the appendix (see Lemma~\ref{lemma:appendix:extension}) that 
any element of ${\rm SL}(F)$ lifts to an element of $O^{+}(d^{\perp})$. 
Together with the surjectivity of the restriction map $O^{+}(\LAM)_{d}\to O^{+}(d^{\perp})$ (see Step 1), 
we obtain that the restriction map $N(F) \ni g \mapsto \rest{g}{F} \in \SL(F)$ is surjective. 
It follows that $N(F)\backslash{\bf P}(F\otimes{\bf C})\cap\overline{\Omega}_{d^{\perp}}^{+} 
= \SL(F)\backslash{\bf P}(F\otimes{\bf C})\cap\overline{\Omega}_{d^{\perp}}^{+}$, 
so ${\mathcal M}_{d^{\perp}}^{*}\setminus{\mathcal M}_{d^{\perp}} = X(1)$. 

>From the construction of the the Baily--Borel  compactification (see \cite[Sect.~4.1]{Sterk}) and 
Lemma~\ref{lemma:isotoropic:sublattices}~(1), 
 the map \eqref{eqn:inclusions:boundary:components} from ${\mathcal M}_{d^{\perp}}^{*}\setminus{\mathcal M}_{d^{\perp}}$to ${\mathcal M}^{*}\setminus{\mathcal M}$ is given 
 by the identity map on $X(1)$.
Since the map \eqref{eqn:inclusions:boundary:components} factors through ${\mathcal D}^{*}\setminus{\mathcal D}$, 
we get ${\mathcal D}^{*}\setminus{\mathcal D}=X(1)$.
\end{pf}

\subsection{The Borcherds $\Phi$-function}
\label{subsec:Borcherds:Phi}
In \cite{Borcherds96}, Borcherds constructed an automorphic form on $\Omega_{\LAM}^{+}$ for $O^{+}(\LAM)$ of weight $4$ vanishing exactly on ${\mathcal H}$. 
We call this automorphic form the {\em Borcherds $\Phi$-function}. 
Let us recall its definition.
For a subset $S\subset{\bf P}(\LAM\otimes{\bf C})$, let $C(S):=\{\eta\in(\LAM\otimes{\bf C})\setminus\{0\};\,[\eta]\in S\}$ be the cone over $S$.
Up to a constant, the Borcherds $\Phi$-function is defined as the holomorphic function $\Phi$ on $C(\Omega_{\LAM}^{+})$ 
with the following properties:
\begin{itemize}
\item
$\Phi(\lambda Z)=\lambda^{-4}\,\Phi(Z)$ for all $\lambda\in{\bf C}^{\times}$.
\item
$\Phi(g(Z))=\chi(g)\,\Phi(Z)$ for all $g\in O^{+}(\LAM)$, where $\chi\in{\rm Hom}(O^{+}(\LAM), \{\pm1\})$.
\item
The zero divisor of $\Phi$ is the cone $C({\mathcal H})$.
\end{itemize}
By choosing a section from $\Omega_{\LAM}^{+}$ to $C(\Omega_{\LAM}^{+})$ and pulling back $\Phi$ by the section, 
$\Phi$ is identified with a holomorphic function on $\Omega_{\LAM}^{+}$. 
In the following, we give 
two such identifications of 
the Borcherds $\Phi$-function, corresponding to the choices of $0$-dimensional cusps of level one and two. 
Then $\Phi$ will be defined without an ambiguity of constant by giving explicit infinite product expressions at those cusps.

\subsubsection{A tube domain realization of $\Omega_{\LAM}^{+}$ with respect to 
the $0$-dimensional cusp of level $\ell$}
\label{subsubsec:tube:domain}

By \cite[Prop.~4.5]{Sterk} and \cite[Sect.~2]{Scattone87} (see also Sect.~\ref{sect:moduli:Enriques:Baily:Borel}), 
the $O({\LAM})$-orbits of primitive isotropic vectors of ${\LAM}$, or equivalently the $0$-dimensional cusps of ${\mathcal M}^{*}$ consist of two points: 
the level $1$ and level $2$ cusps.

Let ${\bf e}_1, {\bf f}_1$ (resp. ${\bf e}_2, {\bf f}_2$) be the standard basis of $\UU$ (resp. $\UU(2)$) of the middle (resp. the left) sublattice of 
$\LAM = \UU(2) \oplus \UU \oplus \EE_8(2)$. Note that ${\bf e}_1$ (resp. ${\bf e}_2$) is a primitive isotropic vector of ${\LAM}$ of level $1$ (resp. level $2$). 
Let $\ell \in \{1, 2\}$. We set 
\begin{equation}
\label{eqn:def:M:ell}
{\mathbb M}_{\ell} \colonequals \UU(2/\ell)\oplus\EE_8(2). 
\end{equation}
Then ${\mathbb M}_{\ell}$ is a Lorentzian lattice of rank $10$ 
and $\LAM = (\ZZ {\bf e}_\ell + \ZZ {\bf f}_\ell) \oplus {\mathbb M}_{\ell}$ in the obvious way. 
Let $\Ccal_{{\mathbb M}_{\ell}} \colonequals \{
x \in {\mathbb M}_{\ell}\otimes \RR \mid (x, x) > 0\}$ be the positive cone of the Lorentzian lattice ${\mathbb M}_{\ell}$. 
Let $\Ccal_{{\mathbb M}_{\ell}}^{+}$ be one of the two connected components of 
$\Ccal_{{\mathbb  M}_{\ell}}$ such that 
${\mathbb M}_{\ell}\otimes\RR + i\,\Ccal_{{\mathbb M}_{\ell}}^{+}$ corresponds to $\Omega_{\LAM}^{+}$ 
via \eqref{eqn:tube:domain:level2}. 
Let $ \overline{\mathcal C}_{{\mathbb M}_{2/\ell}}^{+}$ be the closure of 
${\mathcal C}_{{\mathbb M}_{2/\ell}}^{+}$ in ${\mathbb M}_{2/\ell} \otimes \RR$. 
In what follows, we assume that the basis $\{ {\mathbf e}_{\ell}, {\mathbf f}_{\ell} \}$ is chosen in such a way that 
\begin{equation}
\label{eqn:cond:basis}
{\mathbf e}_{\ell}, {\mathbf f}_{\ell} \in \overline{\mathcal C}_{{\mathbb M}_{2/\ell}}^{+}.
\end{equation} 
Then the tube domain ${\mathbb M}_{\ell}\otimes\RR + i\,\Ccal_{{\mathbb M}_{\ell}}$ is identified with $\O_\LAM$ via the map
\begin{equation}
  \label{eqn:tube:domain:level2}
  j_{\ell} \colon {\mathbb M}_{\ell}\otimes\RR +i\,\Ccal_{{\mathbb M}_{\ell}}  \ni u \mapsto 
  \left[ 
  -(u^{2}/2){\ebf}_{\ell} + ({\fbf}_{\ell}/\ell) + (-1)^{2/\ell}u 
  \right] 
  \in \O_\LAM. 
\end{equation}
Here the sign $(-1)^{2/\ell}$ is due to the condition \eqref{eqn:cond:basis}. 
By an abuse of notation, we also write $j_{\ell}(u) =  -(u^{2}/2){\ebf}_{\ell} + {\fbf}_{\ell}/\ell + (-1)^{2/\ell} u \in C(\Omega_{\LAM})$. 
Thus we obtain an isomorphism $j_{\ell}\colon {\mathbb M}_{\ell}\otimes\RR +i\,\Ccal_{{\mathbb M}_{\ell}}^+ \to 
\O_\LAM^+$. 
On $ {\mathbb M}_{\ell}\otimes\RR +i\,\Ccal_{{\mathbb M}_{\ell}}^+$ acts $O^+(\LAM)$ via the identification~\eqref{eqn:tube:domain:level2}. 
Namely, for $g \in O^+(\LAM)$ and $u \in  {\mathbb M}_{\ell}\otimes\RR +i\,\Ccal_{{\mathbb M}_{\ell}}^+$, we define 
\[
  g \cdot u \colonequals j^{-1}_\ell g(j_\ell(u)). 
\]

\subsubsection{Series $\{c(n)\}_{n \geq -1}$}
\label{subsubsec:c(n)}
Recall that the Dedekind $\eta$-function is defined as  
\begin{equation}
\label{eqn:Dedekind:eta}
\eta(\t)  \colonequals  e^{2\pi i\tau/24}  \prod_{n>0} \left(1-e^{2\pi in\tau}\right).
\end{equation} 
To describe the infinite product expansions of $\Phi$, let 
$\{c(n)\}_{n \geq -1} \subset {\mathbf Z}$ be  the series defined as  the generating function (see \cite[p.~705]{Borcherds96}) 
\[
\eta(\t)^{-8} \eta(2\t)^8 \eta(4\t)^{-8} = \sum_{n \geq -1} c(n)\,e^{2\pi in\tau} = e^{-2\pi i\tau} + 8 + 36e^{2\pi i\tau} +\cdots. 
\]

\subsubsection{The Borcherds $\Phi$-function with respect to the $0$-dimensional cusp of level~$1$}
\label{subsubsec:Phi:level:1}

In \cite[Ex.~3.1]{Borcherds98} (see also \cite[Sect.~2.2.4]{KawaguchiMukaiYoshikawa18}), 
Borcherds introduced  the following infinite product for $z\in{\mathbb M}_{1}\otimes{\bf R}+i\,{\mathcal C}_{{\mathbb M}_{1}}^{+}$ with $(\Ima\, z)^2 \gg 0$: 
\begin{equation}
\label{eqn:Phi:expansion:level1}
\Phi_{1}(z)
 \colonequals 
\prod_{\lambda\in{\mathbb M}_{1}\cap\overline{\mathcal C}_{{\mathbb M}_{1}}^{+}\setminus\{0\}}
\left(
\frac{1-e^{\pi i\langle\lambda,z\rangle_{{\mathbb M}_{1}}}}
{1+e^{\pi i\langle\lambda,z\rangle_{{\mathbb M}_{1}}}}
\right)^{c(\lambda^{2}/2)}.
\end{equation}
Then $\Phi_{1}(z)$ converges absolutely when $(\Ima\, z)^2 \gg 0$ and extends holomorphically to the whole 
${\mathbb M}_{1}\otimes\RR + i\,\Ccal_{{\mathbb M}_{1}}^+$.
By \cite[Ex.~3.1]{Borcherds98}, $\Phi_1 \circ j_{1}^{-1}$ is 
an automorphic form on $\O_\LAM^+$ for $O^+(\LAM)$ of weight $4$ with zero divisor $\Hcal$. 
We call $\Phi_{1}(z)$ 
the {\em Borcherds $\Phi$-function with respect to the $0$-dimensional cusp of level~$1$}. 
Then $\Phi$ is normalized in such a way that $\Phi_{1} = j_{1}^{*}\Phi$. 

\subsubsection{The Borcherds $\Phi$-function with respect to the $0$-dimensional cusp of level~$2$}
\label{subsubsec:Phi:level:2}

Recall that 
$\{ {\mathbf e}_1, {\mathbf f}_1 \}$ is the standard free-basis of $\UU$. 
We set 
\begin{equation}
\label{eqn:positive:roots}
\Pi^+  \colonequals  \{\l \in{\mathbb M}_{2} \mid \langle \l, {\bf e}_1 \rangle_{{\mathbb M}_{2}}> 0,\,\lambda^{2}\geq-2 \}.
\end{equation}
(Our $\Pi^{+}$ is slightly different form the one in \cite[p.~701]{Borcherds96} for the convenience in later use.)
In \cite[Sect.~3]{Borcherds96}, Borcherds introduced the following infinite product for $w\in{\mathbb M}_{2}\otimes{\bf R}+i\,{\mathcal C}_{{\mathbb M}_{2}}^{+}$ with $(\Ima\, w)^2 \gg 0$:
\begin{equation}
\label{eqn:Phi:expansion:level2}
  \Phi_{2}(w) 
   \colonequals  
  2^{8}e^{2\pi i\langle {\bf e}_1, w\rangle_{{\mathbb M}_{2}}}
  \prod_{\lambda \in \ZZ_{>0}{\bf e}\,\cup\,\Pi^+} 
  \left(
  1 - e^{2\pi i \langle \l, w\rangle_{{\mathbb M}_{2}}}
  \right)^{(-1)^{\langle \lambda, 
  {\bf e}_1 - {\bf f}_1\rangle} c(\lambda^2/2)}.
\end{equation}
Here the constant $2^{8}$ comes from the one in \cite[Th.~13.3 (5)]{Borcherds98} (see also \cite[Th.~2.2]{KawaguchiMukaiYoshikawa18}). 
Then $\Phi_{2}(w)$ converges absolutely when $(\Ima\, w)^2 \gg 0$ and extends 
holomorphically to the whole 
${\mathbb M}_{2}\otimes\RR + i\,\Ccal_{{\mathbb M}_{2}}^+$.
By \cite[Sect.~3]{Borcherds96}, $\Phi_2 \circ j_{2}^{-1}$ is 
an automorphic form on $\O_\LAM^+$ for $O^+(\LAM)$ of weight $4$ with zero divisor $\Hcal$. 
We call $\Phi_{2}(w)$ 
the {\em Borcherds $\Phi$-function with respect to the $0$-dimensional cusp of level~$2$}. 
Then we have an equality $\Phi_{2} = j_{2}^{*}\Phi$ of functions on $\Omega_{\LAM}^{+}$. 
Indeed, there is a constant $C'$ with $|C'|=1$ such that $\Phi_{2}=C'j_{2}^{*}\Phi$ by \cite[Th.\,13.3 (5)]{Borcherds98}.
Then $C'=1$ by comparing \cite[(2.12)]{KawaguchiMukaiYoshikawa18} with $C=2^{8}$ and the formula \eqref{eqn:Phi:expansion:level2}.

\subsubsection{Remarks}
\label{subsubsec:Remarks}
Since $\Phi_{\ell}$ is a function on ${\mathbb M}_{2/\ell}\otimes{\mathbf R} + i {\mathcal C}_{{\mathbb M}_{2/\ell}}^{+}$,  
$\Phi_{1} = j_{1}^{*}\Phi$ and $\Phi_{2} =  j_{2}^{*}\Phi$ are two distinct realizations of $\Phi \in {\mathcal O}(C_{\Omega_{\LAM}^{+}})$ 
as a function on the tube domain of a $10$-dimensional affine space. 
For the precise relation between them, see \cite[Sect.~2.2.3]{KawaguchiMukaiYoshikawa18}.
In \cite[Ex.~13.7]{Borcherds98}, Borcherds used the lattice ${\LAM}^{\lor}(2)$ instead of ${\LAM}$. 
For this reason, the infinite product expansion with respect to the level $\ell$ cusp in \cite[Ex.~13.7]{Borcherds98} 
corresponds to that for $\Phi_{2/\ell}$ with respect to the 
level $2/\ell$ cusp.

\subsubsection{Automorphy of the Borcherds $\Phi$-function}
\label{subsubsec:automorphiy:Borchreds}
Let $\ell \in \{1, 2\}$. 
By the automorphy of the Borcherds $\Phi$-function, for any $g \in O^{+}({\LAM})$ and $u \in {\mathbb M}_{\ell}\otimes\RR + i\,\Ccal_{{\mathbb M}_{\ell}}^+$, we have 
\[
\Phi_{\ell}(g \cdot u) = \chi(g) \, J_{\ell}(g, u)^4 \, \Phi_{\ell}(u), 
\]
where $J_{\ell}(g,u) \in {\mathcal O}({\mathbb M}_{\ell}\otimes{\mathbf R} + i {\mathcal C}_{{\mathbb M}_{\ell}}^{+})$ is the automorphic factor defined as
$$
J_{\ell}(g, u) := \langle g( -(u^{2}/2) {\bf e}_{\ell} + ({\bf f}_{\ell}/\ell) + (-1)^{2/\ell} u ), {\bf e}_{\ell} \rangle_{\LAM}.
$$
Since $\chi^2 = 1$ by \cite[proof of Prop.~5.6]{GHS} (see also \cite[Lem.~2.1]{KawaguchiMukaiYoshikawa18}),  
we have 
\begin{equation}
\label{eqn:automorphy:Phi}
\Phi_{\ell}(g \cdot u)^2 = J_{\ell}(g, u)^8 \, \Phi_{\ell}(u)^2. 
\end{equation}

We will use the following invariance of the square of the Borcherds $\Phi$-function. 

\begin{lemma}
\label{lemma:invariance:Phi:stab}
Let $\ell \in \{1, 2\}$. 
Let $g \in O^{+}(\LAM)$ be such that $g({\bf e}_{\ell}) = {\bf e}_{\ell}$. Then $\Phi_{\ell}(g \cdot u)^2 = \Phi_{\ell}(u)^2$ for all $u \in {\mathbb M}_{\ell}\otimes\RR + i\,\Ccal_{{\mathbb M}_{\ell}}^+$. 
\end{lemma}

\begin{pf}
Since $g^{-1}({\bf e}_{\ell}) = {\bf e}_{\ell}$, we have 
\begin{align*}
J_{\ell}(g, u)
& = \langle g( -(u^{2}/2) {\bf e}_{\ell} + ({\bf f}_{\ell}/\ell) + (-1)^{2/\ell} u ), {\bf e}_{\ell} \rangle_{\LAM}
\\
&= 
 \langle  -(u^{2}/2) {\bf e}_{\ell} + ({\bf f}_{\ell}/\ell) + (-1)^{2/\ell} u, g^{-1}({\bf e}_{\ell}) \rangle_{\LAM}\\
& = \langle  -(u^{2}/2) {\bf e}_{\ell} + ({\bf f}_{\ell}/\ell) + (-1)^{2/\ell} u, {\bf e}_{\ell} \rangle_{\LAM}
= \langle {\bf f}_{\ell}/\ell,  {\bf e}_{\ell} \rangle_{\LAM}
= 1. 
\end{align*}
Now the result follows from \eqref{eqn:automorphy:Phi}. 
\end{pf}

\subsubsection{The Petersson norm}
\label{subsubsec:Petersson}
Let $\ell \in \{1, 2\}$. 
The Petersson norm of $\Phi_\ell$ is the $C^{\infty}$ function on 
${\mathbb M}_{\ell}\otimes\RR + i\,\Ccal_{{\mathbb M}_{\ell}}^{+}$ defined as 
\begin{equation}
\label{eqn:def:Petersson}
\|\Phi_\ell(u)\|^{2} \colonequals \langle{\Ima\, u},{\Ima\, u}\rangle_{{\mathbb M}_{\ell}}^{4}|\Phi_\ell(u)|^{2}.
\end{equation}

Since $\langle{\Ima\, (g\cdot u)},{\Ima\, (g\cdot u)}\rangle_{{\mathbb M}_{\ell}} = |J_{\ell}(g, u)|^{-2} \langle{\Ima\, u},{\Ima\, u}\rangle_{{\mathbb M}_{\ell}}$,
it follows from \eqref{eqn:automorphy:Phi} that $\|\Phi_\ell\|^{2}$ is $O^{+}({\LAM})$-invariant, and 
$\|\Phi_\ell\|^{2}$ is regarded as a $C^{\infty}$ function on ${\mathcal M}$.
Moreover, if $j_{1}(u_{1}) = j_{2}(u_{2})$, then $\|\Phi_{1}(u_{1})\| = \|\Phi_{2}(u_{2})\|$. Hence it makes sense to define $\|\Phi(Z)\|$ as
$\|\Phi_{1}(u_{1})\|=\|\Phi_{2}(u_{2})\|$ when $Z=j_{1}(u_{1})=j_{2}(u_{2})$.  
Recall that $\Phi$ is defined as a holomorphic function on $C(S)$.
As a function on $C(S)$, we have $\|\Phi(Z)\|^{2} = 2^{-4} \langle Z, \overline{Z}\rangle_{\LAM}^{4} |\Phi(Z)|^{2}$. 
For an Enriques surface $Y$, we define
\begin{equation}
\label{eqn:Petersson:norm}
\|\Phi(Y)\| := \| \Phi(\overline{\varpi}(Y) \|.
\end{equation}

\section
{Fixed-point-free involutions on Kummer surfaces of product type}
\label{sec:involution}
In this section, we recall Kummer surfaces of product type and fixed-point-free involutions 
on them and their parities. Then we study some properties of the period maps according to their parities. 

\subsection{Kummer surfaces of product type and their periods}
\label{sect:Kummer:prod}

For $\tau\in \Hfrak$, we set 
$$
E_{\tau} \colonequals {\bf C}/{\bf Z}+\tau{\bf Z}.
$$
For $(\tau,\tau')\in\HH\times\HH$, let $f\colon\widetilde{E_{\tau}\times E_{\tau'}} \to E_{\tau}\times E_{\tau'}$ be the blowing-up 
of the points of order $2$ of $E_{\tau}\times E_{\tau'}$.
Then the involution $[-1](x) \colonequals -x$ on $E_{\tau}\times E_{\tau'}$ induces an involution on $\widetilde{E_{\tau}\times E_{\tau'}}$, which is denoted by $\widetilde{[-1]}$. 
The quotient
$$
K_{\tau,\tau'}=\Km(E_{\tau}\times E_{\tau'}) \colonequals \widetilde{E_{\tau}\times E_{\tau'}} / \widetilde{[-1]}
$$
is a $K3$ surface called a {\em Kummer surface of product type}.

We consider the map
\begin{equation}
\label{eqn:phi}
   \phi\colon 
   H^1(E_\t, \ZZ) \otimes H^1(E_{\t'}, \ZZ) 
  {\hookrightarrow} H^2(E_\t\times E_{\t'}, \ZZ) 
  \overset{\widetilde{p}_{!}\circ f^{*}}{\longrightarrow} 
   H^2(K_{\tau,\tau'}, \ZZ), 
\end{equation}
where the first map is given by the cup-product and $\widetilde{p}\colon\widetilde{E_{\tau}\times E_{\tau'}}\to K_{\tau,\tau'}$ in the second map is the projection. We set 
\begin{equation}
\label{eqn:def:K}
  {\mathbf K} \colonequals   \phi\left(
H^1(E_\t, \ZZ) \otimes H^1(E_{\t'}, \ZZ)
\right) \, \subset \, H^2(K_{\tau,\tau'}, \ZZ). 
\end{equation}
By \cite[VIII, Proposition~(5.1)]{BPV84}, we have an isometry of lattices ${\mathbf K} \cong \UU(2)\oplus\UU(2)$. 
Since $H^{0}(K_{\tau, \tau'}, \Omega_{K_{\tau,\tau'}}^{2}) \subset {\mathbf K}\otimes{\mathbf C}$, we have
$T_{K_{\tau,\tau'}} \subset {\mathbf K}$. If $(\tau,\tau')\in{\mathfrak H}\times{\mathfrak H}$ is very general, then
$T_{K_{\tau,\tau'}} = {\mathbf K}$.

Let ${\bf a}$ (resp. ${\bf b}$) be the element of $H_{1}(E_{\tau},{\bf Z})$ corresponding to the line segment
$[0,\tau]$ (resp. $[0,1]$) of ${\bf C}$. 
Let ${\bf a}^{\lor}$ (resp. ${\bf b}^{\lor}$) be the Poincar\'e dual of ${\bf a}$ (resp. ${\bf b}$) such that
$\int_{\bf a}\eta=\int_{E_{\tau}}{\bf a}^{\lor}\wedge\eta$ and $\int_{\bf b}\eta=\int_{E_{\tau}}{\bf b}^{\lor}\wedge\eta$
for any closed $1$-form $\eta$ on $E_{\tau}$.
Then $\{{\bf a}^{\lor},{\bf b}^{\lor}\}$ is a symplectic basis of $H^{1}(E_{\tau},{\bf Z})$ 
such that $[dz]={\bf a}^{\lor}-\tau{\bf b}^{\lor}$ in $H^{1}(E_{\tau},{\bf Z})$. 
Similarly, we define a symplectic basis $\{{\bf a}^{\prime\lor},{\bf b}^{\prime\lor}\}$ of $H^{1}(E_{\tau'},{\bf Z})$.  

Following \cite[Sects.\,6.1. 6.2, 7.3]{KawaguchiMukaiYoshikawa18}, we set
\begin{equation}
\label{eqn:identification:lattices}
\begin{array}{ll}
\Gamma_{12}^{\lor} \colonequals \phi({\bf a}^{\lor}\otimes{\bf a}^{\prime\lor}),
&
\Gamma_{34}^{\lor} \colonequals \phi({\bf b}^{\lor}\otimes{\bf b}^{\prime\lor}),
\\
\Gamma_{14}^{\lor} \colonequals \phi({\bf a}^{\lor}\otimes{\bf b}^{\prime\lor}),
&
\Gamma_{23}^{\lor} \colonequals - \phi({\bf b}^{\lor}\otimes{\bf a}^{\prime\lor}).
\end{array}
\end{equation}
Then  $\{\Gamma_{34}^{\lor}, \Gamma_{12}^{\lor}, \Gamma_{14}^{\lor}, \Gamma_{23}^{\lor}\}$ is a basis of ${\mathbf K}$
with Gram matrix $-\binom{0\,2}{2\,0} \oplus -\binom{0\,2}{2\,0}$. 

\medskip
Suppose that we are given a fixed-point-free involution 
$$
\iota\colon K_{\tau,\tau'} \to K_{\tau,\tau'}.
$$ 
Then it is anti-symplectic. 
Let $H^{2}(K_{\tau,\tau'}, {\mathbf Z})_{\pm}$ be the $\pm1$-eigenspace of $H^{2}(K_{\tau,\tau'}, {\mathbf Z})$ with respect to the $\iota$-action.
We suppose moreover that $\iota$ satisfies the condition
\begin{equation}
\label{eqn:assumption:lattice}
{\mathbf K} \subset H^{2}(K_{\tau,\tau'}, {\mathbf Z})_{-}.
\end{equation}
Since $H^{0}(K_{\tau,\tau'}, \Omega_{K_{\tau,\tau'}}^{2}) \subset {\mathbf K}\otimes{\mathbf C}$, there exist no roots 
$d \in \Delta_{H^{2}(K_{\tau,\tau'}, {\mathbf Z})_{-}}$ with ${\mathbf K} \subset d^{\perp}$ by \eqref{eqn:assumption:lattice}.
Since $\iota$ is anti-symplectic and thus $T_{K_{\tau,\tau'}} \subset H^{2}(K_{\tau,\tau'}, {\mathbf Z})_{-}$,
\eqref{eqn:assumption:lattice} holds if $T_{K_{\tau,\tau'}} = {\mathbf K}$. 
The geometric meaning of \eqref{eqn:assumption:lattice} is given as follows.

\begin{lemma}
\label{lemma:extendable:involution}
Let $U$ be a neighborhood of $(\tau,\tau')$ in ${\mathfrak H}\times{\mathfrak H}'$. Let $\pi \colon {\mathcal K} \to U$ be a family of Kummer surfaces
such that $\pi^{-1}(\sigma, \sigma') \cong K_{\sigma, \sigma'}$ for all $(\sigma, \sigma') \in U$. Let $\theta$ be a fixed-point-free involution on $K_{\tau,\tau'}$.
If $U$ is sufficiently small and contractible, then the following conditions are equivalent.
\begin{itemize}
\item[(1)]
$\theta$ satisfies \eqref{eqn:assumption:lattice}.
\item[(2)]
There exists an involution $\theta_{\mathcal K} \colon {\mathcal K} \to {\mathcal K}$ preserving the fibers of $\pi$ 
such that $\theta_{\mathcal K}|_{K_{\tau, \tau'}} = \theta$.
\end{itemize}
\end{lemma}

\begin{pf}
The result follows from the global Torelli theorem. The details are left to the reader. 
(The fixed-point-free involutions that we will need satisfy the condition~\eqref{eqn:assumption:lattice} (see Theorem~\ref{thm:Mukai:Ohashi}, Sect.~\ref{sec:L} and Sect.~\ref{subsec:realization:odd}).)
\end{pf}

\par
Recall that we have fixed a primitive embedding $\LAM \subset \LL$ throughout this paper and that  
$\{ {\bf e}_1, {\bf f}_1 \}$ (resp. $\{ {\bf e}_2, {\bf f}_2 \}$) is the standard basis of $\UU$ (resp. $\UU(2)$) of the 
middle (resp. the left) sublattice of $\LAM \colonequals \UU(2) \oplus \UU \oplus \EE_8(2)$. 

\begin{definition}
\label{def:normalized:marking}
For a fixed-point-free involution $\iota$ on $K_{\tau,\tau'}$ satisfying \eqref{eqn:assumption:lattice},
let $\ell \in \{1, 2\}$ denote the level of $\Gamma_{34}^{\lor}$ in $H^2(K_{\tau,\tau'}, \ZZ)_{-}$. 
An isometry $\alpha: H^2(K_{\tau,\tau'}, \ZZ) \to \LL$ is called a {\em normalized marking} for 
$(K_{\tau,\tau'}, \iota)$ if it satisfies 
\eqref{eqn:alpha:+:-}, \eqref{eqn:definition:period} and 
\begin{equation}
\label{eqn:marking:normalization}
-\alpha( \Gamma_{34}^{\lor} ) = {\bf e}_{\ell}.
\end{equation}
\end{definition}

\begin{lemma}
\label{lemma:normalized:marking}
If \eqref{eqn:assumption:lattice} holds, then there exists
a normalized marking $\alpha$ for $(K_{\tau, \tau'}, \iota)$.
\end{lemma}

\begin{pf}
Let $\alpha' \colon H^{2}(K_{\tau, \tau'}, {\mathbf Z}) \to \LL$ be a marking satisfying \eqref{eqn:alpha:+:-}. 
By \eqref{eqn:assumption:lattice}, $\Gamma_{34}^{\lor} \in {\mathbf K} \subset H^{2}(K_{\tau, \tau'}, {\mathbf Z})_{-}$.
Then $\ell$ is the level of $\alpha'(\Gamma_{34}^{\lor}) \in \LAM$. Set $I := \alpha' \iota^{*} (\alpha')^{-1}$. 
Then $\LAM$ is exactly the anti-invariant subspace of $\LL$ with respect to the $I$-action. 
Since the $O(\LAM)$-orbit of a primitive isotropic vector of $\LAM$ is determined by its level, there exists $g \in O(\LAM)$ 
such that $g( \alpha'(\Gamma_{34}^{\lor}) ) = -{\bf e}_{\ell}$.

Replacing $g$ by $s \circ g$ where $s \in O({\mathbb M}_{\ell})$ exchanges the components ${\mathcal C}_{{\mathbb M}_{\ell}}^{\pm}$ if necessary, 
we may assume $g \in O^{+}(\LAM)$.  
By \cite[Remark\,1.15]{Namikawa85}, there exists $\widetilde{g} \in O( \LL )$ with $\widetilde{g}\, I = I \,\widetilde{g}$ such that $\widetilde{g}|_{\LAM} = g$. 
Then $\alpha := \widetilde{g} \circ \alpha'$ is a marking on $K_{\tau, \tau'}$ satisfying \eqref{eqn:alpha:+:-}, \eqref{eqn:definition:period}, \eqref{eqn:marking:normalization}.
\end{pf}

Let $\alpha\colon H^2(K_{\tau,\tau'}, \ZZ) \to \LL$ be a normalized marking for $(K_{\tau,\tau'}, \iota)$. 
Then the period of a (normalized) marked Enriques surface $(K_{\tau,\tau'}/ \iota, \alpha)$ is given by that of $(K_{\tau,\tau'}, \alpha)$, i.e.,
\begin{equation}
\label{eqn:period}
  \varpi(K_{\tau,\tau'}/\iota, \alpha)
  = 
  \left[\alpha\left(\phi\left(H^0(E_\t, \O^1_{E_\t}) \otimes 
  H^0(E_{\t^\prime}, \O^1_{E_{\t^\prime}})\right)\right)\right] \in \O_{\LAM}^+.  
\end{equation}
By \eqref{eqn:period}, \eqref{eqn:identification:lattices} and the relation $[dz] = {\mathbf a}^{\lor} - \tau {\mathbf b}^{\lor}$, the period of $(K_{\tau, \tau'}, \alpha)$ 
is concretely expressed as follows 
(see  \cite[Sects.~7.3, 7.4]{KawaguchiMukaiYoshikawa18}):
\begin{equation}
\label{eqn:(8.3)}
\varpi(K_{\tau,\tau'}/\iota, \alpha)
 =
\left[
\tau\tau' \alpha(\Gamma_{34}^{\lor}) + \alpha(\Gamma_{12}^{\lor}) + \tau \alpha(\Gamma_{23}^{\lor}) - \tau'\alpha(\Gamma_{14}^{\lor})
\right]. 
\end{equation}

We make a crucial observation on the value of Borcherds $\Phi$-function. 

\begin{lemma}
\label{lemma:normalized:marking:2}
Let $\iota$ be a fixed-point-free involution on $K_{\tau, \tau'}$ satisfying \eqref{eqn:assumption:lattice}, and 
let $\ell$ be the level of $\Gamma_{34}^{\lor}$ in $H^{2}(K_{\tau, \tau'}, {\mathbf Z})_{-}$.
Let $\alpha, \alpha^\prime$ be normalized markings for $(K_{\tau, \tau'}, \iota)$. 
We put $u \colonequals j_\ell^{-1}\left(\varpi(K_{\tau,\tau'}/\iota, \alpha)
\right)$ and 
$u^\prime \colonequals j_\ell^{-1}\left(\varpi(K_{\tau,\tau'}/\iota, \alpha^\prime)\right)$.  
Then we have 
$
  \Phi_\ell\left(u\right)^2 = \Phi_\ell\left(u^\prime\right)^2. 
$
In other words, the value $\Phi_\ell(u)^2$ is independent of the choice 
of a normalized marking. 
\end{lemma}

\begin{pf}
We put $g \colonequals \alpha^\prime \circ \alpha^{-1}$. Then 
$g \in O(\LL)$ satisfies 
$g(\LAM) =\LAM$, $g(\Omega_{\LAM}^{+}) = \Omega_{\LAM}^+$ and $g({\bf e}_{\ell}) = {\bf e}_{\ell}$. 
It follows that $\rest{g}{\LAM} \in O^+(\LAM)$ and $\rest{g}{\LAM}({\bf e}_{\ell}) = {\bf e}_{\ell}$. 
Since $j_{\ell}(u')= \varpi(K_{\tau, \tau'}/\iota, \alpha') = g(\varpi(K_{\tau, \tau'}/\iota, \alpha)) = g(j_{\ell}(u))$, 
the definition of the action of $O^+(\LAM)$ on ${\mathbb M}_{\ell}\otimes{\mathbf R} + i\,{\mathcal C}_{{\mathbb M}_{\ell}}^{+}$ gives $u' = \rest{g}{\LAM} \cdot u$. 
Hence the result follows from 
Lemma~\ref{lemma:invariance:Phi:stab}.
\end{pf}

By \cite[Sects.\,6.4 and 7.4]{KawaguchiMukaiYoshikawa18}, we can express
\begin{equation}
\label{eqn:period:prod:kummer}
\varpi(K_{\tau, \tau'}, \alpha) = j_{\ell}(u) =
\left[
-(u^{2}/2) {\ebf}_{\ell} + ( {\fbf}_{\ell}/\ell ) + (-1)^{2/\ell} u
\right],
\end{equation}
where $u = u(\tau, \tau'; \alpha) \in {\mathbb M}_{\ell}\otimes{\mathbf R} + i \, {\mathcal C}_{{\mathbb M}_{\ell}}^{+}$ is given explicitly by
\begin{equation}
\label{eqn:formula:period}
(-1)^{2/\ell} u = \frac{1}{2}( A + B \tau + D \tau' ),
\qquad
\begin{cases}
\begin{array}{ll}
A &:= \alpha(\Gamma_{12}^{\lor}) - \langle {\fbf_\ell}/{\ell}
, \alpha(\Gamma_{12}) \rangle_{\LAM} {\ebf}_{\ell} - (2/\ell) {\fbf}_{\ell},
\\
B &:= -\langle {\fbf_\ell}/{\ell}
, \alpha(\Gamma_{23}^{\lor}) \rangle_{\LAM} {\ebf}_{\ell} + \alpha(\Gamma_{23}^{\lor}),
\\
D &:= \langle {\fbf_\ell}/{\ell}
, \alpha(\Gamma_{14}^{\lor}) \rangle_{\LAM} {\ebf}_{\ell} - \alpha(\Gamma_{14}^{\lor}).
\end{array}
\end{cases}
\end{equation}

\begin{remark}
\label{rem:misprint}
There are several misprints in \cite{KawaguchiMukaiYoshikawa18}. In \cite{KawaguchiMukaiYoshikawa18}, (6.7), Lemma~6.2
and its proof, and in Theorem~6.3 and its proof, $z_{\langle J\rangle}(T)$ should be replaced by $(-1)^{2/\ell}z_{\langle J\rangle}(T)$, 
so that $\Im z_{\langle J\rangle}(T)\in {\mathcal C}_{{\mathbb M}_{\ell}}^{+}$.
In the last line of \cite[p.1509]{KawaguchiMukaiYoshikawa18}, 
the formula for $A$ should be replaced by $A = {\fbf}' - (2/\ell){\fbf}_{\ell} - \langle {\fbf}_{\ell}/\ell, {\fbf}' \rangle {\ebf}_{\ell}$. 
For the same reason as above, in \cite{KawaguchiMukaiYoshikawa18}, p.1514, $z_{\langle J\rangle}(\tau_{1},\tau_{2})$ should be 
replaced by $(-1)^{2/\ell}z_{\langle J\rangle}(\tau_{1},\tau_{2})$. In \cite{KawaguchiMukaiYoshikawa18}, the proof of Lemma~7.1,
$B$ and $D$ should be replaced by $(-1)^{2/\ell}B$ and $(-1)^{2/\ell}D$, respectively.
\end{remark}

\subsection{The parity of an involution}
\label{sect:parity:invol}
For $\tau, \tau' \in \HH$, let $K_{\tau,\tau'}$ be a Kummer surface of product type, and let 
$\iota$ be a fixed-point-free involution on $K_{\tau,\tau'}$ satisfying \eqref{eqn:assumption:lattice}. 
We define the {\em patching element} $d_{\iota} \in A_{{\mathbf K}}\setminus\{0\}$ as follows. For simplicity, write $H^{2}_{-}$ for
$H^{2}(K_{\tau,\tau'}, {\mathbf Z})_{-}$. Let ${\mathbf K}^{\perp} = {\mathbf K}^{\perp_{H^{2}_{-}}}$ be the orthogonal complement of
${\mathbf K}$ in $H^{2}_{-}$. By \cite[Proof of Prop.~4.3]{Ohashi07}, ${\mathbf K}^{\perp} \cong {\mathbb E}_{8}(2)$ and 
we have the following inclusions
$$
{\mathbf K} \oplus {\mathbf K}^{\perp} \subset H^{2}_{-} \subset (H^{2}_{-})^{\lor} \subset
{\mathbf K}^{\lor} \oplus ({\mathbf K}^{\perp})^{\lor}.
$$
Since $\dim_{{\bf F}_2} A_{{\mathbf K}} \oplus A_{{\mathbf K}^{\perp}} = 12$ and $\dim_{{\bf F}_2} A_{H^{2}_{-}} = 10$, 
we get $\dim_{{\bf F}_2} H^{2}_{-} / ({\mathbf K} \oplus{\mathbf K}^{\perp}) = 1$. 
Hence there exists $d \in H^{2}_{-}\setminus\{0\}$ such that 
\begin{equation}
\label{eqn:the:choice:of:d}
  H^{2}_{-} = \ZZ d +{\mathbf K} \oplus{\mathbf K}^{\perp}. 
\end{equation}
We write $d = d_1 + d_2$ with $d_1 \in{\mathbf K}^\vee$ and $d_2 \in({\mathbf K}^{\perp})^\vee$. 
Since $d \not\in{\mathbf K}\oplus{\mathbf K}^{\perp}$, the primitivity of the embeddings ${\mathbf K}\subset H^{2}_{-}$ and 
${\mathbf K}^{\perp} \subset H^{2}_{-}$ implies that $d_1 \neq 0$ and $d_2 \neq 0$. 
Since $\delta({\mathbf K}) = \delta( {\mathbf K}^{\perp}) =0$,
we have $d_1^2 \in \ZZ$ and $d_2^2 \in \ZZ$. 
Since the lattice $H^{2}_{-}$ is even and hence $d^2\in2{\bf Z}$, the equality $d^2 = d_1^2 + d_2^2$ implies that $d_1^2\equiv d_2^2\mod 2$. 
Namely, $q_{\mathbf K}(\bar{d}_1) = q_{{\mathbf K}^{\perp}}(\bar{d}_2) \in \ZZ/2\ZZ$. 
Since $d\mod {\mathbf K} \oplus{\mathbf K}^{\perp}$ is determined by $\iota$, 
the value $q_{\mathbf K}(\bar{d}_1) = q_{{\mathbf K}^{\perp}}(\bar{d}_2) \in \ZZ/2\ZZ$ depends only on $\iota$.

\begin{definition}[Patching element and parity]
\label{def:KM:L}
For a fixed-point-free involution $\iota$ on $K_{\tau,\tau'}$ satisfying \eqref{eqn:assumption:lattice},
define the patching element of $\iota$ as $\overline{d}_{\iota} := \overline{d}_{1} \in A_{\mathbf K} \setminus\{0\}$.
Then $\iota$ is said to be of {\em odd} (resp. {\em even}) type 
if $q_{\mathbf K}(\dbar_{\iota}) =1$ (resp. $q_{\mathbf K}(\dbar_{\iota})=0$).
\end{definition}

\begin{remark}
\label{rem:patching:element}
The notion of patching element given in Definition~\ref{def:KM:L} coincides with that of Ohashi \cite[Def.\,4.5]{Ohashi07}.
To see it, write $H_{\pm}^{2}$ and $H^{2}$ for $H^{2}(K_{\tau,\tau'}, {\mathbf Z})_{\pm}$ and $H^{2}(K_{\tau,\tau'},{\mathbf Z})$, respectively. 
Let $S :={\mathbf K}^{\perp_{H^{2}}}$ be the orthogonal complement of ${\mathbf K}$ in $H^{2}\cong\LL$. Then
$S \cong {\mathbb U}\oplus {\mathbb E}_{8}\oplus {\mathbb D}_{4}^{\oplus2}$ and we have
${\mathbf K}^{\perp}=(H^{2}_{+})^{\perp_{S}}$, i.e., ${\mathbf K}^{\perp}$ is also given as the orthogonal complement of $H^{2}_{+}$ in $S$.
There is a subgroup $\Gamma \subset A_{H^{2}_{+}} \oplus A_{{\mathbf K}^{\perp}}$ with $\dim_{{\mathbf F}_{2}}\Gamma=7$ such that
$q_{S} = q_{H^{2}_{+}} \oplus q_{{\mathbf K}^{\perp}} |_{\Gamma^{\perp}/\Gamma}$, where $\Gamma^{\perp}$ is the orthogonal complement of $\Gamma$
with respect to the discriminant bilinear from on $A_{H^{2}_{+}}\oplus A_{{\mathbf K}^{\perp}}$ (cf. \cite[Th.\,4.2 (1)]{Ohashi07}). 
Let $\Gamma_{{\mathbf K}^{\perp}}$ be the image of $\Gamma$
by the obvious projection. Since $\Gamma \cong \Gamma_{{\mathbf K}^{\perp}}$ and $\dim_{{\mathbf F}_{2}} {\mathbf K}^{\perp} =8$, there is a unique vector 
$z \in A_{{\mathbf K}^{\perp}}$ with $\Gamma_{{\mathbf K}^{\perp}} = {\mathbf F}_{2}z$. Then $v:=[(0,z)] \in \Gamma^{\perp}/\Gamma = A_{S}$
is the patching element in \cite{Ohashi07}. Since $(0,z)$ represents $v$ and since 
$q_{H^{2}_{+}} = -q_{S}\oplus q_{{\mathbf K}^{\perp}}|_{\Gamma^{\prime \perp}/\Gamma^{\prime}}$
for some subgroup $\Gamma^{\prime} \subset A_{S}\oplus A_{{\mathbf K}^{\perp}}$ (cf. \cite[Th.\,4.2 (2)]{Ohashi07}), 
we see that the class of $(v,z) \in \Gamma^{\prime \perp} \subset A_{S}\oplus A_{{\mathbf K}^{\perp}}$ coincides with
$0$ in $A_{H^{2}_{+}}$. Since $(A_{S}, -q_{S}) = (A_{\mathbf K}, q_{\mathbf K})$ and $(A_{H^{2}_{+}}, q_{H^{2}_{+}})=(A_{H^{2}_{-}}, q_{H^{2}_{-}})$, 
this implies that $(v,z) \in H^{2}_{-}/A_{\mathbf K}\oplus A_{{\mathbf K}^{\perp}}$. Since $(v,z)\not=(0,0)$ in $A_{\mathbf K}\oplus A_{{\mathbf K}^{\perp}}$,
we get $H^{2}_{-} = {\mathbf Z}(v,z) + {\mathbf K}\oplus {\mathbf K}^{\perp}$. Hence $v \in A_{S}=A_{\mathbf K}$ is the patching element of $\iota$.
\end{remark}

Mukai \cite{Mukai10} and Ohashi \cite[Th.~0.2]{Ohashi07} classified the conjugacy classes of 
fixed-point-free involutions on Kummer surfaces of product type satisfying \eqref{eqn:assumption:lattice}. 
For $n \in \ZZ_{>0}$, let 
$\Gamma(n) :=\{ \binom{a\,b}{c\,d} \in\SL_{2}(\ZZ);\, a\equiv d\equiv1,\,b\equiv c\equiv 0\mod n\} \subset \SL_2(\ZZ)$ 
denote the principal  congruence subgroup of level $n$. 
Let $\varDelta_{\mathfrak H}$ be the diagonal locus of ${\mathfrak H}\times{\mathfrak H}$. We set
\begin{equation}
\label{eqn:frakD}
{\mathfrak D}
 \colonequals 
\bigcup_{\gamma\in \Gamma(2)}
(\gamma\times 1)\varDelta_{\mathfrak H}
=
\bigcup_{\gamma\in \Gamma(2)}
(1\times\gamma)\varDelta_{\mathfrak H}.
\end{equation}

\begin{theorem}[\cite{Mukai10}, \cite{Ohashi07}]
\label{thm:Mukai:Ohashi}
For $\tau, \tau' \in {\mathfrak H}$, let $K_{\tau,\tau'}$ be a Kummer surface of product type and let $\omega_{K_{\tau,\tau'}}$ be its non-zero canonical form.
\begin{enumerate}
\item
Let $(\tau, \tau') \in {\mathfrak H}\times{\mathfrak H}\setminus {\mathfrak D}$. 
Then there exist $15$ conjugacy classes of fixed-point-free involutions on $K_{\tau,\tau'}$ satisfying \eqref{eqn:assumption:lattice} 
and there exists a bijection between the set of these $15$ conjugacy classes of involutions and $A_{{\mathbf K}}\setminus\{0\}$
given by the assignment $\iota \mapsto \overline{d}_{\iota}$, where $\overline{d}_{\iota} \in A_{\mathbf K}\setminus\{0\}$ is the patching element of $\iota$. 
\item
If $(\tau, \tau')$ is very general, i.e., $T_{K_{\tau, \tau'}} = {\mathbf K}$ 
and ${\rm Aut}(T_{K_{\tau,\tau'}}, \omega_{K_{\tau,\tau'}}) = \{\pm1\}$, 
then the $15$ fixed-point-free involutions in~\textup{(1)} 
are, up to conjugacy,  the only fixed-point-free involutions on $K_{\tau,\tau'}$. 
\end{enumerate}
\end{theorem}

We will recall geometric descriptions of these $15$ involutions in Sect.~\ref{sec:L} and Sect.~\ref{subsec:realization:odd}. 
Note that $A_{{\mathbf K}}\setminus\{0\}$ consists of $6$ elements with $x^2 = 1 \in \ZZ/2\ZZ$ and $9$ elements $x$ with norm $x^2 = 0 \in \ZZ/2\ZZ$. 
Thus these $15$ involutions are divided into $6$ odd involutions and $9$ even involutions. 
\par
By Lemma~\ref{lemma:normalized:marking}, the following definition makes sense.

\begin{definition}
\label{def:Phi:gamma}
Let $\iota \colon K_{\tau, \tau'} \to K_{\tau, \tau'}$ be a fixed-point-free involution with \eqref{eqn:assumption:lattice} 
as in Theorem~\ref{thm:Mukai:Ohashi}.
Let $\gamma = \overline{d}_{\iota} \in A_{\mathbf K} \setminus \{0\}$ be the patching element of $\iota$. Define 
$$
\Phi_{\gamma}(\tau, \tau')^{2} := \Phi_{\ell}\left( j_{\ell}^{-1}\varpi( K_{\tau, \tau'}/\iota, \alpha )
 \right)^{2},
$$
where $\alpha$ is a normalized marking for $(K_{\tau, \tau'}, \iota)$ and $\ell$ is the level of $\Gamma_{34}^{\lor}$ in $H^{2}(K_{\tau,\tau'}, {\mathbf Z})_{-}$.
\end{definition}

Let $\gamma \in A_{\mathbf K} \setminus \{0\}$. As we see 
in Sections~\ref{sec:L} and \ref{subsec:realization:odd} below, 
there is a family of Kummer surfaces of product type $\pi \colon {\mathcal K} \to {\frak H} \times {\frak H}$ and fixed-point-free involutions 
$\iota_{\gamma} \colon {\mathcal K}|_{{\frak H} \times {\frak H} \setminus {\frak D}} \to {\mathcal K}|_{{\frak H} \times {\frak H} \setminus {\frak D}}$ 
preserving the fibers of $\pi$ and satisfying \eqref{eqn:assumption:lattice} 
such that the set of representatives of the $15$ conjugacy classes of Theorem~\ref{thm:Mukai:Ohashi} (1) is given by 
$\{\iota_{\gamma}|_{K_{\tau, \tau'}} \}_{\gamma \in A_{\mathbf K}\setminus\{0\}}$.
Since ${\frak H} \times {\frak H}$ is contractible, by choosing a reference point $(\tau_{0}, \tau'_{0}) \in {\frak H} \times {\frak H} \setminus {\frak D}$
and choosing a normalized marking $\alpha_{\gamma, 0}$ for $(K_{\tau_{0},\tau'_{0}}, \iota_{\gamma})$, we have a marking 
$\alpha_{\gamma} \colon R^{2}\pi_{*}{\mathbf Z} \cong \LL$ for $\pi \colon {\mathcal K} \to {\frak H} \times {\frak H}$
extending $\alpha_{\gamma,0}$ such that $\alpha_{\gamma}|_{K_{\tau,\tau'}}$ is a normalized marking
for $(K_{\tau, \tau'}, \iota_{\gamma})$ for all $(\tau, \tau') \in {\frak H} \times {\frak H} \setminus {\frak D}$. Write $\varpi_{\gamma}$ for the period mapping
for $(\pi \colon {\mathcal K} \to {\frak H} \times {\frak H}, \alpha_{\gamma})$:
$
\varpi_{\gamma}( \tau, \tau' ) := \varpi( K_{\tau, \tau'}/\iota_{\gamma}, \alpha_{\gamma} )
.
$
Then $\varpi_{\gamma}$ is a holomorphic map from ${\frak H} \times {\frak H}$ to $\Omega_{\alpha_{\gamma}({\mathbf K})}$ such that 
$$
\Phi_{\gamma}(\tau, \tau')^{2} = \Phi_{\ell}\left( j_{\ell}^{-1}\varpi_{\gamma}( \tau, \tau' ) \right)^{2}
$$ 
for all $(\tau, \tau') \in {\frak H} \times {\frak H} \setminus {\frak D}$. 
Since $\Phi_{\ell}^{2}\circ j_{\ell}^{-1}\circ \varpi_{\gamma}$ is a holomorphic function on ${\frak H} \times {\frak H}$, so is $\Phi_{\gamma}^{2}$.

\subsection{Automorphy of $\Phi_{\gamma}^{2}$}
\label{subsec:actions:period:map}
In Sect.~\ref{subsec:actions:period:map}, we prove that $\Phi_{\gamma}^{2}$ is an automorphic form of weight $8$ for the principal 
congruence subgroup of level $2$ of ${\rm SL}_{2}({\mathbf Z}) \times {\rm SL}_{2}( {\mathbf Z})$. We keep the notation in Sect.~\ref{sect:parity:invol}.
For $\gamma \in A_{\mathbf K}\setminus\{0\}$, we set
$$
{\mathbb K}_{\gamma} := \alpha_{\gamma}({\mathbf K}),
\qquad
{\mathbb E}_{\gamma} := \alpha_{\gamma}({\mathbf K}^{\perp_{H^{2}_{-}}}).
$$
Since $H^{1}(E_{\tau},{\bf Z})$ is endowed with the basis $\{{\bf a}^{\lor},{\bf b}^{\lor}\}$, ${\rm SL}_{2}({\bf Z})$ acts on $H^{1}(E_{\tau},{\bf Z})$ by 
$g\cdot(m{\bf a}^{\lor}+n{\bf b}^{\lor}) \colonequals ({\bf a}^{\lor},{\bf b}^{\lor})g\binom{m}{n}$ for $g \in {\rm SL}_{2}({\bf Z})$. 
Since this action preserves the cup-product, so does the induced
${\rm SL}_{2}({\bf Z})\times{\rm SL}_{2}({\bf Z})$-action on 
${\mathbf K} = \phi(H^{1}(E_{\tau},{\bf Z})\otimes H^{1}(E_{\tau'},{\bf Z}))$.

We define a map $\rho\colon{\rm SL}_{2}({\bf Z})\times{\rm SL}_{2}({\bf Z})\to O^{+}({\mathbf K})$ as  
\[
\rho(g,h)\cdot\phi({\bf u}\otimes{\bf v}) \colonequals \phi({}^{t}g^{-1}\cdot({\bf u})\otimes{}^{t}h^{-1}\cdot({\bf v}))
\]
for any $g, h \in \SL_2(\ZZ)$ and 
${\bf u}\in H^{1}(E_{\tau},{\bf Z})$, ${\bf v}\in H^{1}(E_{\tau'},{\bf Z})$. 
Then $\rho$ is a group homomorphism such that the period map 
$$
\varpi\colon{\mathfrak H}\times{\mathfrak H} 
\ni (\tau,\tau') \to
\left[ \tau\tau'\Gamma_{34}^{\lor} + \Gamma_{12}^{\lor} + \tau\Gamma_{23}^{\lor} + \tau'\Gamma_{14}^{\lor} \right] \in 
\Omega_{{\mathbf K}}^{+}
$$ 
is ${\rm SL}_{2}({\bf Z})\times{\rm SL}_{2}({\bf Z})$-equivariant, i.e., 
\begin{equation}
\label{eqn:equivariance:SL*SL-action}
\rho(g, g')\cdot\varpi(\tau,\tau')=\varpi(g\cdot\tau, g'\cdot\tau'),
\end{equation}
for any $g, g'\in{\rm SL}_{2}({\bf Z})$ and $\tau, \tau' \in {\mathfrak H}$. 
To summarize, for any involution $\iota_{\gamma}$ on ${\mathcal K}|_{{\frak H}\times{\frak H}\setminus{\frak D}}$ satisfying \eqref{eqn:assumption:lattice}
as in Theorem~\ref{thm:Mukai:Ohashi}, under the identification $\alpha_{\gamma}\colon {\mathbf K} \cong \KK_{\gamma}$ via a normalized marking for
 $({\mathcal K}, \iota_{\gamma})$, 
$\widetilde{\rho} := \alpha_{\gamma}\rho\alpha_{\gamma}^{-1} \colon{\rm SL}_{2}({\bf Z})\times{\rm SL}_{2}({\bf Z})\to O^{+}(\KK_{\gamma})$ 
is a group homomorphism such that the period map 
$\varpi_{\gamma} = \alpha_{\gamma}\circ\varpi\colon{\mathfrak H}\times{\mathfrak H}\to\Omega_{\KK_{\gamma}}^{+}$ in \eqref{eqn:(8.3)} 
is ${\rm SL}_{2}({\bf Z})\times{\rm SL}_{2}({\bf Z})$-equivariant. 

However, this does {\em not} imply that the period map \eqref{eqn:(8.3)} extends to an ${\rm SL}_{2}({\bf Z})\times{\rm SL}_{2}({\bf Z})$-equivariant holomorphic map 
from ${\mathfrak H}\times{\mathfrak H}$ to $\Omega_{\LAM}^{+}$, because $O^{+}({\mathbb K}_{\gamma})$ is {\em not} a subgroup of $O^{+}(\LAM)$.

Let $\widetilde{O}^{+}({\mathbb K}_{\gamma})$ be the kernel of the canonical homomorphism 
$O^{+}({\mathbb K}_{\gamma})\to O(q_{{\mathbb K}_{\gamma}})$.
Since $g\oplus1_{{\mathbb E}_{\gamma}}\in O^{+}({\mathbb K}_{\gamma}\oplus{\mathbb E}_{\gamma})$ preserves $\LAM$ 
for any $g\in\widetilde{O}^{+}({\mathbb K}_{\gamma})$ by the expression \eqref{eqn:the:choice:of:d}, we get the inclusion 
$\widetilde{O}^{+}({\mathbb K}_{\gamma})\ni g \mapsto g\oplus1_{{\mathbb E}_{\gamma}}\in O^{+}(\LAM)$.

Recall that $\Gamma(2) \subset \SL_2(\ZZ)$ is the principal  congruence subgroup of level $2$. 
The image of $\Gamma(2)\times\Gamma(2)$ under 
$\widetilde{\rho}\colon {\rm SL}_{2}({\bf Z})\times{\rm SL}_{2}({\bf Z})\to O^{+}({\mathbb K}_{\gamma})$ 
is contained in $\widetilde{O}^{+}({\mathbb K}_{\gamma})$. 
It follows that the period map 
$
\varpi\colon{\mathfrak H}\times{\mathfrak H}\to\Omega_{\LAM}^{+}, 
$
where the target space is now $\Omega_{\LAM}^{+}$, is $\Gamma(2)\times\Gamma(2)$-equivariant 
with respect to the homomorphism $\Gamma(2)\times\Gamma(2)\to O^{+}(\LAM)$.

\begin{lemma}
\label{lemma:automorphy:Phi:gamma}
For any $\gamma \in A_{\mathbf K} \setminus\{0\}$, $\Phi_{\gamma}^{2}$ is an automorphic form on ${\frak H}\times{\frak H}$ of weight $8$ for $\Gamma(2)\times\Gamma(2)$. 
Namely, the following functional equation holds for all $(g,g') = ( \binom{a\,b}{c\,d}, \binom{a'\,b'}{c'\,d'} )\in \Gamma(2)\times\Gamma(2)$ and $(\tau, \tau') \in {\frak H}\times{\frak H}$: 
$$
\Phi_{\gamma}(g\cdot\tau, g'\cdot\tau')^{2} = (c\tau+d)^{8}(c'\tau'+d')^{8} \Phi_{\gamma}(\tau, \tau')^{2}.
$$
\end{lemma}

\begin{pf}
Write $\varpi_{\gamma}(\tau, \tau') = \varpi( K_{\tau, \tau'}/\iota, \alpha )
 = j_{\ell}(u) =\left[ -(u^{2}/2){\ebf}_{\ell} +({\fbf}_{\ell}/\ell) + (-1)^{2/\ell} u \right]$,
where $\ell$ is the level of $\Gamma_{34}^{\lor}$ in $H^{2}(K_{\tau, \tau'}, {\mathbf Z})_{-}$.
Then $\Phi_{\gamma}(\tau, \tau')^{2} = \Phi_{\ell}( u )^{2}$ by the definition of $\Phi_{\gamma}$. Let $(g,g') \in \Gamma(2)\times\Gamma(2)$. 
By the $\Gamma(2)\times\Gamma(2)$-equivariance of $\widetilde{\rho}$, we get 
$$
\Phi_{\gamma}(g\cdot\tau, g'\cdot\tau')^{2} = \Phi_{\ell}( j_{\ell}^{-1}\varpi_{\gamma}(g\cdot\tau, g'\cdot\tau') )^{2} 
= \Phi_{\ell}( j_{\ell}^{-1}\widetilde{\rho}(g,g')\varpi_{\gamma}(\tau, \tau') )^{2} 
= \Phi_{\ell}( \widetilde{\rho}(g,g') \cdot u ).
$$
By the automorphy of $\Phi_{\ell}$ (see \eqref{eqn:automorphy:Phi}), we get 
\begin{equation}
\label{eqn:automorphy:1}
\Phi_{\gamma}(g\cdot\tau, g'\cdot\tau')^{2} = J_{\ell}(\widetilde{\rho}(g,g'), u)^{8}\,\Phi_{\ell}( u )^{2} 
=
J_{\ell}(\widetilde{\rho}(g,g'), u)^{8}\,\Phi_{\gamma}( \tau, \tau' )^{2}.
\end{equation}
Since 
$\langle \tau\tau'\alpha(\Gamma_{34}^{\lor}) + \alpha(\Gamma_{12}^{\lor}) + \tau\alpha(\Gamma_{23}^{\lor}) - \tau'\alpha(\Gamma_{14}^{\lor}), {\ebf}_{\ell}\rangle = -2$
by \eqref{eqn:marking:normalization}, we get
$$
\begin{aligned}
-(u^{2}/2){\ebf}_{\ell} +({\fbf}_{\ell}/\ell) + (-1)^{2/\ell} u 
&=
\frac{ \tau\tau'\alpha(\Gamma_{34}^{\lor}) + \alpha(\Gamma_{12}^{\lor}) + \tau\alpha(\Gamma_{23}^{\lor}) - \tau'\alpha(\Gamma_{14}^{\lor})}
{\langle \tau\tau'\alpha(\Gamma_{34}^{\lor}) + \alpha(\Gamma_{12}^{\lor}) + \tau\alpha(\Gamma_{23}^{\lor}) - \tau'\alpha(\Gamma_{14}^{\lor}), {\ebf}_{\ell}\rangle_{\LAM}}
\\
&=
\frac{1}{2}\{ \tau\tau'\alpha(\Gamma_{34}^{\lor}) + \alpha(\Gamma_{12}^{\lor}) + \tau\alpha(\Gamma_{23}^{\lor}) - \tau'\alpha(\Gamma_{14}^{\lor}) \},
\end{aligned}
$$
where we used \eqref{eqn:(8.3)} to get the first equality. 
By \eqref{eqn:identification:lattices}, \eqref{eqn:marking:normalization}, 
$$
\begin{aligned}
\,&
J_{\ell}( \widetilde{\rho}(g, g'), u )
=
\langle \widetilde{\rho}(g, g')\{ -(u^{2}/2){\ebf}_{\ell} +({\fbf}_{\ell}/\ell) + (-1)^{2/\ell} u \}, {\ebf}_{\ell} \rangle_{\LAM}
\\
&=
-\frac{1}{2}
\langle \rho(g,g') ( \phi( ({\bf a}^{\lor} - \tau {\bf b}^{\lor} ) \otimes ({\bf a}^{\prime \lor} - \tau' {\bf b}^{\prime \lor} ) ), \Gamma_{34}^{\lor} \rangle
\end{aligned}
$$
Write $g = \binom{a\,b}{c\,d}$, $g' = \binom{a'\,b'}{c'\,d'}$. By the definition of the $\Gamma(2)\times\Gamma(2)$-action, we get 
\begin{equation}
\label{eqn:automorphy:2}
\begin{aligned}
\,&
J_{\ell}( \widetilde{\rho}(g, g'), u )
=
-\frac{1}{2}
\langle 
\phi \left( ({\bf a}^{\lor}, {\bf b}^{\lor}) {}^{t}g^{-1}\binom{1}{-\tau} \otimes ({\bf a}^{\prime\lor}, {\bf b}^{\prime\lor}) {}^{t}g^{\prime -1}\binom{1}{-\tau'} \right), 
\Gamma_{34}^{\lor} 
\rangle
\\
&=
-\frac{1}{2}
\langle
(c\tau+d)(c'\tau'+d') \Gamma_{12}^{\lor} - (c\tau+d)(a'\tau'+b') \Gamma_{14}^{\lor} 
\\
&\qquad 
+ (a\tau+b)(c'\tau'+d') \Gamma_{23}^{\lor} + (a\tau+b)(a'\tau'+b') \Gamma_{34}^{\lor},
\Gamma_{34}^{\lor} \rangle
=
(c\tau+d)(c'\tau'+d'). 
\end{aligned}
\end{equation}
Substituting \eqref{eqn:automorphy:2} into \eqref{eqn:automorphy:1}, we get the desired functional equation.
\end{pf}

For $u\in{\mathbb M}_{\ell}\otimes{\bf R}+i\,{\mathcal C}_{{\mathbb M}_{\ell}}^{+}$ and $(\tau,\tau')\in{\mathfrak H}\times{\mathfrak H}$, we set
$$
B_{{\mathbb M}_{\ell}}(u) \colonequals \langle\Ima\, u,\Ima\, u\rangle_{{\mathbb M}_{\ell}},
\qquad
B_{{\KK}}(\tau,\tau') \colonequals \Ima\,\tau\cdot\Ima\,\tau'.
$$
We define the Petersson norm of $\Phi_{\gamma}$ as 
$$
\| \Phi_{\gamma}(\tau, \tau') \|^{2} := (\Ima\,\tau\cdot\Ima\,\tau')^{4} |\Phi_{\gamma}(\tau, \tau')|^{2}.
$$ 
When $u = j_{\ell}^{-1}\varpi_{\gamma}(\tau, \tau') \in {\mathbb M}_{\ell}\otimes{\mathbf R} + i\,{\mathcal C}_{{\mathbb M}_{\ell}}^{+}$ is given by \eqref{eqn:formula:period}, 
since $\langle B, B \rangle = \langle D, D \rangle = 0$ and $\langle B, D \rangle =2$ by \eqref{eqn:marking:normalization}, 
we have $\langle\Ima\, u,\Ima\, u\rangle_{{\mathbb M}_{\ell}} = \Ima\,\tau\cdot\Ima\,\tau'$. 
By \eqref{eqn:def:Petersson}, we have 
\begin{equation}
\label{eqn:def:Petersson:2}
\|\Phi_\ell ( \varpi_{\gamma}(\tau, \tau'))\|^{2} 
= 
(\Ima\,\tau\cdot\Ima\,\tau')^{4} \left| \Phi_\ell ( j_{\ell}^{-1}\varpi_{\gamma}(\tau, \tau')) \right|^{2} 
= \| \Phi_{\gamma}(\tau, \tau') \|^{2}. 
\end{equation}
By \eqref{eqn:def:Petersson:2} and the $O^{+}(\LAM)$-invariance of $\| \Phi_{\gamma} \|^2$ or by the automorphy of $\Phi_{\gamma}$,
$\| \Phi_{\gamma} \|^2$ is a $\Gamma(2) \times \Gamma(2)$-invariant $C^{\infty}$ function on ${\frak H} \times {\frak H}$.

\subsection{The period map and the discriminant locus}
\label{sect:period:disc:locus}
In this subsection, we study those involutions $\iota$ in Theorem~\ref{thm:Mukai:Ohashi} with patching element $\gamma \in A_{{\mathbf K}}\setminus\{0\}$
and satisfying $\Hcal \cap \Omega_{\KK_{\gamma}}^+ \neq \emptyset$. Further, 
when $H_{\delta} \cap \Omega_{\KK_{\gamma}}^+ \neq \emptyset$, we determine which $d \in \Delta_{\LAM}$ 
satisfies $H_{\delta} \cap \Omega_{\KK_{\gamma}}^+ = H_d \cap \Omega_{\KK_{\gamma}}^+$. 
In \eqref{eqn:the:choice:of:d}, we write 
$$
\LAM = \ZZ (d_1 + d_2) + {\mathbb K}_{\gamma} \oplus{\mathbb E}_{\gamma} \subset {\mathbb K}_{\gamma}^\vee\oplus{\mathbb E}_{\gamma}^{\vee},
\qquad
d_1 \in {\mathbb K}_{\gamma}^{\vee}\setminus{\mathbb K}_{\gamma},
\quad
d_2 \in{\mathbb E}_{\gamma}^{\vee}\setminus{\mathbb E}_{\gamma}.
$$
Since $\KK_{\gamma}$ 
and $\EE_{\gamma}$ 
are $2$-elementary, we have $2 d_1 \in \KK_{\gamma}$ and $2 d_2 \in \EE_{\gamma}$. 

\begin{lemma}
\label{lemma:upper}
Let $\delta \in \D_{\LAM}$ be any root of $\LAM$. 
If $\delta = \delta_{{\mathbb K}_{\gamma}} + \delta_{{\mathbb E}_{\gamma}} \in {\mathbb K}_{\gamma}^{\vee}\oplus{\mathbb E}_{\gamma}^{\vee}$ 
is the orthogonal decomposition, then $\delta_{{\mathbb K}_{\gamma}} -  d_1 \in{\mathbb K}_{\gamma}$ and $\delta_{{\mathbb E}_{\gamma}} - d_2 \in{\mathbb E}_{\gamma}$. 
In particular, $\LAM = \ZZ \delta + {\mathbb K}_{\gamma} \oplus{\mathbb E}_{\gamma}$. 
\end{lemma}

\begin{pf}
Since $\delta \in \LAM$, we can write $\delta = m (d_{1} + d_{2}) + k + e$, where $m \in \ZZ$, $k \in {\mathbb K}_{\gamma}$ and $e \in{\mathbb E}_{\gamma}$. 
Then $\delta_{{\mathbb K}_{\gamma}} = m d_{1} + k$ and $ \delta_{{\mathbb E}_{\gamma}} = md_{2} + e$.
Suppose that $\delta_{{\mathbb K}_{\gamma}} \in {\mathbb K}_{\gamma}$.  
Then $m$ is even whence $\delta_{{\mathbb E}_{\gamma}}  \in{\mathbb E}_{\gamma}$. 
It follows from the isometries ${\mathbb K}_{\gamma} \cong  \UU(2)\oplus\UU(2)$ 
and ${\mathbb E}_{\gamma} \cong \EE_8(2)$ that 
$\delta_{{\mathbb K}_{\gamma}}^2\in 4 \ZZ$ and $\delta_{{\mathbb E}_{\gamma}}^{2}\in4{\bf Z}$. 
This contradicts the equality $-2 = \delta^{2} = \delta_{{\mathbb K}_{\gamma}}^{2} + \delta_{{\mathbb E}_{\gamma}}^{2}$. 
Thus $\delta_{{\mathbb K}_{\gamma}} \not\in {\mathbb K}_{\gamma}$ and $m$ is odd. 
This proves the result.
\end{pf}

For $\delta \in \Delta_\LAM$, 
we write 
\[
\LAM = \ZZ (\delta_{{\mathbb K}_\gamma} + \delta_{{\mathbb E}_\gamma}) + {\mathbb K}_{\gamma} \oplus{\mathbb E}_{\gamma} 
\subset 
{\mathbb K}_{\gamma}^\vee\oplus{\mathbb E}_{\gamma}^{\vee},
\qquad
\delta_{{\mathbb K}_\gamma} \in {\mathbb K}_{\gamma}^{\vee}\setminus{\mathbb K}_{\gamma},
\quad
\delta_{{\mathbb E}_\gamma} \in{\mathbb E}_{\gamma}^{\vee}\setminus{\mathbb E}_{\gamma}.
\]

\begin{lemma}
\label{lemma:BP}
If $\iota$ is of odd type with patching element $\gamma$,
then there exists $\delta \in \Delta_{\LAM}$ such that 
$\delta_{{\mathbb K}_{\gamma}}^{2} = -1$ and $\delta_{{\mathbb E}_{\gamma}}^{2} = -1$.
\end{lemma}

\begin{pf}
Any element of $A_{\UU(2)\oplus\UU(2)}$ of odd norm is represented by one of the following vectors 
of $(\UU(2)\oplus\UU(2))^\vee$ of norm $-1$
\begin{gather*}
\left(1/2, -1/2, 0, 0\right), 
\left(1/2, -1/2, 1/2, 0\right), 
\left(1/2, -1/2, 0, 1/2 \right), \\
\left(0, 0, 1/2, -1/2\right), 
\left(1/2, 0, 1/2, -1/2\right), 
\left(0, 1/2, 1/2, -1/2\right).
\end{gather*}
Thus there exists $d_1' \in \KK_{\gamma}^\vee$ such that ${d'}_1^2 = -1$ and $d_1' -d_1 \in \KK_{\gamma}$. 
Similarly, since any element of $A_{\EE_8(2)}$ of odd norm is represented by a vector of $\EE_8(2)^\vee$ of norm $-1$ by \cite[Lemma~1.4, Cor.~1.5]{BP83},
there exists $d_2' \in \EE_{\gamma}^\vee$ such that ${d'}_2^2 = -1$ and $d_2' -d_2 \in \EE_{\gamma}$. 
We set $\delta \colonequals d_1' + d_2'$. 
It follows from $\delta - (d_1 + d_2) \in {\mathbb K}_{\gamma} \oplus{\mathbb E}_{\gamma}$ 
that  $\delta \in \LAM$. Further, $\delta^2 = {d'}_1^2 + {d'}_2^2 = -2$, so 
$\delta \in \Delta_\LAM$. 
\end{pf}

\begin{proposition}
\label{prop:discriminant:fixed:point}
The following hold:
\begin{enumerate}
\item
If $\iota$ is an involution of odd type with patching element $\gamma$, then 
$\Hcal \cap \Omega_{{\mathbb K}_{\gamma}}^+\not=\emptyset$.
\item
If $\iota$ is an involution of even type with patching element $\gamma$, then 
$\Hcal \cap \Omega_{{\mathbb K}_{\gamma}}^+=\emptyset$.
\end{enumerate}
\end{proposition}

\begin{pf}
(1) 
It suffices to prove the existence of $\delta\in\Delta_{\LAM}$ with $H_{\delta}\cap\Omega_{{\mathbb K}_{\gamma}}^{+}\not=\emptyset$.
Since $\iota$ is an involution of odd type, we take $\delta \in \Delta_{\LAM}$ 
as in Lemma~\ref{lemma:BP}.
Since $\delta\in\Delta_{\LAM}$, we get
$H_\delta\cap\Omega_{{\mathbb K}_{\gamma}}^{+}= H_{\delta_{{\mathbb K}_{\gamma}}}\cap\Omega_{{\mathbb K}_{\gamma}}^{+}\not=\emptyset$,
where the non-emptiness follows from $\delta_{{\mathbb K}_{\gamma}}^{2}=-1<0$ (see~\eqref{eqn:Hd:d2}).
This proves (1).

(2) To derive a contradiction, we assume that there exists a root $\delta \in \Delta_{\LAM}$ with $H_\delta\cap\Omega_{{\mathbb K}_{\gamma}}^{+}\not=\emptyset$.
Since $\iota$ is of even type, we get $\delta_{{\mathbb K}_{\gamma}}^2 \in 2 \ZZ$. 
Since $H_\delta\cap\Omega_{{\mathbb K}_{\gamma}}^+\not=\emptyset$, we have $\delta_{{\mathbb K}_{\gamma}}^2 < 0$ (see~\ref{eqn:Hd:d2}).  
On the other hand, 
since ${\mathbb E}_{\gamma}\cong{\mathbb E}_{8}(2)$ is negative-definite, we have $\delta_{{\mathbb E}_{\gamma}}^{2}\leq0$.
It then follows from the equality $-2 = \delta^2 = \delta_{{\mathbb K}_{\gamma}}^2 +  \delta_{{\mathbb E}_{\gamma}}^2$ that 
$-2 \leq \delta_{{\mathbb K}_{\gamma}}^2 < 0$. Thus 
$\delta_{{\mathbb K}_{\gamma}}^2 = -2$ and $\delta_{{\mathbb E}_{\gamma}}^2 = 0$. 
Since ${\mathbb E}_{\gamma}\cong{\mathbb E}_8(2)$ is negative-definite, 
we get  $\delta_{{\mathbb E}_{\gamma}} = 0$, so $\delta = \delta_{{\mathbb K}_{\gamma}}$. 
The primitivity of the embedding ${\mathbb K}_{\gamma} \subset \LAM$ then gives $\delta \in\Delta_{{\mathbb K}_{\gamma}}$. 
Since $\KK_\iota \cong \UU(2)\oplus\UU(2)$, 
this contradicts the fact $\Delta_{\UU(2)\oplus\UU(2)}=\emptyset$.
\end{pf}

Let $\mu\colon{\rm Km}(E\times E')\to{\rm Km}(E\times E')$ be the anti-symplectic holomorphic involution induced from 
the one on $E\times E'$ defined as  $(x,y)\mapsto(-x,y)$. (We remark that $\mu$ is not fixed-point-free.) 
Let $I_{\mu}\in O({\mathbb L}_{K3})$ be the involution defined as $I_{\mu}=\alpha\mu^{*}\alpha^{-1}$.
By \cite[Prop.~6]{Mukai10}, we get the following: 
\begin{equation}
\label{lemma:involution:varsigma}
{\mathbb K}_{\gamma}=\{\lambda\in{\mathbb L}_{K3}\mid I_{\mu}(\lambda)=-\lambda\},
\qquad
{\mathbb K}_{\gamma}^{\perp_{\LL}}=\{\lambda\in{\mathbb L}_{K3}\mid I_{\mu}(\lambda)=\lambda\}.
\end{equation}

\begin{proposition}
\label{prop:characterization:zero:Kondo-Mukai}
Suppose that $\iota$ is an involution of odd type with patching element $\gamma$. 
Let $d,\delta\in\Delta_{\LAM}$ with $\delta_{\KK_{\gamma}}^2 < 0$. Then 
$$
H_d \cap\Omega_{{\mathbb K}_{\gamma}}^+=H_\delta\cap\Omega_{{\mathbb K}_{\gamma}}^+
\qquad\text{if and only if}\qquad
d\in\{\pm\delta,\pm I_{\mu}(\delta)\}.
$$
\end{proposition}

\begin{pf}
In the proof, for simplicity, we write ${\mathbb K}$ (resp. ${\mathbb E}$) for ${\mathbb K}_{\gamma}$ (resp. ${\mathbb E}_{\gamma}$).
We write $d = d_{\mathbb K} + d_{\mathbb E}$ and $\d = \d_{\mathbb K} + \d_{\mathbb E}$, 
where $d_{\mathbb K}, \d_{\mathbb K}\in{\mathbb K}^{\vee}$ and $d_{\mathbb E}, \d_{\mathbb E}\in{\mathbb E}^{\vee}$. 
Note that, since $\delta_{\KK_\iota}^2 < 0$, we have $\delta^{\perp}\cap\Omega_{\mathbb K}^+ \neq \emptyset$ (see~\eqref{eqn:Hd:d2}).  
First we show the ``if'' part. We assume $d\in\{\pm\delta,\pm I_{\mu}(\delta)\}$.
Since $I_{\mu}(\d) = -\delta_{\mathbb K} + \delta_{\mathbb E}$, we get 
$H_\d \cap\Omega_{\mathbb K}^+ = H_{\delta_{\mathbb K}} \cap \O_{\mathbb K}^+ 
= H_{I_{\mu}(\d)} \cap\O_{\mathbb K}^+$. Since $H_{-\d} = H_{\d}$ and 
$H_{I_{\mu}(\d)} = H_{-I_{\mu}(\d)}$, we obtain the assertion. 

To prove the ``only if'' part, assume $d\not=\pm\delta$ and 
$d^{\perp}\cap\Omega_{\mathbb K}^+=\delta^{\perp}\cap\Omega_{\mathbb K}^+$. 
Then ${\bf Z}d+{\bf Z}\delta\subset\LAM$ is a sublattice of rank $2$.

\medskip
{\sl Step 1.}\;
We show $\langle d,\delta\rangle=0$. 
Since 
$H_d \cap H_\delta \cap\Omega_{\LAM}^+ = H_{\d}\cap\Omega_{\LAM}^+ \not=\emptyset$, 
the sublattice ${\bf Z}d+{\bf Z}\delta\subset\LAM$ is negative-definite by \eqref{eqn:Hd:d2}.
Thus $(d \pm \delta)^2 < 0$, so $-2 < \langle d,\delta\rangle < 2$. 
Since $\langle d,\delta\rangle\in2{\bf Z}$ for all $d,\delta\in\Delta_{\LAM}$ (cf. \cite[proof of Th.\,4.7]{Yoshikawa09}), 
we get $\langle d,\delta\rangle=0$.

\medskip
{\sl Step 2.}\;
We show $d_{\mathbb K} = \pm \delta_{\mathbb K}$. 
Since $H_d\cap\Omega_{\mathbb K}^+=\Omega_{d_{\mathbb K}^{\perp}\cap{\mathbb K}}^+$ and
$H_{\delta}\cap\Omega_{\mathbb K}^+=\Omega_{\delta_{\mathbb K}^{\perp}\cap{\mathbb K}}^+$, we get
$\Omega_{d_{\mathbb K}^{\perp}\cap{\mathbb K}}^+=\Omega_{\delta_{\mathbb K}^{\perp}\cap{\mathbb K}}^+ \neq \emptyset$ by the assumption.
Hence $d_{\mathbb K}=t\,\delta_{\mathbb K}$ for some $t\in{\bf Q}$.
Since $\Omega_{d_{\mathbb K}^{\perp}\cap{\mathbb K}}^+\not=\emptyset$ and
$\Omega_{\delta_{\mathbb K}^{\perp}\cap{\mathbb K}}^+\not=\emptyset$, 
we get $d_{\mathbb K}^{2}<0$ and $\delta_{\mathbb K}^{2}<0$ by \eqref{eqn:Hd:d2}.
Since $q_{\mathbb K}$ is $\ZZ/2\ZZ$-valued, we have 
$d_{\mathbb K}^{2}\in{\bf Z}_{<0}$ and $\d_{\mathbb K}^{2}\in{\bf Z}_{<0}$.

On the other hand, it follows from $\Delta_{\mathbb K}=\emptyset$ that $d_{\mathbb E}\not=0$ and $\delta_{\mathbb E}\not=0$. 
Since ${\mathbb E}$ is negative-definite, we get $d_{\mathbb E}^{2}<0$ and $\delta_{\mathbb E}^{2}<0$.
By the conditions
$-2=d^{2}=(d_{\mathbb K})^{2}+(d_{\mathbb E})^{2}$,
$d_{\mathbb K}^{2}\in{\bf Z}_{<0}$, 
$d_{\mathbb E}^{2}\in{\bf Z}_{<0}$,
we get $d_{\mathbb K}^{2}=d_{\mathbb E}^{2}=-1$. Similarly, $\delta_{\mathbb K}^{2}=\delta_{\mathbb E}^{2}=-1$.
Since $d_{\mathbb K}=t\delta_{\mathbb K}$, we get $d_{\mathbb K}=\pm\delta_{\mathbb K}$
by the equality $d_{\mathbb K}^{2}=\delta_{\mathbb K}^{2}=-1$.

\medskip
{\sl Step 3.}\;
We show $d_{\mathbb E}=\pm \delta_{\mathbb E}$. 
Indeed, since
$0=\langle d,\delta\rangle=\langle d_{\mathbb K},\delta_{\mathbb K}\rangle+\langle d_{\mathbb E},\delta_{\mathbb E}\rangle$ 
and $d_{\mathbb K}=\pm\delta_{\mathbb K}$, 
we get $\langle d_{\mathbb E},\delta_{\mathbb E}\rangle=\mp d_{\mathbb K}^{2}=\pm 1$ and thus 
\begin{equation}
\label{eqn:matrix}
\begin{pmatrix}
d_{\mathbb E}^{2}&\langle d_{\mathbb E},\delta_{\mathbb E}\rangle
\\
\langle d_{\mathbb E},\delta_{\mathbb E}\rangle& \delta_{\mathbb E}^{2}
\end{pmatrix}
=
\begin{pmatrix}
-1&\pm1
\\
\pm1&-1
\end{pmatrix}.
\end{equation}
If $d_{\mathbb E}\not=\pm\delta_{\mathbb E}$,
then ${\bf Q}d_{\mathbb E}+{\bf Q}\delta_{\mathbb E}\subset{\mathbb E}\otimes{\bf Q}$ is a $2$-dimensional subspace, 
which is {\it not} negative-definite by \eqref{eqn:matrix}. 
This contradicts the fact that ${\mathbb E}$ is negative-definite. Hence $d_{\mathbb E}=\pm\delta_{\mathbb E}$.
Since we assume $d\not=\pm\delta$, we get $d=\pm I_{\mu}(\delta)$ (see \eqref{lemma:involution:varsigma}).
This proves $d=\pm I_{\mu}(\delta)$.
\end{pf}

\section{Involutions of even type}
\label{sec:L}

In this section, we consider involutions of even type on $K_{\tau,\tau'}$. The main result of this section is 
Theorem~\ref{prop:c:-3} below, which is essentially shown in \cite[Cor.~7.6]{KawaguchiMukaiYoshikawa18} 
as a consequence of an algebraic expression of the Borcherds $\Phi$-function. 
For details of this section, we refer the reader to \cite[Sect.~7.2]{KawaguchiMukaiYoshikawa18}. 

For $\tau \in \HH$, let 
$\theta_{2}(\tau) := \sum_{n\in{\mathbf Z}} e^{\pi i (n+\frac{1}{2})^{2}}$,
$\theta_{3}(\tau) := \sum_{n\in{\mathbf Z}} e^{\pi i n^{2}}$
be the theta constants. 
For $\tau,\tau'\in{\mathfrak H}$, we set 
$$
\lambda \colonequals \theta_{2}(\tau)^{4}/\theta_{3}(\tau)^{4}, 
\qquad
\lambda' \colonequals \theta_{2}(\tau')^{4}/\theta_{3}(\tau')^{4}.
$$
Let ${\mathcal X} \to {\frak H}\times{\frak H}$ be the family of surfaces over ${\frak H}\times{\frak H}$ defined as
$$
{\mathcal X} :=
\left\{
(x,(\tau,\tau')) \in {\bf P}^{5} \times{\frak H}\times{\frak H}
\;\left|\;
\begin{array}{rl}
(1-\lambda)x_{1}^{2}+\lambda x_{2}^{2}-x_{3}^{2}
&=
0,
\\
\lambda' x_{1}^{2}-\lambda' x_{2}^{2} - x_{4}^{2} + x_{6}^{2}
&=
0,
\\
x_{1}^{2} - x_{2}^{2} - x_{4}^{2}+x_{5}^{2}
&=
0
\end{array}
\right.
\right\}.
$$
Let $X_{\lambda,\lambda'}$ be its fiber over $(\tau,\tau')\in{\mathfrak H}\times{\mathfrak H}$.
Then $\varSigma = \{ (0:0:0:1:\pm1:\pm1),\, (1:\pm1:\pm1:0:0:0)\}$ is the singular locus of $X_{\lambda,\lambda'}$, which consists of ordinary double points, and 
the minimal resolution $\widetilde{X}_{\lambda,\lambda'}$ of $X_{\lambda,\lambda'}$ is isomorphic to the Kummer surface $K_{\tau,\tau'}={\rm Km}(E_{\tau}\times E_{\tau'})$. 
Let ${\mathcal K} \to {\mathcal X}$ be the blowing-up of ${\mathcal X} \subset {\mathbf P}^{5}\times{\frak H}\times{\frak H}$ along the loci 
$\amalg_{p\in \varSigma} \{p\} \times{\frak H}\times{\frak H}$ and let $\pi \colon {\mathcal K} \to {\frak H}\times{\frak H}$ be the obvious projection. 
By construction, 
$
\pi \colon {\mathcal K} \to {\frak H}\times{\frak H}
$
is a family of Kummer surfaces such that $\pi^{-1}(\tau, \tau') \cong K_{\tau,\tau'}$ for all $(\tau,\tau')\in{\frak H}\times{\frak H}$. 

For $J=\{j_{1},j_{2},j_{3}\}\subset\{1,\ldots,6\}$ with $j_{1}<j_{2}<j_{3}$, let $\langle J\rangle$ denote the partition $\{1,\ldots,6\}=J\amalg J^{c}$. 
We write $J^{c}=\{j_{4},j_{5},j_{6}\}$ with $j_{4}<j_{5}<j_{6}$. 
Note that $\langle J\rangle = \langle J^c\rangle$.  
We also denote $\langle J\rangle$ by $\binom{j_{1}j_{2}j_{3}}{j_{4}j_{5}j_{6}}$. 
For a partition $\langle J\rangle=\binom{j_{1}j_{2}j_{3}}{j_{4}j_{5}j_{6}}$,
we define the involution $\iota_{\langle J\rangle}$ on ${\bf P}^{5}$ as 
$$
\iota_{\langle J\rangle}(x_{j_{1}},x_{j_{2}},x_{j_{3}},x_{j_{4}},x_{j_{5}},x_{j_{6}})
 \colonequals 
(x_{j_{1}},x_{j_{2}},x_{j_{3}},-x_{j_{4}},-x_{j_{5}},-x_{j_{6}}).
$$
Then $\iota_{\langle J\rangle}$ acts on $X_{\lambda,\lambda'}$. 
For $\langle J\rangle\not=\binom{123}{456}$, $\iota_{\langle J\rangle}$ is fixed-point-free on $X_{\lambda,\lambda'}$.
Since $K_{\tau,\tau'}=\widetilde{X}_{\lambda,\lambda'}$, the involution $\iota_{\langle J\rangle}$ acts on ${\mathcal K}$
and preserves the fibers of $\pi$. 
By Lemma~\ref{lemma:extendable:involution}, $\iota_{\langle J\rangle}$ satisfies \eqref{eqn:assumption:lattice}. 
By \cite[Sect.\,7.2]{KawaguchiMukaiYoshikawa18}, the involution on $K_{\tau,\tau'}$ induced by $\iota_{\langle J\rangle}$ 
is an 
involution of even type, which is 
again denoted by $\iota_{\langle J\rangle}$,
and these $9$ involutions $\{\iota_{\langle J\rangle}\}_{\langle J\rangle\not=\binom{123}{456}}$ 
give the $9$ conjugacy classes of involutions of even type on $K_{\tau,\tau'}$ (cf.~Theorem~\ref{thm:Mukai:Ohashi}). 
Since $\iota_{\langle J\rangle}$ satisfies \eqref{eqn:assumption:lattice}, 
we define $\ell(J)\in\{1,2\}$ as the level of the primitive isotropic vector $\Gamma_{34}^{\lor}$ in $H^{2}(K_{\tau, \tau'}, {\mathbf Z})_{-}$. 

By Lemma~\ref{lemma:normalized:marking}, $K_{\tau, \tau'}$ admits a normalized marking.
By the triviality of the local system $R^{2}\pi_{*}{\bf Z}$, this marking extends to the one for the family $\pi\colon{\mathcal K}\to{\mathfrak H}\times{\mathfrak H}$,
which is fiberwise normalized with respect to $\iota_{\langle J\rangle}$.
Let $\alpha_{\langle J\rangle}\colon R^{2}\pi_{*}{\bf Z}\cong{\mathbb L}_{K3}$ be a normalized marking obtained in this way.
We set ${\mathbb K}_{\langle J\rangle} \colonequals {\mathbb K}_{\iota_{\langle J\rangle}}=\alpha_{\langle J\rangle}({\mathbf K})$. Then
$\Omega_{{\mathbb K}_{\langle J\rangle}}=\{[\omega] \in\Omega_{\LAM}^{+}\mid \omega\in{\mathbb K}_{\langle J\rangle}\otimes{\bf C}\}$ is the period domain for 
the marked family of Enriques surfaces $(\pi\colon\widetilde{\mathcal X}/\iota_{\langle J\rangle}\to{\mathfrak H}\times{\mathfrak H},\alpha_{\langle J\rangle})$.
Let 
\begin{equation}
\label{eqn:varpi:J}
\varpi_{\langle J\rangle}\colon{\mathfrak H}\times{\mathfrak H}\to\Omega_{{\mathbb K}_{\langle J\rangle}}^{+}
\end{equation}
be its period map. 
Let $\gamma \in A_{\mathbf K} \setminus\{ 0 \}$ be the patching element of $\iota_{\langle J \rangle}$. Then
$$
\Phi_{\gamma}(\tau, \tau')^{2} 
= 
\Phi_{\ell(J)}\left(  j_{\ell(J)}^{-1}(\varpi_{\langle J\rangle}(\tau, \tau')) \right)^{2}
=
\Phi_{\ell(J)}\left(  j_{\ell(J)}^{-1}\varpi( K_{\tau, \tau'}/\iota_{\langle J\rangle}, \alpha_{\langle J\rangle}) \right)^{2}
$$
by Definition~\ref{def:Phi:gamma}. By Lemma~\ref{lemma:normalized:marking} (2), $\Phi_{\gamma}^{2}$ is independent of the choice of a normalized marking
$\alpha_{\langle J \rangle}$ and is a holomorphic function on ${\frak H}\times{\frak H}$, which is an automorphic form for $\Gamma(2)\times\Gamma(2)$ of weight $8$ 
by Lemma~\ref{lemma:automorphy:Phi:gamma}. 
By \cite{KawaguchiMukaiYoshikawa18}, we have the following

\begin{theorem}
\label{prop:c:-3}
For any $(\tau, \tau') \in {\mathfrak H}\times{\mathfrak H}$, we have 
$$
\prod_{\gamma\,{\rm even}} \Phi_{\gamma}(\tau, \tau')^{2} 
=
\prod_{\langle J\rangle\not=\binom{123}{456}}\Phi_{\ell(J)}\left( j_{\ell(J)}^{-1} \varpi_{\langle J\rangle}(\tau,\tau') \right)^2
=
2^{96}\eta^{144}(\tau)\eta^{144}(\tau').
$$
\end{theorem}

\begin{pf}
The result follows from \cite[Cor.\,7.6]{KawaguchiMukaiYoshikawa18}.
\end{pf}

\section{Involutions of odd type} 
\label{sec:KM}
In this section, for involutions of odd type, we give their geometric realizations and study the 
extended period map. Then we study behavior of the Borcherds $\Phi$-function along 
the boundary of the extended period map.

\subsection{Realization of the involutions of odd type}
\label{subsec:realization:odd}

For $\lambda \in {\mathbf C}\setminus \{0,1\}$, let $E(\lambda)$ be the elliptic curve defined as 
the double covering of ${\mathbf P}^{1}$ with ordered four branch points 
$(0, 1, \lambda, \infty)$. 
Let $\lambda, \lambda' \in {\mathbf C}\setminus \{0,1\}$ and set
$p_{1} := (0,0)$, $p_{2} := (1,1)$, $p_{3} := (\lambda, \lambda')$, $p_{4} := (\infty, \infty) \in {\mathbf P}^{1}\times{\mathbf P}^{1}$.
Then ${\mathbf P}^{1}\times{\mathbf P}^{1}$ with ordered four points $p_{1},p_{2},p_{3},p_{4}$ is isomorphic to the quadric $Q_{\lambda,\lambda'} $
of ${\mathbf P}^{3}$ with 
ordered four points  
$(1:0:0:0)$, $(0:1:0:0)$, $(0:0:1:0)$, $(0:0:0:1)$, 
where $Q_{\lambda,\lambda'}$ is defined by the equation
\begin{equation}
\label{eqn:defining:eq}
\lambda\lambda' x_{2}x_{3} + x_{1}x_{3} + (1-\lambda)(1-\lambda') x_{1}x_{2} + x_{4} ( x_{1} + x_{2} + x_{3} ) =0.
\end{equation}

When $\lambda \not= \lambda'$, 
$Q_{\lambda, \lambda'}$ is a non-singular quadric. 
When $\lambda= \lambda'$, 
$Q_{\lambda, \lambda}$ is a singular quadric with a unique ordinary double point.
\par
Let ${\mathcal W}$ be the subvariety of ${\mathbf P}(1:1:1:1:2) \times ({\mathbf C}\setminus\{0,1\})^{2}$ defined by the equations 
\eqref{eqn:defining:eq} and $w^{2} = x_{1}x_{2}x_{3}x_{4}$, which is endowed with the obvious projection
${\mathcal W} \to ({\mathbf C}\setminus\{0,1\})^{2}$.
Let $W_{\lambda, \lambda'}$ be the fiber of ${\mathcal W}$ over $(\lambda, \lambda') \in ({\mathbf C}\setminus\{0,1\})^{2}$.
Then $W_{\lambda,\lambda'}$ is the double covering of $Q_{\lambda,\lambda'}$ with branch divisor $x_{1}x_{2}x_{3}x_{4}=0$.
When $\lambda\not=\lambda'$, ${\rm Sing}\,W_{\lambda,\lambda'}$ consists of four $D_{4}$-singularities. When $\lambda=\lambda'$,
${\rm Sing}\,W_{\lambda,\lambda'}$ consists of four $D_{4}$-singularities and two $A_{1}$-singularities. As is easily verified, we can take a simultaneous resolution
of ${\mathcal W}\to ({\mathbf C}\setminus\{0,1\})^{2}$,
which we write
$$
\pi \colon {\mathcal Y} \to ({\mathbf C}\setminus\{0,1\})^{2}.
$$ 
We set $Y_{\lambda,\lambda'} := \pi^{-1}(\lambda,\lambda')$. 
When $\lambda \neq \lambda'$, we have $Y_{\lambda,\lambda'} \cong {\rm Km}( E(\lambda) \times E(\lambda') )$
by \cite[Lemma 11]{Mukai10}.
Thus for any $(\lambda, \lambda') \in ({\mathbf C}\setminus\{0,1\})^{2}$, the weak Torelli theorem implies that 
$Y_{\lambda,\lambda'} \cong {\rm Km}( E(\lambda) \times E(\lambda') )$.
\par
Let $Z \subset ({\mathbf C}\setminus\{0,1\})^{2}$ be the diagonal locus. 
Let $\iota$ 
be the involution on ${\mathcal Y}|_{({\mathbf C}\setminus\{0,1\})^{2}\setminus Z}$ induced by  
the rational involution on ${\mathcal W}|_{({\mathbf C}\setminus\{0,1\})^{2}\setminus Z}$
\begin{equation}
\label{eqn:iota:l:l'}
(x, w) 
\mapsto 
\left( \frac{\lambda\lambda'}{x_{1}}, \frac{1}{x_{2}}, \frac{(1-\lambda)(1-\lambda')}{x_{3}}, 
\frac{\lambda\lambda'(1-\lambda)(1-\lambda')}{x_{4}}, \frac{\lambda\lambda'(1-\lambda)(1-\lambda')}{w} \right).
\end{equation}
Then $\iota$ preserves the fibers of $\pi $.
By \cite[Lemma 16 (2)]{Mukai10}, $\iota$ is an involution of odd type in the sense of Definition~\ref{def:KM:L}.
By Lemma~\ref{lemma:extendable:involution}, $\iota$ satisfies \eqref{eqn:assumption:lattice}. 

\begin{remark}
\label{rmk:invol:not:extend}
When $\lambda = \lambda'$, it follows from \eqref{eqn:iota:l:l'} that 
$\iota|_{W_{\lambda,\lambda}}$ has two fixed points, i.e., 
the two $A_{1}$-singularities of $W_{\lambda,\lambda}$ lying over ${\rm Sing}\,Q_{\lambda, \lambda}$. Hence the involution on $Y_{\lambda, \lambda}$ 
induced by the rational involution \eqref{eqn:iota:l:l'} has non-empty fixed locus consisting of two disjoint $(-2)$-curves. 
This implies that $\iota$ 
does not lift to an involution on ${\mathcal Y}$.
\end{remark}

\par
Recall that $E(\lambda)$ is endowed with the ordered four points $(0, 1, \lambda,\infty)$. We set $z_{1} = 0$, $z_{2} = 1$, $z_{3} = \lambda$. 
It is classical that if we change the order of the points $z_{1},z_{2},z_{3}$ by a permutation $\varrho \in {\frak S}_{3}$, 
then the pair $(E(\lambda), (z_{\varrho(1)}, z_{\varrho(2)}, z_{\varrho(3)}, \infty))$ is isomorphic to $\left(E(\varrho(\lambda)), (0, 1, \varrho(\lambda),\infty)\right)$, 
where $\varrho(\lambda) := (z_{\varrho(3)}-z_{\varrho(1)})/(z_{\varrho(2)}-z_{\varrho(1)})$ is a linear fractional transformation preserving 
${\mathbf C} \setminus\{ 0, 1 \}$.
Since $E(\varrho(\lambda)) \cong E(\lambda)$, there is an isomorphism 
$\psi_{\varrho,\varrho'} \colon Y_{\varrho(\lambda), \varrho'(\lambda')} \to {\rm Km}(E(\lambda)\times E(\lambda'))$, 
so that $\psi_{\varrho,\varrho'} \iota \psi_{\varrho,\varrho'}^{-1}$ is a fixed-point-free involution on ${\rm Km}(E(\lambda)\times E(\lambda'))$ 
when $\varrho(\lambda)\not=\varrho'(\lambda')$. 
By \eqref{eqn:defining:eq}, \eqref{eqn:iota:l:l'}, we easily see that
$(W_{\sigma(\lambda), \sigma(\lambda')}, \iota) \cong (W_{\lambda, \lambda'}, \iota)$ for all $\sigma \in {\mathfrak S}_{3}$.
Hence we have
\begin{equation}
\label{eqn:S3:invariance}
(Y_{\sigma(\lambda), \sigma(\lambda')}, \iota) \cong (Y_{\lambda, \lambda'}, \iota)
\qquad
(\forall\,\lambda,\lambda'\in {\mathbf C}\setminus\{0,1\},\,
\forall\,\sigma\in{\mathfrak S}_{3}).
\end{equation}
So we have the following isomorphism for all $\varrho, \varrho', \sigma \in {\frak S}_{3}$:
$$
( {\rm Km}(E(\lambda)\times E(\lambda')), \psi_{\varrho,\varrho'} \iota \psi_{\varrho,\varrho'}^{-1} )
\cong
( {\rm Km}(E(\lambda)\times E(\lambda')), \psi_{\sigma\varrho,\sigma\varrho'} \iota \psi_{\sigma\varrho,\sigma\varrho'}^{-1} ).
$$
We set $\psi_{\varrho} := \psi_{1, \varrho}$. 
By \cite[Lemma 16 (2)]{Mukai10}, \cite[Th.~4.1]{Ohashi07}, 
the involutions $\{ \psi_{\varrho} \iota \psi_{\varrho}^{-1} \}_{\varrho \in {\frak S}_{3}}$ 
give complete representatives of 
involutions of odd type on ${\rm Km}(E(\lambda)\times E(\lambda'))$. 
\par
For $\varrho\in{\mathfrak S}_{3}$, let $\pi_{\varrho} \colon {\mathcal Y}_{\varrho} \to ({\bf C}\setminus\{0,1\})^{2}$ be the family induced from
${\mathcal Y}$ by the map ${\rm id}\times \varrho \colon ({\bf C}\setminus\{0,1\})^{2} \to ({\bf C}\setminus\{0,1\})^{2}$.
Hence $\pi_{\varrho}^{-1}(\lambda,\lambda') = Y_{\lambda,\varrho(\lambda')}$. We set
$$
Z_{\varrho} \colonequals \{(\lambda,\lambda')\in({\bf C}\setminus\{0,1\})^{2}\mid \lambda = \varrho(\lambda') \}.
$$
Then the family 
$\pi_{\varrho} \colon {\mathcal Y}_{\varrho}|_{({\bf C}\setminus\{0,1\})^{2} \setminus Z_{\varrho}} \to ({\bf C}\setminus\{0,1\})^{2} \setminus Z_{\varrho}$
is endowed with the fixed-point-free involution $\iota_{\varrho}$ such that 
\begin{equation}
\label{eqn:def:iota:rho}
(Y_{\lambda, \lambda'}, \iota_{\varrho}) \cong (Y_{\lambda, \varrho(\lambda')}, \iota).
\end{equation}

\subsection{The period map}
\label{subsec:period:KM}

For $n\in{\bf Z}_{>0}$, recall 
that $\Gamma(n)\subset \SL_{2}({\bf Z})$ is the principal congruence subgroup of level $n$, 
which acts projectively on ${\mathfrak H}$. Note that $\Gamma(1) = \SL_2(\ZZ)$. 
Recall from \eqref{eqn:frakD} that ${\mathfrak D}
= 
\bigcup_{\gamma\in \Gamma(2)}
(\gamma\times 1)\varDelta_{\mathfrak H}
=
\bigcup_{\gamma\in \Gamma(2)}
(1\times\gamma)\varDelta_{\mathfrak H}$. 
For $\varrho\in{\mathfrak S}_{3}$, we define the reduced divisor on ${\mathfrak H}\times{\mathfrak H}$ as 
\begin{equation}
\label{eqn:frakD:gamma}
{\mathfrak D}_{\varrho} \colonequals (1\times\varrho){\mathfrak D} = (\varrho^{-1}\times1){\mathfrak D},
\end{equation}
where ${\mathfrak S}_{3}$ acts on ${\frak H}/\Gamma(2)$ via the standard isomorphism ${\mathfrak S}_3 \cong \Gamma(1)/\Gamma(2)$.

Recall that $\lambda(\tau) = \theta_{2}(\tau)^{4}/\theta_{3}(\tau)^{4}$.  
Let $\varPi \colon {\mathfrak H} \times {\mathfrak H} \to ( {\mathbf C} \setminus \{0,1\} )^{2}$ be the holomorphic map defined as  
$\varPi(\tau,\tau') := ( \lambda(\tau), \lambda(\tau') )$. Then $\varPi^{-1}( Z_{\varrho} ) = {\mathfrak D}_{\varrho}$. 
For $\varrho\in{\mathfrak S}_{3}$, we define a family of Kummer surfaces 
$$
\pi_{\varrho} \colon {\mathcal K}_{\varrho} \to {\mathfrak H}\times{\mathfrak H}
$$ 
as the pullback of $\pi_{\varrho} \colon {\mathcal Y}_{\varrho} \to ( {\mathbf C} \setminus \{0,1\} )^{2}$ by $\varPi$. 
Then it is a ``universal family'' in the sense that 
$\pi_{\varrho}^{-1}(\tau,\tau') = Y_{\lambda(\tau), \varrho(\lambda(\tau'))} \cong {\rm Km}\left(E(\lambda(\tau))\times E\left(\varrho(\lambda(\tau'))\right)\right) 
\cong K_{\tau, \tau'}$ for all $(\tau,\tau') \in {\mathfrak H} \times {\mathfrak H}$.
By Sect.~\ref{subsec:realization:odd}, ${\mathcal K}_{\varrho}|_{({\mathfrak H}\times{\mathfrak H})\setminus{\mathfrak D}_{\varrho}}$ is endowed with the involution $\iota_{\varrho}$, which 
does not extend to an 
involution on ${\mathcal K}_{\varrho}$ by Remark~\ref{rmk:invol:not:extend}.

As in  Sect.~\ref{sec:L}, we take a normalized
marking $\alpha_{\varrho}\colon R^{2}\pi_{*}{\bf Z}\cong{\mathbb L}_{K3}$ for $({\mathcal K}_{\varrho}, \iota_{\varrho})$.
Using the map  $\phi$ in \eqref{eqn:phi}, we set 
$$
{\mathbb K}_{\varrho} 
 \colonequals  
\alpha_{\varrho}\left( {\mathbf K} \right) \subset \LAM,
\qquad
\EE_{\varrho} \colonequals \KK_{\varrho}^{\perp_{\LAM}}.
$$

For $\varrho\in{\mathfrak S}_{3}$, let 
\begin{equation}
\label{eqn:varpi:rho}
\varpi_{\varrho}\colon{\mathfrak H}\times{\mathfrak H}\to \Omega_{\KK_{\varrho}}^{+} \subset \Omega_{\LAM}^{+}
\end{equation}
be the period map \eqref{eqn:period} for the marked family 
$(\pi_{\varrho}\colon{\mathcal K}_{\varrho}\to{\mathfrak H}\times{\mathfrak H},\alpha_{\varrho})$. 
Since $\varpi_{\varrho}$ is given explicitly by \eqref{eqn:(8.3)} under the identification \eqref{eqn:identification:lattices}, 
$\varpi_{\varrho}$
induces an isomorphism between ${\mathfrak H}\times{\mathfrak H}$ and 
$\Omega_{{\KK_{\varrho}}}^{+}$. 
Since
$\varpi_{\varrho}(\tau,\tau')=\varpi(K_{\tau,\tau'}/\iota_{\varrho},\alpha_{\varrho})\in\Omega_{\LAM}^{+}\setminus{\mathcal H}$,
we get 
$$
\varpi_{\varrho}^{-1}({\mathcal H}\cap\Omega_{{\mathbb K}_{\varrho}}^{+})\subset{\mathfrak D}_{\varrho}.
$$
\par
Let $\gamma \in A_{\mathbf K} \setminus\{ 0 \}$ be the patching element of $\iota_{\varrho}$. 
Let $\ell(\varrho)$ be the level of $\Gamma_{34}^{\lor}$ in $H^{2}(K_{\tau, \tau'}, {\mathbf Z})_{-}$. 
By Lemma~\ref{lemma:automorphy:Phi:gamma},
\begin{equation}
\label{eqn:Phi:gamma:varrho}
\Phi_{\gamma}(\tau, \tau')^{2} 
= 
\Phi_{\ell(\varrho)}\left(  j_{\ell(\varrho)}^{-1}(\varpi_{\varrho}(\tau, \tau') )\right)^{2}
=
\Phi_{\ell(\varrho)}\left(  j_{\ell(\varrho)}^{-1}\varpi(K_{\tau,\tau'}/\iota_{\varrho},\alpha_{\varrho}) \right)^{2}
\end{equation}
is an automorphic form for $\Gamma(2)\times\Gamma(2)$ of weight $8$.
In particular, $\| \Phi_{\gamma} \|^{2}$ is a $\Gamma(2)\times\Gamma(2)$-invariant $C^{\infty}$ function
on ${\frak H}\times{\frak H}$. Indeed, we have the following equality of functions on ${\mathfrak H}\times{\mathfrak H}$
for all $\varrho \in {\mathfrak S}_{3}$: 
\begin{equation}
\label{eqn:norm:Phi:gamma}
\| \Phi_{\gamma}(\tau,\tau') \|^{2}
=
\| \Phi( {\rm Km}\left( E(\lambda(\tau)) \times E(\varrho(\lambda(\tau')))/\iota \right) \|^{2}.
\end{equation}
In this section, we study the Petersson norm $\| \Phi_{\gamma} \|^{2}$. For this sake, we study the behavior of the period mapping at the boundary.

\subsection{Behavior of the period map at the boundary}
\label{subsec:KM:bdy}
For $n \in \ZZ_{> 0}$, we set
$Y(n) := \Gamma(n)\backslash{\mathfrak H}$,
$X(n) := (\Gamma(n)\backslash{\mathfrak H})^{*}$ and
$$
B_{X(2)\times X(2)} := (X(2)\times X(2) )\setminus (Y(2)\times Y(2)), 
$$
where the asterisk $*$ denotes the Baily--Borel compactification. 
It is classical that the modular $\lambda$-function induces an ${\mathfrak S}_{3}$-equivariant isomorphism from $Y(2)$ to ${\mathbf C}\setminus\{0,1\}$,
where the ${\mathfrak S}_{3}$-action on $Y(2)$ is given by the isomorphism ${\mathfrak S}_{3}\cong \Gamma(1)/\Gamma(2)$. 
Since $X(2)$ has three $0$-dimensional cusps, 
$B_{X(2)\times X(2)}$ is the union of $9$ $\PP^1$'s. 

Let $p\colon Y(2)\to Y(1)$ be the projection. 
Then $p\colon Y(2)\to Y(1)$ is a Galois covering with Galois group ${\mathfrak S}_{3}\cong\Gamma(1)/\Gamma(2)$, 
which induces a ${\mathfrak S}_{3}\times{\mathfrak S}_{3}$-action on $Y(2)\times Y(2)$. 
Let $\varDelta_{Y(1)}$ (resp. $\varDelta_{Y(2)}$) be the diagonal locus of $Y(1)\times Y(1)$ 
(resp. $Y(2)\times Y(2)$) and define the divisor $Z\subset Y(2)\times Y(2)$ as
$$
Z \colonequals (p\times p)^{-1}(\varDelta_{Y(1)})=\sum_{\varrho\in{\mathfrak S}_{3}}(1\times\varrho)\varDelta_{Y(2)}.
$$
Under the identification $Y(2)\cong{\bf C}\setminus\{0,1\}$ via the $\lambda$-invariant, 
we have
$$
Z_{\varrho} = (1\times\varrho)\varDelta_{Y(2)}
$$
for $\varrho\in{\mathfrak S}_{3}$.
Since $\varDelta_{Y(2)}={\mathfrak D}/\Gamma(2)\times\Gamma(2)$ (cf.~\eqref{eqn:frakD}), we have
$Z_{\varrho}={\mathfrak D}_{\varrho}/\Gamma(2)\times\Gamma(2)$ for $\varrho\in{\mathfrak S}_{3}$ (cf.~\eqref{eqn:frakD:gamma}),  and thus 
$Z=\bigcup_{\varrho\in{\mathfrak S}_{3}}Z_{\varrho}$.

By Sect.~\ref{subsec:actions:period:map}, the period map 
$\varpi_{\varrho}$ 
is a $\Gamma(2)\times\Gamma(2)$-equivariant holomorphic map from ${\mathfrak H}\times{\mathfrak H}$ to $\Omega_{\LAM}^{+}$. 
Thus $\varpi_{\varrho}$ descends to a morphism of modular varieties $\overline{\varpi}_{\varrho}\colon Y(2)\times Y(2)\to{\mathcal M}$ such that
$\overline{\varpi}_{\varrho}(  \lambda(\tau), \lambda(\tau') )= \overline{\varpi}(K_{\tau,\tau'}/\iota_{\varrho})$.
Since $\overline{\varpi}_{\varrho}(Y(2)\times Y(2)\setminus Z_{\varrho})\subset{\mathcal M}\setminus{\mathcal D}$, 
it follows from Borel's extension theorem \cite{Borel72} that 
$\overline{\varpi}_{\varrho}$ extends to a holomorphic 
map from $X(2)\times X(2)$ to ${\mathcal M}^*$, which is again denoted by $\overline{\varpi}_{\varrho}$.
Recall from Sect.~\ref{sect:moduli:Enriques:Baily:Borel} that $\Mcal^*\setminus \Mcal$ consists of two modular curves $X(1)$ and $X^1(2)$.

\begin{lemma}
\label{lemma:boundary:KM}
The following inclusion holds:
$$
\overline{\varpi}_{\varrho}(B_{X(2)\times X(2)}) \subset X^1(2).
$$
\end{lemma}

\begin{pf}
Via \eqref{eqn:(8.3)}, $\overline{\varpi}_{\varrho}$ sends the boundary $B_{X(2)\times X(2)}$ to the boundary $\Mcal^*\setminus{\mathcal M}$.  
Assume that there is an irreducible component $C$ of $B_{X(2)\times X(2)}$ with $\overline{\varpi}_{\varrho}(C)\subset X(1)$. 
Since $X(1)={\Dcal}^*\setminus{\Dcal}$ by Lemma~\ref{lemma:boundary:component}, 
we get $\overline{\varpi}_{\varrho}(C)\subset{\Dcal}^*\setminus{\Dcal}$. 
In lattice-theoretical terms, by  \cite[Sect.~2]{Scattone87}, 
an irreducible component of $B_{X(2)\times X(2)}$ corresponds to a primitive totally isotropic sublattice of ${\mathbb K}_{\varrho}$ of rank $2$. 
Thus there is a primitive totally isotropic sublattice $L\subset{\mathbb K}_{\varrho}$ of rank $2$ and a root $\d \in \Delta_{\LAM}$ such that $L \subset \d^{\perp_\LAM}$. 
By Lemma~\ref{lemma:upper}, 
we have $\LAM = \ZZ \delta +{\mathbb K}_{\varrho}\oplus{\mathbb E}_{\varrho}$. 
We write $\delta = \delta_1 + \delta_2$ with $\delta_{1} \in{\mathbb K}_{\varrho}^\vee$ and $\delta_{2} \in{\mathbb E}_{\varrho}^\vee$. 
Since $\iota_{\varrho}$ is of odd type, we get $\delta_1^2 \equiv 1\mod 2$, 
so $\d_1^2 \neq 0$.  
Thus $\d^{\perp_{\LAM}} \cap{\mathbb K}_{\varrho}=\delta_{1}^{\perp_{\KK_{\varrho}}}\cap{\mathbb K}_{\varrho}$ has signature $(1,2)$ or $(2,1)$,
which cannot contain a totally isotropic sublattice of rank $2$.
This contradicts the assumption $L \subset \d^{\perp_\LAM}\cap{\mathbb K}_{\varrho}$. 
Thus $\overline{\varpi}_{\varrho}(C)\subset X^{1}(2)$ for any component $C$ of $B_{X(2)\times X(2)}$.
\end{pf}

\subsection{Involution of odd type and the Borcherds $\Phi$-function}
\label{sect:invol:odd:Borcherds}

Let $\Sigma^{2}Y(1)\colonequals Y(1)\times Y(1)/{\mathfrak S}_{2}$ (resp. $\Sigma^{2}X(1) \colonequals X(1)\times X(1)/{\mathfrak S}_{2}$) denote 
the second symmetric product of $Y(1)$ (resp. $X(1)$), where ${\mathfrak S}_{2}$ acts on $Y(1)\times Y(1)$ (resp. $X(1)\times X(1)$) as the permutation of coordinates.
Since $Y(1)\cong{\bf C}$ and $X(1)\cong{\bf P}^{1}$, we have isomorphisms $\Sigma^{2}Y(1)\cong{\bf C}^{2}$ and $\Sigma^{2}X(1)\cong{\bf P}^{2}$.
Let $\varDelta\subset\Sigma^{2}X(1)$ be the image of the diagonal locus $\varDelta_{X(1)}$ by the projection $X(1)\times X(1)\to\Sigma^{2}X(1)$, i.e.,
\[
\varDelta \colonequals \varDelta_{X(1)}/{\mathfrak S}_{2}\cong{\bf P}^{1}.
\]

We set
\[
B \colonequals \Sigma^{2}X(1)\setminus\Sigma^{2}Y(1).
\]
Under the identifications $X(1)\cong{\bf P}^{1}$ and $Y(1)\cong{\bf P}^{1}\setminus\{\infty\}$,
$B$ is the divisor at infinity of ${\bf P}^{2}$:
\[
B=(\{\infty\}\times X(1)\amalg X(1)\times\{\infty\})/{\mathfrak S}_{2}
\cong
X(1)
\cong
{\bf P}^{1}.
\]

In the rest of this section, we are going to compute $-dd^{c}\left[\log\prod_{\varrho\in{\mathfrak S}_{3}}\overline{\varpi}_{\varrho}^{*}\|\Phi\|^{2}\right]$ 
as a current on $\Sigma^2 X(1)$. 
Since $\overline{\varpi}_{\varrho}( Y(2)\times Y(2)\setminus Z_{\varrho} ) \subset {\mathcal M}\setminus{\mathcal D}$ and ${\rm div}(\Phi) = {\mathcal H}$, 
$\prod_{\varrho\in{\mathfrak S}_{3}}\overline{\varpi}_{\varrho}^{*}\|\Phi\|^{2}$ is a nowhere vanishing $C^{\infty}$ function on $(Y(2)\times Y(2))\setminus Z$.
By \eqref{eqn:S3:invariance}, \eqref{eqn:def:iota:rho}, \eqref{eqn:norm:Phi:gamma}, $\prod_{\varrho\in{\mathfrak S}_{3}}\overline{\varpi}_{\varrho}^{*}\|\Phi\|^{2}$ 
is invariant under the actions of ${\mathfrak S}_{3}\times{\mathfrak S}_{3}$ and ${\mathfrak S}_{2}$. 
Thus we regard $\prod_{\varrho\in{\mathfrak S}_{3}}\overline{\varpi}_{\varrho}^{*}\|\Phi\|^{2}$ as a function on $\Sigma^{2}Y(1)\setminus\varDelta$.

Let $\omega_{{\mathfrak H}\times{\mathfrak H}}$ be the K{\"a}hler form of the Poincar{\'e} metric on ${\mathfrak H}\times{\mathfrak H}$,
i.e.,
\[
\omega_{{\mathfrak H}\times{\mathfrak H}}
 \colonequals 
-dd^{c}\log{\rm Im}\,\tau-dd^{c}\log{\rm Im}\,\tau'.
\]
Let $\omega_{\Sigma^{2}Y(1)}$ (resp. $\omega_{\Sigma^{2}Y(2)}$) be the K{\"a}hler form on $\Sigma^{2}Y(1)$ 
(resp. $\Sigma^{2}Y(2)$) in the sense of orbifolds induced from $\omega_{{\mathfrak H}\times{\mathfrak H}}$.
Since the area of $\Sigma^{2}Y(1)$ (resp. $\Sigma^{2}Y(2)$) with respect to $\omega_{\Sigma^{2}Y(1)}$ 
(resp. $\omega_{\Sigma^{2}Y(2)}$) is finite, 
the $(1,1)$-from $\omega_{\Sigma^{2}Y(1)}$ (resp. $\omega_{\Sigma^{2}Y(2)}$)
extends trivially to a closed positive $(1,1)$-current $\widetilde{\omega_{\Sigma^{2}Y(1)}}$ on $\Sigma^{2}X(1)$. 

\begin{proposition}
\label{prop:a:b}
The following equation of currents on $\Sigma^{2}X(1)$ holds\textup{:}
\[
-dd^{c}
[\log\prod_{\varrho\in{\mathfrak S}_{3}}\overline{\varpi}_{\varrho}^{*}\|\Phi\|^{2}]
=
24\,\widetilde{\omega_{\Sigma^{2}Y(1)}}-\delta_{\varDelta}.
\]
\end{proposition}

\begin{pf}
{\sl Step 1.}\;
Set $F := \log\prod_{\varrho\in{\mathfrak S}_{3}}\overline{\varpi}_{\varrho}^{*}\|\Phi\|^{2} = \sum_{\varrho\in{\mathfrak S}_{3}}\overline{\varpi}_{\varrho}^{*}\log\|\Phi\|^{2}$. 
By \eqref{eqn:def:Petersson:2}, 
we have the equality of $(1,1)$-forms on $\Sigma^{2}Y(1)$
\begin{equation}
\label{eqn:curevature:norm:Phi:2}
-dd^{c} F
=
24\,\omega_{\Sigma^{2}Y(1)}.
\end{equation}
Since $\log\|\Phi\|^{2}$ has logarithmic singularities along $\mathcal D\cup X(1)$,
$F$ has at most logarithmic singularities along $\varDelta\cup B$, 
which, together with \eqref{eqn:curevature:norm:Phi:2} and the irreducibility of $\Delta$ and $B$, 
yields the following equation of currents on $\Sigma^{2}X(1)$
\begin{equation}
\label{eqn:curevature:norm:Phi:3}
-dd^{c} F
=
24\,\widetilde{\omega_{\Sigma^{2}Y(1)}}+\alpha\,\delta_{\varDelta}+\beta\,\delta_{B}, 
\end{equation}
where $\alpha,\beta\in{\bf R}$ are some constants. We are going to verify that $\alpha=-1$ and $\beta=0$. 

{\sl Step 2.}\;
Let $\varpi_{\varrho}^{*}({\mathcal H})$ denote the pull-back of the divisor ${\mathcal H}$. 
Since ${\mathcal H}\cap\varpi_{\varrho}(({\mathfrak H}\times{\mathfrak H})\setminus{\mathfrak D}_{\varrho})=\emptyset$, 
we have ${\rm Supp}(\varpi_{\varrho}^{*}({\mathcal H}))\subset \Dfrak_{\varrho}$.
Since the map $\varpi_{\varrho}$ is $\Gamma(2)\times\Gamma(2)$-equivariant and since the divisor
$Z_{\varrho}={\mathfrak D}_{\varrho}/(\Gamma(2)\times\Gamma(2))$ of $Y(2)\times Y(2)$ is irreducible, the inclusion
${\rm Supp}(\varpi_{\varrho}^{*}({\mathcal H}))\subset \Dfrak_{\varrho}$ implies the existence of an integer
$\nu\in{\bf Z}_{\geq0}$ with
\begin{equation}
\label{eqn:pull:back:discriminant:odd:1}
\varpi_{\varrho}^{*}({\mathcal H})=\nu\,  \Dfrak_{\varrho}.
\end{equation}
Let $i \colon \Omega_{\KK_{\varrho}}^{+} \hookrightarrow \Omega_{\LAM}^{+}$ be the inclusion.
By Proposition~\ref{prop:characterization:zero:Kondo-Mukai}, 
$E \colonequals \sum_{d\in\Delta_{\Lambda}/\{\pm1,\pm I_{\mu}\}}H_d\cap\Omega_{{\mathbb K}_{\varrho}}^{+}$ is a {\em reduced} divisor 
on $\Omega_{{\mathbb K}_{\varrho}}^{+}$ with
\[
i^{*}{\mathcal H}
=
\sum_{d\in\Delta_{\Lambda}/\{\pm1,\pm I_{\mu}\}}
(H_d\cap\Omega_{{\mathbb K}_{\varrho}}^{+}+H_{I_{\mu}(d)}\cap\Omega_{{\mathbb K}_{\varrho}}^{+})
=2E.
\]
Since $\varpi_{\varrho}$ is an isomorphism from ${\mathfrak H}\times{\mathfrak H}$ to $\Omega_{{\mathbb K}_{\varrho}}^{+}$, 
we obtain the equality of divisors
\begin{equation}
\label{eqn:pull:back:discriminant:odd:2}
\varpi_{\varrho}^{*}({\mathcal H})=2\varpi_{\varrho}^{*}(E),
\end{equation}
where $\varpi_{\varrho}^{*}(E)$ is a reduced divisor on ${\mathfrak H}\times{\mathfrak H}$.
Comparing \eqref{eqn:pull:back:discriminant:odd:1} and \eqref{eqn:pull:back:discriminant:odd:2},
we get the equality of divisors on ${\mathfrak H}\times{\mathfrak H}$
\begin{equation}
\label{eqn:pull:back:discriminant:odd:3}
\varpi_{\varrho}^{*}({\mathcal H})=2\,  \Dfrak_{\varrho}.
\end{equation}
Since ${\rm div}(\Phi)={\mathcal H}$, we deduce from \eqref{eqn:pull:back:discriminant:odd:3} that
$\varpi_{\varrho}^{*}\|\Phi\|$ has zeros of order $2$ along ${\mathfrak D}_{\varrho}$. 
Hence $F$ 
has zeros of order $2$ 
along ${\mathfrak D}=\sum_{\varrho\in{\mathfrak S}_{3}}{\mathfrak D}_{\varrho}$. 
This, together with \eqref{eqn:curevature:norm:Phi:2}, implies the following equation of currents on ${\mathfrak H}\times{\mathfrak H}$:
\begin{equation}
\label{eqn:curevature:norm:Phi:4}
-dd^{c} F
=
24\,\omega_{{\mathfrak H}\times{\mathfrak H}}-2\,\delta_{\mathfrak D}.
\end{equation}
Since the projection $p\colon Y(1)\times Y(1)\to\Sigma^{2}Y(1)$ is doubly ramified along $\varDelta_{Y(1)}$, 
we get $\alpha=-2/\deg p=-1$ by \eqref{eqn:curevature:norm:Phi:4}.

{\sl Step 3.}\;
By Lemma~\ref{lemma:boundary:KM}, 
$\overline{\varpi}_{\varrho}(B_{X(2)\times X(2)})$ is not contained in the boundary component $X(1) = {\Dcal}^*\setminus {\Dcal}$. 
Thus the pullback $F$ 
does not vanish identically on $B$. 
Since $B$ is an irreducible divisor on $\Sigma^{2}X(1)$, we get $\beta = 0$. 
\end{pf}

\section{Involutions of odd type: the leading term of $\Phi$}
\label{sec:comparison}

Let $(\pi \colon {\mathcal K} \to {\frak H}\times{\frak H}, \alpha)$ be the marked family of Kummer surfaces defined in Sect.~\ref{sec:KM} 
such that $\alpha$ is normalized for $({\mathcal K}, \iota)$.
So far, we have fixed the Enriques lattice $\LAM$ in $\LL$, and we have considered sublattices $\alpha({\mathbf K})$ 
of $\LAM$ associated to the involutions $\iota$ on $K_{\tau, \tau'}$.
Here $\alpha({\mathbf K})$ is isometric to  $\UU(2) \oplus \UU(2)$. 
In this section, we compute the leading term of $\Phi_{\gamma}$ for odd $\gamma$. 
To do this, 
however, 
descriptions become much simpler if we fix the lattice $\UU(2) \oplus \UU(2) \oplus {\mathbb E}_{8}(2)$ and vary the lattice
$H^{2}(K_{\tau,\tau'}, {\mathbf Z})_{-}$ inside its dual lattice. 

In Sect.~\ref{sect:set:up}, we set the stage.  Then 
we compute the the leading term of $\Phi_{\gamma}(\tau, \tau')$ near $(\tau, \tau') = (+i\infty, +i\infty)$
for odd $\gamma$. 
To ease the notation, we write 
\[
  \KK \colonequals \UU(2) \oplus \UU(2), 
\]
and we set 
\[
{\bf v} \colonequals (1,0,0,0)\in{{\KK}}.
\]
Recall that ${\bf K}$ is the sublattice of $H^2(K_{\tau, \tau'}, \ZZ)$ defined in \eqref{eqn:def:K}. 
We define the isometry $\psi\colon {\bf K} \to \KK$ as 
\begin{align*}
\label{eqn:identify:K:K}
& \psi(-\Gamma_{34}^\vee) = (1,0,0,0), \quad  \psi(\Gamma_{23}^\vee)  = (0,0,1,0), \\
\notag
& \psi(\Gamma_{12}^\vee) = (0,1,0,0), \quad  \psi( -\Gamma_{14}^\vee) = (0,0,0,1). 
\end{align*}
In what follows, we identify ${\mathbf K}$ with $\KK$ via $\psi$.

\subsection{Set up}
\label{sect:set:up}

Let $\gamma\in A_{{\KK}}\setminus\{0\}$ and $\delta_{\gamma}\in A_{{\mathbb E}_{8}(2)}\setminus\{0\}$ 
be such that $\gamma^{2} =  \delta_{\gamma}^{2}$ in $\ZZ/2\ZZ$. 
Let $d_{1}\in{\KK}^{\lor}$ and $d_{2}\in{\mathbb E}_{8}(2)^{\lor}$ be vectors such that
$\bar{d}_{1} = \gamma$ and $\bar{d}_{2} = \delta_{\gamma}$. 
As in Sect.~\ref{sect:period:disc:locus}, the anti-invariant sublattice with respect to 
the fixed-point-free involution corresponding to $\gamma$ is realized as
\begin{equation}
\label{eqn:def:LAM:gamma}
\LAM_{\gamma}={\bf Z}(d_{1},d_{2})+{\KK}\oplus{\mathbb E}_{8}(2).
\end{equation}

Let $(\tau, \tau') \in {\mathfrak H}\times{\mathfrak H}\setminus {\mathfrak D}$. 
Let $\iota\colon K_{\tau,\tau'} \to K_{\tau,\tau'}$ be an involution 
satisfying \eqref{eqn:assumption:lattice} with patching element $\psi^{-1}(\gamma) \in A_{{\mathbf K}}\setminus\{0\}$. 

\begin{lemma}
\label{lem:ext:psi}
There exists an isometry $\widetilde{\psi}\colon H^2(K_{\tau, \tau'}, \ZZ)_{-} \to \LAM_{\gamma}$ extending $\psi$. 
\end{lemma}

\begin{pf}
Let $\{ e_{1},\ldots, e_{8} \}$ (resp. $\{ e'_{1},\ldots, e'_{8} \}$) be a basis of ${\mathbf K}^{\perp}={\mathbf K}^{\perp_{H^{2}_{-}}}$ (resp. $\EE_{8}(2)$) 
whose Gram matrix is the (negative-definite) Cartan matrix of type $E_{8}$. We define an isometry $\psi' \colon {\mathbf K}^{\perp} \to \EE_{8}(2)$
as $\psi'(e_{i}) := e'_{i}$ $(i=1,\ldots,8)$ and we set $\LAM' := (\psi\oplus\psi')( H^2(K_{\tau, \tau'}, \ZZ)_{-} )$.
Since $\LAM' \cong \LAM$, we can express $\LAM' = {\mathbf Z}(d'_{1}, d'_{2}) + \KK \oplus \EE_{8}(2)$ with $d_{1}\in \KK^{\lor}$, $d_{2}\in \EE_{8}(2)^{\lor}$.
Since $\psi^{-1}(\overline{\gamma}) \in A_{{\mathbf K}}\setminus\{0\}$ is the patching element of $\iota$, we have 
$\overline{d'_{1}} = \overline{\gamma} \in A_{\KK}\setminus\{0\}$. If $a\in O(\EE_{8}(2))$ and if we define $\psi'_{a} \colon {\mathbf K}^{\perp} \cong \EE_{8}(2)$
as  $\psi'_{a}( e_{i} ) := a( e'_{i} )$ $(i=1,\ldots,8)$ and set $\LAM'_{a} := (\psi\oplus\psi'_{a})( H^2(K_{\tau, \tau'}, \ZZ)_{-} )$, then we have
$\LAM'_{a} = {\mathbf Z}(\overline{\gamma}, \overline{a}(\overline{d'_{2}})) + \KK \oplus \EE_{8}(2)$. Since the natural homomorphism 
$O(\EE_{8}(2)) \to O(q_{\EE_{8}(2)})$ is surjective by e.g. \cite[l.-7]{BP83}, 
we have $\overline{a}(\overline{d'_{2}}) = \delta_{\gamma}$ by choosing $a \in O(\EE_{8}(2))$ suitably. 
Then $\widetilde{\psi} := \psi \oplus \psi'_{a}$ is the desired isometry. 
\end{pf}

Let $\ell \in \{1, 2\}$ be the level of $\Gamma_{34}^\vee$ in $H^2(K_{\tau, \tau'}, \ZZ)_-$. 
Via $\psi$, $\ell$ is then the level of ${\bf v}$ in $\LAM_{\gamma}$. 
Since the $O(\Lambda)$-orbit of a primitive isotropic vector of $\Lambda_{\gamma}$ is determined by its level, 
we can take 
a primitive isotropic vector ${\bf v'}  \in {\LAM}_{\gamma}$ of level $\ell$ with 
$$
\langle{\bf v},{\bf v}'\rangle= \ell.
$$
We set $\UU(\ell)_{\gamma} := {\mathbf Z} {\mathbf v} + {\mathbf Z} {\mathbf v}'$ and 
$$
{\mathbb M}_{\gamma}
\colonequals 
{\bf v}^{\perp}\cap{{\bf v}'}^{\perp}\cap{\LAM}_{\gamma} 
= \UU(\ell)_{\gamma}^{\perp}\cap\LAM_{\gamma}.
$$
Then $\LAM_{\gamma}=\UU(\ell)_{\gamma} \oplus{\mathbb M}_{\gamma}$, where  
$\UU(\ell)_{\gamma} \cong \UU(\ell)$ and ${\mathbb M}_{\gamma} \cong {\mathbb M}_{\ell}$ (cf. \eqref{eqn:def:M:ell}). 
Recall that $\UU(1), \UU(2) \subset \LAM$ are endowed with the basis $\{ {\mathbf e}_{1}, {\mathbf f}_{1}\}$, $\{ {\mathbf e}_{2}, {\mathbf f}_{2}\}$, respectively.

\begin{lemma}
\label{lemma:noralization:lattice}
Let $\rho, \rho'$ be primitive isotropic vectors of ${\mathbb M}_{\gamma}$ 
such that 
$$
\langle \rho, \rho' \rangle =2/\ell,
\qquad
(-1)^{2/\ell}\langle \rho, ((0,0,1,0),{\mathbf 0})>0,
\qquad
(-1)^{2/\ell}\langle \rho, ((0,0,0,1),{\mathbf 0})>0.
$$
Then there exists a normalized marking $\alpha\colon H^2(K_{\tau, \tau'}, \ZZ) \to {\mathbb L}_{K3}$ for $(K_{\tau, \tau'}, \iota)$ 
such that, if we define the isometry $\beta_\gamma\colon \LAM_\gamma \to \LAM$ as  $\beta_\gamma \colonequals \alpha \circ \widetilde{\psi}^{-1}$, then  
\begin{equation}
\label{eqn:normalization:lattice:2}
\beta_{\gamma}({\mathbf v}) = {\ebf}_{\ell}, \quad 
\beta_{\gamma}({\mathbf v}') = {\fbf}_{\ell}, \quad 
\beta_{\gamma}(\rho) = {\ebf}_{2/\ell}, \quad 
\beta_{\gamma}(\rho') = {\fbf}_{2/\ell}.
\end{equation}
\end{lemma}

\begin{pf}
Let $\alpha'$ be a normalized marking for $(K_{\tau,\tau'}, \iota)$, and we set 
$\beta_\gamma' \colonequals \alpha' \circ \widetilde{\psi}^{-1}$.
Let $g \in O(\LAM)$ be such that $g({\ebf}_{\ell}) = {\ebf}_{\ell}$ and we set 
$\beta_\gamma \colonequals  g\beta_\gamma'$ . Let us see that by choosing $g$ appropriately, $\beta_\gamma$ 
satisfies \eqref{eqn:normalization:lattice:2}. 
We set 
\[
{\mathfrak e}_{\ell} \colonequals \beta'_{\gamma}( {\mathbf v}) =  {\ebf}_{\ell}, \quad  
{\mathfrak f}_{\ell} \colonequals \beta_{\gamma}'({\mathbf v}'), \quad
{\mathfrak e}_{2/\ell} \colonequals \beta_{\gamma}'(\rho), \quad 
{\mathfrak f}_{2/\ell} \colonequals \beta_{\gamma}'(\rho').
\] 
Then ${\mathfrak e}_{\ell}, {\mathfrak f}_{\ell}$ are primitive isotropic vectors of $\LAM$ of level $\ell$ and 
${\mathfrak e}_{2/\ell}$, ${\mathfrak f}_{2/\ell}$ are primitive isotropic vectors of $\LAM$ of level $2/\ell$ such that the Gram matrix of
$\{ {\mathfrak e}_{\ell}, {\mathfrak f}_{\ell}, {\mathfrak e}_{2/\ell}, {\mathfrak f}_{2/\ell} \}$ is given by 
$\binom{0\,\ell}{\ell\,0} \oplus \binom{0\, \frac{2}{\ell}}{\frac{2}{\ell}\,0}$.
Set $L_{1} \colonequals {\mathbf Z}{\mathfrak e}_{\ell} + {\mathbf Z}{\mathfrak f}_{\ell} + {\mathbf Z}{\mathfrak e}_{2/\ell} + {\mathbf Z}{\mathfrak f}_{2/\ell}$ and 
$L_{2} \colonequals L_{1}^{\perp_{\LAM}}$.
It is easy to see that $\LAM = L_{1} \oplus L_{2}$. 
Hence $L_{2}$ is an even $2$-elementary lattice of rank $8$ with $\dim_{{\mathbf F}_{2}} A_{L_{2}} = 8$ and $\delta(L_{2}) =0$, 
which implies $L_{2} \cong \EE_{8}(2)$. Let $\theta \colon L_{2} \to \EE_{8}(2)$ be an isometry of lattices.
We define $g \in O(\LAM)$ as 
\[
g(a {\mathfrak e}_{\ell} + b {\mathfrak f}_{\ell} + c {\mathfrak e}_{2/\ell} + d {\mathfrak f}_{2/\ell} + {\mathbf x} ) 
:= 
a {\ebf}_{\ell} + b {\fbf}_{\ell} + c {\ebf}_{2/\ell} + d {\fbf}_{2/\ell} + \theta( {\mathbf x} ).
\] 
Then 
$\beta_{\gamma} = g \beta_{\gamma}'$ satisfies \eqref{eqn:normalization:lattice:2}. 
\par
By \cite[Cor.\,2.6]{Namikawa85}, there exists $\widetilde{g}\in O({\mathbb L}_{K3})$ with $\widetilde{g}(\LAM)=\LAM$ such that
$\widetilde{g}|_{\LAM}=g$. Set $\alpha:=\widetilde{g}\alpha'$. Then $\alpha(H^{2}(K_{\tau,\tau'},{\mathbf Z})_{-})=\LAM$. 
Since $\alpha\circ\widetilde{\psi}^{-1} = \widetilde{g}\alpha'\circ\widetilde{\psi}^{-1}=g\beta'_{\gamma}=\beta_{\gamma}$ 
and since $\beta_{\gamma}$ satisfies \eqref{eqn:normalization:lattice:2}, we get 
$\alpha(-\Gamma_{34}^{\lor})=\beta_{\gamma}\circ\widetilde{\psi}(-\Gamma_{34}^{\lor})=\beta_{\gamma}({\mathbf v})={\mathbf e}_{\ell}$.
Let us verify \eqref{eqn:definition:period}. Recall that $u$, $B$, $D$ were defined in \eqref{eqn:formula:period}.
Then
$$
\begin{aligned}
\langle (-1)^{2/\ell}B, {\mathbf e}_{2/\ell}\rangle
&=
\langle (-1)^{2/\ell}\alpha(\Gamma_{23}^{\lor}), {\mathbf e}_{2/\ell}\rangle
=
\langle (-1)^{2/\ell}\beta_{\gamma}\circ\widetilde{\psi}(\Gamma_{23}^{\lor}), {\mathbf e}_{2/\ell}\rangle
\\
&=
\langle (-1)^{2/\ell}\beta_{\gamma}\circ\widetilde{\psi}(\Gamma_{23}^{\lor}), \beta_{\gamma}(\rho)\rangle
=
\langle (-1)^{2/\ell}\widetilde{\psi}(\Gamma_{23}^{\lor}), \rho\rangle
\\
&=
(-1)^{2/\ell}\langle ((0,0,1,0),{\mathbf 0}), \rho\rangle>0.
\end{aligned}
$$
Similarly, 
$\langle (-1)^{2/\ell}D, {\mathbf e}_{2/\ell}\rangle=(-1)^{2/\ell}\langle ((0,0,0,1),{\mathbf 0}), \rho\rangle>0$.
Since ${\mathbf e}_{2/\ell}\in\overline{\mathcal C}_{{\mathbb M}_{\ell}}^{+}$ by \eqref{eqn:cond:basis}, we get
$(-1)^{2/\ell}B\in{\mathcal C}_{{\mathbb M}_{\ell}}^{+}$ and $(-1)^{2/\ell}D\in {\mathcal C}_{{\mathbb M}_{\ell}}^{+}$.
By \eqref{eqn:formula:period}, this implies $\Im u \in {\mathcal C}_{{\mathbb M}_{\ell}}^{+}$.
Since $\varpi(K_{\tau,\tau'},\alpha)=[-(u^{2}/2){\mathbf e}_{\ell}+({\mathbf f}_{\ell}/\ell)+(-1)^{2/\ell}u]$ 
by \eqref{eqn:period:prod:kummer} and since $\Im u \in {\mathcal C}_{{\mathbb M}_{\ell}}^{+}$, 
we see that $\alpha$ satisfies \eqref{eqn:definition:period}. This completes the proof.
\end{pf}

In the situation of Lemma~\ref{lemma:noralization:lattice}, via $\beta_\gamma$, let 
$\Omega_{{\LAM}_{\gamma}}^{+}$ and ${\mathcal C}_{{\mathbb M}_{\gamma}}^{+}$ correspond to 
$\Omega_{{\LAM}}^{+}$ and ${\mathcal C}_{{\mathbb M}_{\ell}}^{+}$, respectively. Then, 
similar to \eqref{eqn:tube:domain:level2}, 
the tube domain ${\mathbb M}_{\gamma}\otimes{\bf R}+i\,{\mathcal C}_{{\mathbb M}_{\gamma}}^{+}$
is identified with $\Omega_{{\LAM}_{\gamma}}^{+}$ through the map
$$
j_{\gamma}\colon
{\mathbb M}_{\gamma}\otimes{\bf R}+i\,{\mathcal C}_{{\mathbb M}_{\gamma}}^{+}
\ni z\to
\left[
-(z^{2}/2){\bf v}+{\bf v}'/\ell + (-1)^{2/\ell} z
\right]
\in\Omega_{{\LAM}_{\gamma}}^{+}.
$$
Let 
\begin{equation}
\label{eqn:h:gamma}
h_{\gamma}\colon
\Omega_{{\KK}}\ni[\eta]\to[(\eta,0)]\in
\Omega_{{\LAM}_{\gamma}}
\end{equation}
be the embedding of domains induced by the inclusion of lattices ${{\KK}}\subset\LAM_{\gamma}$, 
and let $\Omega_{{\KK}}^+$ be one of the two (isomorphic) connected components of 
$\Omega_{{\KK}}$ such that $h_{\gamma}(\Omega_{{\KK}}^+) \subset \Omega_{{\LAM}_{\gamma}}^+$. 
Then we have the following commutative diagram:
$$
\begin{CD}
{\frak H}\times{\frak H} @>{\varpi}>> \Omega_{\LAM}^{+} @>{j_{\ell}^{-1}}>> {\mathbb M}_{\ell}\otimes{\mathbf R} + i\,{\mathcal C}_{{\mathbb M}_{\ell}}^{+}
\\
@V{j}VV  @A{\beta_{\gamma}}AA  @A{\beta_{\gamma}}AA
\\
\Omega_{\KK}^{+}  @>{h_{\gamma}}>>  \Omega_{\LAM_{\gamma}}^{+}  @>{j_{\gamma}^{-1}}>> {\mathbb M}_{\gamma}\otimes{\mathbf R} + i\,{\mathcal C}_{{\mathbb M}_{\gamma}}^{+}
\end{CD}
$$
where $\varpi$ is the period mapping for $(\pi \colon {\mathcal K} \to {\frak H}\times{\frak H}, \alpha)$ given by \eqref{eqn:(8.3)}, 
$j$ is the isomorphism given by $j(\tau,\tau') := \left[ \left( -\frac{\tau\tau'}{2}, \frac{1}{2}, \frac{\tau}{2}, \frac{\tau'}{2} \right) \right]$, 
and 
$\beta_{\gamma} \colon {\mathbb M}_{\gamma}\otimes{\mathbf R} + i\,{\mathcal C}_{{\mathbb M}_{\gamma}}^{+} \to
{\mathbb M}_{\ell}\otimes{\mathbf R} + i\,{\mathcal C}_{{\mathbb M}_{\ell}}^{+}$ is the restriction of the linear map
$\beta_{\gamma} \colon {\mathbb M}_{\gamma}\otimes{\mathbf C} \to {\mathbb M}_{\ell}\otimes{\mathbf C}$. Now, we define a map 
$\varphi_{\gamma} \colon {\frak H} \times {\frak H} \to {\mathbb M}_{\gamma}\otimes{\mathbf R} + i\,{\mathcal C}_{{\mathbb M}_{\gamma}}^{+}$ as 
\begin{equation}
\label{eqn:varphi}
\varphi_{\gamma} := \beta_{\gamma}^{-1}\circ j_{\ell}^{-1}\circ\varpi = j_{\gamma}^{-1} \circ h_{\gamma} \circ j.
\end{equation}
For $(\tau, \tau') \in {\frak H}\times{\frak H}$, in Definition~\ref{def:Phi:gamma}, we defined 
$\Phi_{\gamma}(\tau, \tau')^{2} \colonequals \Phi_{\ell} \left(j_\ell^{-1}(\varpi(K_{\tau, \tau'}/\iota, \alpha))\right)^{2}$. 
>From the above commutative diagram, Lemma~\ref{lemma:normalized:marking:2} and Lemma~\ref{lemma:noralization:lattice}, we have 
\begin{equation}
\label{eqn:Phi:gamma}
\Phi_{\gamma}(\tau, \tau')^{2} = 
\Phi_{\ell}\left( \beta_{\gamma}\circ\varphi_{\gamma}(\tau, \tau') \right)^{2}. 
\end{equation}

\medskip
Recall that the elements of $A_{{\KK}}\setminus\{0\}$ of odd norm are represented by the following vectors of ${\KK}^{\lor}$ of norm $-1$ (see the proof of Lemma~\ref{lemma:BP}):
\begin{gather*}
\left(1/2, -1/2, 0, 0\right), 
\left(1/2, -1/2, 1/2, 0\right), 
\left(1/2, -1/2, 0, 1/2 \right), \\
\left(0, 0, 1/2, -1/2\right), 
\left(1/2, 0, 1/2, -1/2\right), 
\left(0, 1/2, 1/2, -1/2\right). 
\end{gather*}

\begin{lemma}
\label{lemma:list:representative}
Regard ${\bf v} = (1, 0, 0, 0)$ as a primitive isotropic vector of ${\LAM}_{\gamma}$.
\begin{itemize}
\item[(1)]
If $\gamma\equiv(\frac{1}{2},\frac{1}{2},\frac{1}{2},0),
(\frac{1}{2},\frac{1}{2},0,0),
(\frac{1}{2},\frac{1}{2},0,\frac{1}{2}),
(0,\frac{1}{2},\frac{1}{2},\frac{1}{2})\mod{{\KK}}$,
then ${\bf v}$ has level $1$ in ${\LAM}_{\gamma}$.
\item[(2)]
If $\gamma\equiv(0,0,\frac{1}{2},\frac{1}{2}),
(\frac{1}{2},0,\frac{1}{2},\frac{1}{2})\mod{{\KK}}$, 
then ${\bf v}$ has level $2$ in ${\LAM}_{\gamma}$.
\end{itemize}
\end{lemma}

\begin{pf}
In \eqref{eqn:def:LAM:gamma}, we may assume that $d_{1}\in{\KK}^{\lor}$ is one of the $6$ vectors in (1), (2) as above.
Since $\langle{\bf v},{{\KK}}\rangle_{{\KK}}=2{\bf Z}$, 
the generator of $\langle{\bf v},{\LAM}_{\gamma}\rangle_{\LAM_{\gamma}}$ has the same parity as that of $\langle{\bf v},d_1\rangle_{\LAM_{\gamma}}$.
>From this, we can verify the assertion.
\end{pf}

For $\tau,\tau'\in{\mathfrak H}$, we set $p^{1/2} \colonequals e^{\pi i\tau}$, $q^{1/2} \colonequals e^{\pi i\tau'}$.
Let ${\bf Z}\{p^{1/2},q^{1/2}\}$ denote the ring of covergent series of $p^{1/2},q^{1/2}$ with coefficients in $\ZZ$.  
Let ${\mathfrak m}$ be the ideal generated by $p^{1/2},q^{1/2}$: 
\[
{\mathfrak m} \colonequals p^{1/2}{\bf Z}\{p^{1/2},q^{1/2}\}+q^{1/2}{\bf Z}\{p^{1/2},q^{1/2}\}.
\]

\subsection{The leading term of $\Phi_{\gamma}$ for odd $\gamma$: level $1$ case}
\label{subsec:leading:level:1}
The rest of this section is devoted to 
determining the leading term of $\Phi_{\gamma}(\tau, \tau')$ near $(+i \infty, +i \infty)$ 
for each odd $\gamma$ as above in Lemma~\ref{lemma:list:representative}. 
In this subesection, we assume that the level of ${\mathbf v}$ 
in $\LAM_{\gamma}$ is $1$. Thus, modulo $\KK$, $\gamma$ is one of the four vectors in Lemma~\ref{lemma:list:representative}~(1).

\begin{lemma}
\label{lemma:odd:leading:level:1}
If $\gamma\equiv(\frac{1}{2},\frac{1}{2},\frac{1}{2},0),
(\frac{1}{2},\frac{1}{2},0,0),
(\frac{1}{2},\frac{1}{2},0,\frac{1}{2}),
(0,\frac{1}{2},\frac{1}{2},\frac{1}{2})\mod\KK$,
then
$$
\Phi_{\gamma}(\tau,\tau')\equiv 1\mod{\mathfrak m}.
$$
\end{lemma}

\begin{pf}
Since ${\bf v}\in\LAM_{\gamma}$ is an isotropic vector of level $1$ by Lemma~\ref{lemma:list:representative},
the proof is the same as that of \cite[Lemma 7.1]{KawaguchiMukaiYoshikawa18}. 
\end{pf}

\subsection
{The leading term of $\Phi_{\gamma}$ for odd $\gamma$: level $2$ case}
\label{subsec:leading:level:2}

In this subesection, we assume that the level of ${\mathbf v} = (1, 0, 0, 0)$ in $\LAM_{\gamma}$ is $2$, i.e., $\ell =2$.
Thus $\gamma\equiv(0,0,\frac{1}{2},\frac{1}{2})$ or $(\frac{1}{2},0,\frac{1}{2},\frac{1}{2})$ $\mod{\KK}$.
We take a vector ${\bf r}\in{\mathbb E}_{8}(2)^{\lor}$ with ${\bf r}^{2}=-1$ and define 
$$
{\LAM}_{\gamma} \colonequals 
\begin{cases}
\begin{array}{lll}
{\bf Z}((0,0,\frac{1}{2},\frac{1}{2}), {\bf r})+{\KK}\oplus{\mathbb E}_{8}(2)
&{\rm if}
&\gamma\equiv(0,0,\frac{1}{2},\frac{1}{2})\mod\KK,
\\
{\bf Z}((\frac{1}{2},0,\frac{1}{2},\frac{1}{2}),{\bf r})+{\KK}\oplus{\mathbb E}_{8}(2)
&{\rm if}
&\gamma\equiv(\frac{1}{2},0,\frac{1}{2},\frac{1}{2})\mod\KK.
\end{array}
\end{cases}
$$
We set 
$$
{\bf v}' \colonequals 
\begin{cases}
\begin{array}{lll}
((0,1,0,0),0)&{\rm if}
&\gamma\equiv(0,0,\frac{1}{2},\frac{1}{2})\mod\KK,
\\
((0,1,1,0),0)&{\rm if}
&\gamma\equiv(\frac{1}{2},0,\frac{1}{2},\frac{1}{2})\mod\KK.
\end{array}
\end{cases}
$$
Then ${\bf v}'\in {\LAM}_{\gamma}$ is a primitive isotropic vector of $\LAM_{\gamma}$  
with $\langle{\bf v},{\bf v}'\rangle= 2$ and $\langle {\mathbf v}', {\mathbf r}\rangle=0$. 
We define ${\mathbf w}\in{\bf v}^{\perp}\cap{{\bf v}'}^{\perp}\cap{\LAM}_{\gamma}$ and
${\mathbf w}'\in{\bf v}^{\perp}\cap{{\bf v}'}^{\perp}\cap{\LAM}_{\gamma}^{\lor}$ as 
$$
{\mathbf w} \colonequals ((0,0, -1,0),0),
\qquad
{\mathbf w}' \colonequals 
\begin{cases}
\begin{array}{lll}
((0,0, \frac{1}{2}, -\frac{1}{2}),0)&{\rm if}
&\gamma\equiv(0,0,\frac{1}{2},\frac{1}{2})\mod\KK,
\\
((\frac{1}{2},0, \frac{1}{2}, -\frac{1}{2}),0)&{\rm if}
&\gamma\equiv(\frac{1}{2},0,\frac{1}{2},\frac{1}{2})\mod\KK.
\end{array}
\end{cases}
$$
Then ${\mathbf w}^{2}=0$, ${\mathbf w}'^{2}=-1$, $\langle {\mathbf w}, {\mathbf w}'\rangle=1$ 
and $\langle {\mathbf w}, {\mathbf r}\rangle = \langle {\mathbf w}', {\mathbf r}\rangle =0$. 

We are going to express elements of 
${\mathbb M}_{\gamma} = {\bf v}^{\perp}\cap{{\bf v}'}^{\perp}\cap{\LAM}_{\gamma}$ 
concretely. We note that an element of 
${\bf v}^{\perp}\cap{{\bf v}'}^{\perp}$ in ${\KK}^\vee\oplus{\mathbb E}_{8}(2)^\vee$ 
is written as
\[
(a+b){\mathbf w}+2b{\mathbf w}'+(0,x) 
=
\begin{cases}
\begin{array}{lll}
((0,0, -a, -b),x)
&{\rm if}&\gamma\equiv(0,0,\frac{1}{2},\frac{1}{2})\mod\KK,
\\
((b,0, -a, -b),x)
&{\rm if}&\gamma\equiv(\frac{1}{2},0,\frac{1}{2},\frac{1}{2})\mod\KK,
\end{array}
\end{cases}
\]
where $a, b \in (1/2)\ZZ$ and $x \in \EE_8(2)^{\lor}$. 

We define $\rho,\rho'\in{\bf v}^{\perp}\cap{{\bf v}'}^{\perp}\cap{\LAM}_{\gamma}$ as 
\begin{align*}
\rho & \colonequals {\mathbf w}+{\mathbf w}'+(0,{\bf r})
= \begin{cases}
\begin{array}{lll}
((0,0,-\frac{1}{2},-\frac{1}{2}), {\bf r})
&{\rm if}
&\gamma\equiv(0,0,\frac{1}{2},\frac{1}{2})\mod\KK,
\\
((\frac{1}{2},0, -\frac{1}{2}, -\frac{1}{2}),{\bf r})
&{\rm if}
&\gamma\equiv(\frac{1}{2},0,\frac{1}{2},\frac{1}{2})\mod\KK,
\end{array}
\end{cases}
\\
\rho' & \colonequals {\mathbf w}.
\end{align*}
Then $\rho$ and $\rho'$ are primitive isotropic vectors of ${\mathbb M}_{\gamma}$ with level $1$,  
which satisfy the assumptions   $\langle\rho,\rho'\rangle=1$, 
$-\langle \rho, ((0,0,1,0),0)\rangle>0$, and $-\langle \rho, ((0,0,0,1),0)\rangle>0$ in Lemma~\ref{lemma:noralization:lattice}. 
Since 
\begin{center}
$
{\bf v}^{\perp}\cap{{\bf v}'}^{\perp}\cap({\KK}\oplus{\mathbb E}_{8}(2))
=
\{(a+b){\mathbf w}+2b{\mathbf w}'+(0,x)\mid a,b\in{\bf Z},\,x\in{\mathbb E}_{8}(2)\}
$
\end{center}
and ${\LAM}_{\gamma}={\bf Z}\rho+{\KK}\oplus{\mathbb E}_{8}(2)$,
we get by the identification ${\mathbb M}_{\gamma}={\bf v}^{\perp}\cap{{\bf v}'}^{\perp}\cap{\LAM}_{\gamma}$: 
\begin{equation}
\label{eqn:M:gamma}
\begin{aligned}
{\mathbb M}_{\gamma}
&=
{\bf Z}\rho
+
{\bf v}^{\perp}\cap{{\bf v}'}^{\perp}\cap
({\KK}\oplus{\mathbb E}_{8}(2))
\\
&=
\{(a+b)\,{\mathbf w} + 2b\,{\mathbf w}' + (0,x)\mid 
a,b\in(1/2){\bf Z},\,a-b\in{\bf Z},\,
x\in2a {\bf r}+{\mathbb E}_{8}(2)\}.
\end{aligned}
\end{equation}

By Lemma~\ref{lemma:noralization:lattice}, we can take an isometry of lattices $\beta_{\gamma} \colon \LAM_{\gamma} \to \LAM$ 
satisfying \eqref{eqn:normalization:lattice:2} 
for the vectors ${\mathbf v}$, ${\mathbf v}'$, $\rho$, $\rho'$ as above. Then $\varphi_{\gamma}$ is given as follows.

\begin{lemma}
\label{lemma:varphi:gamma}
The following equality holds:
\begin{equation}
\label{eqn:Phi:gamma:t:t'}
\begin{aligned}
\varphi_{\gamma}(\tau,\tau') 
&=
\begin{cases}
\begin{array}{lll}
((0,0, \frac{-\tau}{2}, \frac{-\tau'}{2}),0)&{\rm if}
&\gamma\equiv(0,0,\frac{1}{2},\frac{1}{2})\mod\KK,
\\
((\frac{\tau'}{2},0, \frac{-\tau+1}{2}, \frac{-\tau'}{2}),0)&{\rm if}
&\gamma\equiv(\frac{1}{2},0,\frac{1}{2},\frac{1}{2})\mod\KK,
\end{array}
\end{cases}
\\
&=
\begin{cases}
\begin{array}{lll}
\frac{\tau+\tau'}{2}{\mathbf w} + \tau'\,{\mathbf w}'
&{\rm if}&\gamma\equiv(0,0,\frac{1}{2},\frac{1}{2})\mod\KK,
\\
\frac{\tau+\tau'-1}{2}{\mathbf w} + \tau'\,{\mathbf w}'
&{\rm if}&\gamma\equiv(\frac{1}{2},0,\frac{1}{2},\frac{1}{2})\mod\KK.
\end{array}
\end{cases}
\end{aligned}
\end{equation}
\end{lemma}

\begin{pf}
To compute $\varphi_{\gamma}(\tau,\tau')$, write $\varpi(K_{\tau,\tau'}, \alpha) = j_{\ell}(u)$ with 
$u \in {\mathbb M}_{\ell}\otimes{\mathbf R} + i\,{\mathcal C}_{{\mathbb M}_{\ell}}^{+}$. 
By \eqref{eqn:formula:period}, we can express 
$\beta_{\gamma}^{-1}(u) = -(\widetilde{A} + \widetilde{B}\tau + \widetilde{D}\tau')/2$, where
$$
\begin{aligned}
\widetilde{A} 
&= 
\beta_{\gamma}^{-1}(A) 
= 
\beta_{\gamma}^{-1}\alpha(\Gamma_{12}^{\lor}) 
- 
\langle \beta_{\gamma}^{-1}({\fbf}_{2}/2), \beta_{\gamma}^{-1}\alpha(\Gamma_{12}) \rangle \beta_{\gamma}^{-1}{\ebf}_{2} 
- 
\beta_{\gamma}^{-1}{\fbf}_{2}
\\
&
=
\widetilde{\psi}(\Gamma_{12}^{\lor}) - \langle \beta_{\gamma}^{-1}({\mathbf f}_{2}/2), \widetilde{\psi}(\Gamma_{12})\rangle {\mathbf v}
-{\mathbf v}' 
\\
&=
((0,1,0,0),0) - \langle {\mathbf v}'/2, ((0,1,0,0),0)\rangle {\mathbf v} -{\mathbf v}' 
\\
&=
\begin{cases}
\begin{array}{lll}
((0,0,0,0),0)&{\rm if}
&\gamma\equiv(0,0,\frac{1}{2},\frac{1}{2})\mod\KK,
\\
((0,0,-1,0),0)&{\rm if}
&\gamma\equiv(\frac{1}{2},0,\frac{1}{2},\frac{1}{2})\mod\KK.
\end{array}
\end{cases}
\end{aligned}
$$
Here we used \eqref{eqn:normalization:lattice:2} to get the fourth equality. Similarly, we get
\begin{align*}
\widetilde{B} 
&= 
\beta_{\gamma}^{-1}(B)
=
-\langle \beta_{\gamma}^{-1}({\fbf}_{2}/2), \beta_{\gamma}^{-1}\alpha(\Gamma_{23}^{\lor}) \rangle \beta_{\gamma}^{-1}{\ebf}_{2} 
+ 
\beta_{\gamma}^{-1}\alpha(\Gamma_{23}^{\lor})
\\
&=
-\langle {\mathbf v}'/2, ((0,0,1,0),0) \rangle {\mathbf v} + ((0,0,1,0),0)
=
((0,0,1,0),0),
\\
\widetilde{D} 
&= 
\beta_{\gamma}^{-1}(D)
=
\langle \beta_{\gamma}^{-1}({\fbf}_{2}/2), \beta_{\gamma}^{-1}\alpha(\Gamma_{14}^{\lor}) \rangle \beta_{\gamma}^{-1}{\ebf}_{2} 
- 
\beta_{\gamma}^{-1}\alpha(\Gamma_{14}^{\lor})
\\
&=
\langle {\mathbf v}'/2, ((0,0,0,-1),0) \rangle {\mathbf v} + ((0,0,0,1), 0)
\\
&=
\begin{cases}
\begin{array}{lll}
((0,0,0,1),0)&{\rm if}
&\gamma\equiv(0,0,\frac{1}{2},\frac{1}{2})\mod\KK,
\\
((-1,0,0,1),0)&{\rm if}
&\gamma\equiv(\frac{1}{2},0,\frac{1}{2},\frac{1}{2})\mod\KK.
\end{array}
\end{cases}
\end{align*}
Substituting these formulae into the following equation (cf. \eqref{eqn:varphi})
$$
\varphi_{\gamma}(\tau,\tau') 
= 
\beta_{\gamma}^{-1}\circ j_{\ell}^{-1} \circ \varpi (\tau, \tau') =
\beta_{\gamma}^{-1}(u) = -(\widetilde{A} + \widetilde{B}\tau + \widetilde{D}\tau')/2,
$$
we get the result.
\end{pf}

We recall from  Sect.~\ref{subsubsec:Phi:level:2} and the relation $\beta_{\gamma}^{-1}({\ebf}_{1}) = \rho$ that 
\[
 \Pi^{+}_{\gamma}
 :=
\beta_{\gamma}^{-1}(\Pi^{+}) = \left\{
 \lambda \in \mathbb{M}_\gamma \mid \langle \lambda, \rho \rangle > 0, \, \lambda^2 \geq -2 
 \right\}.
\]
By \eqref{eqn:Phi:gamma}, $\Phi_{\gamma}(\tau,\tau') = \pm\Phi_{\ell}( \beta_{\gamma}\varphi_{\gamma}(\tau,\tau'))$ is given by 
{\allowdisplaybreaks
\begin{align}
\label{eqn:Phi:expansion}
  & \pm \Phi_{\gamma}(\tau, \tau') 
  \\
  \notag
   & \quad =
  2^{8}e^{2\pi i\langle {\ebf}_{1},  \beta_{\gamma}\varphi_{\gamma}(\tau,\tau')\rangle} 
  \prod_{\mu \in \ZZ_{>0}{\ebf}_{1}\,\cup\,\Pi^+} 
  \left(
  1 - e^{2\pi i \langle \mu,  \beta_{\gamma}\varphi_{\gamma}(\tau,\tau') \rangle}
  \right)^{(-1)^{\langle \mu, 
  {\ebf}_{1} - {\fbf}_{1}\rangle} c(\mu^2/2)}
  \\
  \notag
  & \quad =   
  2^{8}e^{2\pi i\langle\r, \varphi_{\gamma}(\tau, \tau')\rangle} 
  \prod_{\lambda \in \ZZ_{>0}\rho\,\cup\,\Pi^+_{\gamma}} 
  \left(
  1 - e^{2\pi i \langle \l, \varphi_{\gamma}(\tau, \tau')\rangle}
  \right)^{(-1)^{\langle \lambda, 
  \r - \r^\prime\rangle} c(\lambda^2/2)}
\end{align}
}
for $\tau, \tau' \in \HH$, where we used 
$\beta_{\gamma}^{-1}({\ebf}_{1})=\rho$, $\beta_{\gamma}^{-1}({\fbf}_{1})=\rho'$ to get the second equality.

We will compute the leading term of $\Phi_{\gamma}(\tau, \tau')$ near 
$(\tau, \tau') = (+i\infty, + i\infty)$ by dividing $\prod_{\lambda \in \ZZ_{>0}\rho\,\cup\,\Pi^+_{\gamma}}$ 
into three parts: $\prod_{\lambda \in \Pi^+_{\gamma}, \, H_\lambda \cap \Omega_{\KK}^+ \neq \emptyset}$, 
$\prod_{\lambda \in \Pi^+_{\gamma}, \, H_\lambda \cap \Omega_{\KK}^+ = \emptyset}$,  
and $\prod_{\lambda \in \ZZ_{>0}\rho}$, 
which will be respectively treated in the following. 

First we consider the part $\prod_{\lambda \in \Pi^+_{\gamma}, \, H_\lambda \cap \Omega_{\KK}^+ \neq \emptyset}$ of \eqref{eqn:Phi:expansion}.

\begin{lemma}
\label{lemma:odd:leading:part:1:pre}
If $\lambda\in\Pi^{+}_{\gamma}$ satisfies $H_\lambda\cap\Omega_{{\KK}}^{+}\not=\emptyset$, then $\lambda=\pm{\mathbf w}'- {\mathbf r}$.
\end{lemma}

\begin{pf}
Let $\lambda\in\Pi^{+}_{\gamma}$ be such that $H_\lambda\cap\Omega_{{\KK}}^{+}\not=\emptyset$.
Then $\lambda^{2}<0$ (see \eqref{eqn:Hd:d2}). Since $\lambda^2 \in 2 \ZZ$ and $\lambda \in \Pi^{+}_{\gamma}$, 
we obtain $\lambda^{2}=-2$, i.e., $\lambda\in \Delta_{{\mathbb M}_{\gamma}} \subset \Delta_{\LAM_\gamma}$. 
We write $\lambda=\lambda_{1}+\lambda_{2}$, where $\lambda_{1}\in{\KK}^{\lor}$, 
$\lambda_{2}\in{\mathbb E}_{8}(2)^{\lor}$. 
As in the proof of Step~2 of Proposition~\ref{prop:characterization:zero:Kondo-Mukai}, 
we have $\lambda_{1}^{2}=\lambda_{2}^{2}=-1$. 
By Lemma~\ref{lemma:upper}, 
${\LAM}_{\gamma}={\bf Z}\lambda+{{\KK}}\oplus{\mathbb E}_{8}(2)$, 
so $\lambda +\rho\in{\KK}\oplus{\mathbb E}_{8}(2)$.
Since $\langle\lambda,{\bf v}\rangle=0$, we can write $\lambda_{1}=(c,0,\frac{a}{2},\frac{b}{2})$ 
with $a\equiv b\equiv 1\mod 2$ and $c\in (1/2){\bf Z}$. 
Since $-1=\lambda_{1}^{2}=ab$, we get $\lambda_{1}=c'\,{\bf v}\pm{\mathbf w}'$ with 
$c' \in (1/2){\bf Z}$.  
Since $\langle\lambda,{\bf v}'\rangle=0$, we get $c'=0$, so $\lambda_{1}=\pm{\mathbf w}'$. 
It follows from $\lambda + \rho\in{\KK}\oplus{\mathbb E}_{8}(2)$ that 
$\lambda_2 + {\mathbf r} \in {\mathbb E}_{8}(2)$. Since $2 {\bf r} \in  {\mathbb E}_{8}(2)$, 
we can write $\lambda_{2}= - {\bf r}+x$ with $x\in{\mathbb E}_{8}(2)$. 
Since $\lambda_{2}^{2}= -1$ and ${\bf r}^{2}=-1$, we get $\langle  {\bf r}, x\rangle=  x^{2}/2$. 
Noting that $\lambda \in \Pi_\gamma^+$, we have 
\[
0<\langle\rho,\lambda\rangle
=
\langle
{\mathbf w} + {\mathbf w}' + {\bf r},  \pm{\mathbf w}' - {\bf r}+x
\rangle
= 
\pm \langle {\mathbf w} + {\mathbf w}' , {\mathbf w}' \rangle 
+ 1 + \langle {\bf r} ,x\rangle
=
1 + x^{2}/2. 
\]
Since $x^{2}/2\in2{\bf Z}_{\leq0}$, 
we have $x = 0$. Thus $\lambda = \pm{\mathbf w}' - {\bf r}$. 
\end{pf}

\begin{lemma}
\label{lemma:odd:leading:part:1}
The following equality holds:
\begin{equation}
\label{eqn:lemma:odd:leading:part:1}
\begin{aligned}
\,&
\prod_{\lambda\in\Pi^{+}_{\gamma},\, H_\lambda\cap\Omega_{{\KK}}^{+}\not=\emptyset}
\left(
1-e^{2\pi i \langle\lambda,\varphi_{\gamma}(\tau,\tau')\rangle}
\right)^{(-1)^{\langle \lambda, 
  \r - \r^\prime\rangle} c(\lambda^2/2)}
\\
&=
\begin{cases}
\begin{array}{lll}
(1-p^{-1/2}q^{1/2})(1-p^{1/2}q^{-1/2}),
&{\rm if}&\gamma\equiv(0,0,\frac{1}{2},\frac{1}{2})\mod\KK,
\\
(1+p^{-1/2}q^{1/2})(1+p^{1/2}q^{-1/2}),
&{\rm if}&\gamma\equiv(\frac{1}{2},0,\frac{1}{2},\frac{1}{2})\mod\KK.
\end{array}
\end{cases}
\end{aligned}
\end{equation}
\end{lemma}

\begin{pf}
We have 
$
\langle \lambda,  \r - \r^\prime\rangle
= 
\langle \pm{\mathbf w}' - {\bf r},  {\mathbf w}' + {\bf r}\rangle
= \pm {{\mathbf w}'}^2 - {\bf r}^2 
= \pm 1 + 1
$ and $c(\lambda^2/2) = c(-1) = 1$, so 
$
(-1)^{\langle \lambda, 
  \r - \r^\prime\rangle} c(\lambda^2/2) = 1
$. 
Since 
$$
\langle\lambda,\varphi_{\gamma}(\tau,\tau')\rangle
=
\begin{cases}
\begin{array}{lll}
\pm(\frac{\tau-\tau'}{2})
&
{\rm if}
&
\lambda=\pm{\mathbf w}' - {\bf r},
\quad
\gamma\equiv(0,0,\frac{1}{2},\frac{1}{2})\mod\KK,
\\
\pm(\frac{\tau-\tau'-1}{2})
&
{\rm if}
&
\lambda=\pm{\mathbf w}' - {\bf r}, 
\quad
\gamma\equiv(\frac{1}{2},0,\frac{1}{2},\frac{1}{2})\mod\KK
\end{array}
\end{cases}
$$
by \eqref{eqn:Phi:gamma:t:t'}, we get the result.  
\end{pf}

Next we consider the part $\prod_{\lambda \in \Pi^+_{\gamma}, \, H_\lambda \cap \Omega_{\KK}^+ = \emptyset}$ of \eqref{eqn:Phi:expansion}.

\begin{lemma}
\label{lemma:odd:leading:part:2:pre}
Let $\lambda\in\Pi^{+}_{\gamma}$ be such that $H_\lambda\cap\Omega_{{\KK}}^{+}=\emptyset$. 
We write $\lambda=\lambda_{1}+\lambda_{2}\in\Pi^{+}_{\gamma}$, 
where $\lambda_{1}=(c,0,a,b)\in{\KK}^{\lor}$, $a,b,c\in (1/2){\bf Z}$, $\lambda_{2}\in{\mathbb E}_{8}(2)^{\lor}$. 
Then $a \leq 0$, $b \leq 0$, and $(a,b)\not=(0,0)$. 
\end{lemma}

\begin{pf}
Since $H_{\lambda_{1}}\cap\Omega_{{\KK}}^{+}=\emptyset$, we have $\lambda_{1}^2 \geq 0$ (see \eqref{eqn:Hd:d2}). It follows that $ab\geq0$. 
Hence ``$a \leq 0$ and $b\leq 0$'' or ``$a \geq 0$ and 
$b\geq 0$.'' 
To derive a contradiction, 
we assume that $a\geq 0$ and $b\geq 0$. 
We set $a' \colonequals -a$, $b' \colonequals -b$. Then $a' \leq 0$ and $b'\leq 0$. 
We set  
\begin{align*}
\nu \colonequals \langle\rho,\lambda\rangle
& = \begin{cases}
\begin{array}{lll}
\langle(0,0,-\frac{1}{2},-\frac{1}{2})+ {\bf r}, (c,0,a,b)+\lambda_{2} \rangle
&{\rm if}
&\gamma\equiv(0,0,\frac{1}{2},\frac{1}{2})\mod\KK,
\\
\langle(\frac{1}{2},0, -\frac{1}{2},-\frac{1}{2})+ {\bf r}, (c,0,a,b)+\lambda_{2} \rangle
&{\rm if}
&\gamma\equiv(\frac{1}{2},0,\frac{1}{2},\frac{1}{2})\mod\KK,
\end{array}
\end{cases}\\
& = 
-a - b + \langle {\bf r}, \lambda_{2}\rangle. 
\end{align*}
We also set $2k \colonequals \lambda^{2}$. 
Since $\lambda \in \Pi^+_{\gamma}$, we have 
$\nu > 0$ and $k \geq -1$. 
Note that $\nu \in \ZZ$ and $k \in \ZZ$. 
Since the sublattice ${\bf Z} {\bf r}+{\bf Z}\lambda_{2}\subset{\mathbb E}_{8}(2)^{\lor}$ is negative-definite, the matrix
$$
\begin{pmatrix}
{\bf r}^{2}&
\langle {\bf r},\lambda_{2}\rangle
\\
\langle {\bf r},\lambda_{2}\rangle&
\lambda_{2}^{2}
\end{pmatrix}
=
\begin{pmatrix}
-1& \nu -a'-b'
\\
\nu -a'-b' &2k-4a'b'
\end{pmatrix}
$$
is semi-negative, so $(\nu -a'-b')^{2}\leq4a'b'-2k$. Thus we get the inequality
\begin{equation}
\label{eqn:lemma:odd:leading:part:2:pre}
(a'-b')^{2}+ \nu^{2} -2\nu(a'+b')\leq-2k.
\end{equation}
Since $a' \leq 0$, $b'\leq 0$ and $\nu >0$, we get $k<0$. 
It follows from $k \in \ZZ$ and $k \geq -1$ that $k=-1$. 
Then, by \eqref{eqn:lemma:odd:leading:part:2:pre}, we conclude that 
$\nu =1$ and $a'=b'=0$. 
Since $a'=b'=0$ and $\lambda\in{\mathbb M}_{\gamma}\subset{{\bf v}'}^{\perp}$, we get $0=\langle\lambda,{\bf v}'\rangle=c$. 
Hence $\lambda_{1}=0$, so $\lambda=\lambda_{2}\in{\mathbb E}_{8}(2)$.
Since $\lambda^{2}=2k=-2$, this contradicts $\Delta_{{\mathbb E}_{8}(2)}=\emptyset$.
Thus we conclude $a \leq 0$, $b\leq 0$, and $(a,b)\not=(0,0)$. 
\end{pf}

\begin{lemma}
\label{lemma:odd:leading:part:2}
The following holds:
\begin{equation}
\label{eqn:lemma:odd:leading:part:2}
\prod_{\lambda\in\Pi^{+}_{\gamma},\,H_\lambda\cap\Omega_{{\KK}}^{+}=\emptyset}
\left(
1-e^{2\pi i \langle\lambda,\varphi_{\gamma}(\tau,\tau')\rangle}
\right)^{(-1)^{\langle \lambda, 
  \r - \r^\prime\rangle} c(\lambda^2/2)}
\in 1+{\mathfrak m}.
\end{equation}
\end{lemma}

\begin{pf}
It suffices to show that for any  $\lambda\in\Pi^{+}_{\gamma}$ with $H_\lambda\cap\Omega_{\KK}^{+}=\emptyset$, 
we have 
\[
1-e^{2\pi i\langle\lambda,\varphi_{\gamma}(\tau,\tau')\rangle_{{\mathbb M}_{\gamma}}} 
\in 1+{\mathfrak m} \subset {\bf Z}\{p^{1/2},q^{1/2}\}. 
\]
We write $\lambda=(c,0,a,b)+\lambda_{2}= c  {\mathbf v}  -(a+b){\mathbf w} - 2b{\mathbf w}' + \lambda_{2}\in{\mathbb M}_{\gamma}$
as in Lemma~\ref{lemma:odd:leading:part:2:pre}. 
By \eqref{eqn:M:gamma}, \eqref{eqn:Phi:gamma:t:t'}, we get 
$$
\langle\lambda,\varphi_{\gamma}(\tau,\tau')\rangle
=
\begin{cases}
\begin{array}{lll}
-b\tau-a\tau'&{\rm if}
&\gamma\equiv(0,0,\frac{1}{2},\frac{1}{2})\mod\KK,
\\
-b\tau-a\tau'+b&{\rm if}
&\gamma\equiv(\frac{1}{2},0,\frac{1}{2},\frac{1}{2})\mod\KK. 
\end{array}
\end{cases}
$$
Since $a \leq 0$ and $b\leq 0$ with $(a, b) \neq (0, 0)$ 
by Lemma~\ref{lemma:odd:leading:part:2:pre}, 
we have $e^{2 \pi i \langle\lambda,\varphi_{\gamma}(\tau,\tau')\rangle} \in {\mathfrak m}$ and 
we get the result. 
\end{pf}

Finally we consider the part $\prod_{\lambda \in \ZZ_{>0}\rho}$ of \eqref{eqn:Phi:expansion}.

\begin{lemma}
\label{lemma:odd:leading:part:3}
The following equality holds:
\begin{equation}
\label{eqn:lemma:odd:leading:part:3}
e^{2\pi i\langle\rho,\varphi_{\gamma}(\tau,\tau')\rangle}
=
\begin{cases}
\begin{array}{lll}
p^{1/2}q^{1/2}&{\rm if}
&\gamma\equiv(0,0,\frac{1}{2},\frac{1}{2})\mod\KK,
\\
-p^{1/2}q^{1/2}&{\rm if}
&\gamma\equiv(\frac{1}{2},0,\frac{1}{2},\frac{1}{2})\mod\KK.
\end{array}
\end{cases}
\end{equation}
\end{lemma}

\begin{pf}
Since
$\langle\rho,\varphi_{\gamma}(\tau,\tau')\rangle_{{\mathbb M}_{\gamma}}=\frac{\tau+\tau'}{2}$ 
if $\gamma\equiv(0,0,\frac{1}{2},\frac{1}{2})$ 
and since
$\langle\rho,\varphi_{\gamma}(\tau,\tau')\rangle_{{\mathbb M}_{\gamma}}=\frac{\tau+\tau'-1}{2}$ 
if $\gamma\equiv(0,0,\frac{1}{2},\frac{1}{2})$ by the definition of $\rho$ and \eqref{eqn:Phi:gamma:t:t'}, 
we get the result.
\end{pf}

\begin{lemma}
\label{lemma:leading:rho}
The following holds:
\begin{align}
\label{eqn:lemma:leading:rho}
& e^{2\pi i\langle\rho,\varphi_{\gamma}(\tau,\tau')\rangle}
\prod_{n>0}
\left(
1-e^{2\pi i\langle n \rho,\varphi_{\gamma}(\tau,\tau')\rangle}
\right)^{(-1)^{\langle n \rho, 
  \r - \r^\prime\rangle} c((n \rho)^2/2)}
\\
\nonumber
& \qquad \in 
\begin{cases}
\begin{array}{lll}
p^{1/2}q^{1/2} (1+{\mathfrak m}) &{\rm if}
&\gamma\equiv(0,0,\frac{1}{2},\frac{1}{2})\mod\KK,
\\
-p^{1/2}q^{1/2}(1+{\mathfrak m}) &{\rm if}
&\gamma\equiv(\frac{1}{2},0,\frac{1}{2},\frac{1}{2})\mod\KK.
\end{array}
\end{cases}
\end{align}
\end{lemma}

\begin{pf}
We have 
$\langle n \rho, \r - \r^\prime\rangle = \langle n \rho, - \r^\prime\rangle = -n$ and 
$c((n \rho)^2/2) = c(0) = 8$ (see Sect.~\ref{subsubsec:c(n)}), so 
$(-1)^{\langle n \rho, 
  \r - \r^\prime\rangle} c((n \rho)^2/2) = 8 \cdot (-1)^n$. 
By  Lemma~\ref{lemma:odd:leading:part:3}, we get 
 $\left(
1-e^{2\pi i\langle n \rho,\varphi_{\gamma}(\tau,\tau')\rangle}
\right)^{8 \cdot (-1)^n} \in 1 + {\mathfrak m}$ and the assertion \eqref{eqn:lemma:leading:rho}. 
\end{pf}

All together, we obtain the leading term of \eqref{eqn:Phi:expansion}.

\begin{lemma}
\label{lemma:odd:leading:level:2}
The following holds:
$$
\Phi_{\gamma}(\tau,\tau')
\in 
\begin{cases}
\begin{array}{lll}
-2^{8}(p^{1/2}-q^{1/2})^{2}
(1 + {\mathfrak m}) 
&{\rm if}
&\gamma\equiv(0,0,\frac{1}{2},\frac{1}{2})\mod\KK,
\\
-2^{8}(p^{1/2}+q^{1/2})^{2}
(1 + {\mathfrak m}) 
&{\rm if}
&\gamma\equiv(\frac{1}{2},0,\frac{1}{2},\frac{1}{2})\mod\KK.
\end{array}
\end{cases}
$$
\end{lemma}

\begin{pf}
Suppose that $\gamma\equiv(0,0,\frac{1}{2},\frac{1}{2})$. Then 
substituting \eqref{eqn:lemma:odd:leading:part:1}, \eqref{eqn:lemma:odd:leading:part:2}, 
\eqref{eqn:lemma:leading:rho} into \eqref{eqn:Phi:expansion}, 
we obtain 
\begin{equation*}
\Phi_{\gamma}(\tau,\tau')
\in 
2^{8}
(1-p^{-1/2}q^{1/2})(1-p^{1/2}q^{-1/2}) p^{1/2}q^{1/2} (1 + {\mathfrak m})
= 
-2^{8}(p^{1/2}-q^{1/2})^{2} (1 + {\mathfrak m}). 
\end{equation*}
This proves the first case. 
The second case is shown similarly. 
\end{pf}

\subsection{The leading term of $\prod_{\gamma\, \text{odd}}\Phi_{\gamma}$}
\label{subsec:leading:level:1:and:2}
\begin{proposition}
\label{prop:odd:leading:level:1:and:2}
The following holds:
$$
\prod_{\gamma\,{\rm odd}}\Phi_{\gamma}(\tau,\tau')
\in 2^{16}(p-q)^{2} (1 + {\mathfrak m}). 
$$
\end{proposition}

\begin{pf}
The result follows from Lemmas ~\ref{lemma:odd:leading:level:1} and~\ref{lemma:odd:leading:level:2}.
\end{pf}

\section{Proof of Theorem~\ref{thm:MainTheorem}}
\label{sec:proof}
We first show \eqref{eqn:Main:Theorem} up to a constant, and then determine the constant.

\subsection{The formula \eqref{eqn:Main:Theorem} up to a constant}
Recall that the $j$-invariant $j(\tau)$ is the ${\rm SL}_{2}({\bf Z})$-invariant holomorphic function on ${\mathfrak H}$ defined as 
\[
j(\tau) 
= 
\frac{\left(1+ 240 \sum_{n>0} \sigma_3(n) p^n\right)^3}{p \prod_{n> 0} (1-p^n)^{24}}
=
\frac{1}{p}+744+196884p+\cdots, 
\]
where $p = \exp(2 \pi i \tau)$ and 
$\sigma_3(n) = \sum_{d \vert n} d^3$. Then $Y(1) = \SL_2(\ZZ)\backslash {\mathfrak H}$ is isomorphic to 
$\CC$ and $X(1) = Y(1)^*$ is isomorphic to $\PP^1 = \CC\cup\{\infty\}$ via~$j$.

Let ${\rm pr}_{i}\colon X(1)\times X(1)\to X(1)$ be the projection to the $i$-th factor.
Under the identification $X(1) \cong \PP^1$ via $j$, we get the equality of divisors on $X(1)\times X(1)$
$$
{\rm div}({\rm pr}_{1}^{*}j - {\rm pr}_{2}^{*}j)=\varDelta_{X(1)\times X(1)} - (\{\infty\}\times\PP^1) -(\PP^1\times \{\infty\}).
$$

Recall from Sect.~\ref{sect:invol:odd:Borcherds} that 
$\Sigma^{2}X(1) \colonequals (X(1)\times X(1))/{\mathfrak S}_2 \cong \PP^2$, that 
$\varDelta$ is the image of the diagonal of $X(1)\times X(1)$ in $\Sigma^{2}X(1)$,  
and that $B$ is the line at infinity of $\Sigma^{2}X(1) \cong  \PP^2$. 
Since the projection $X(1) \times X(1) \to \Sigma^{2}X(1)$ is a double covering with ramification divisor $\varDelta$, we get 
$$
{\rm div}({\rm pr}_{1}^{*}j - {\rm pr}_{2}^{*}j)^{2}=\varDelta-2B
$$ 
on $\Sigma^{2}X(1)$. Thus we get the following equation of currents on $\Sigma^{2}X(1)$:
\begin{equation}
\label{eqn:5:2}
-dd^{c}\log
\left|
({\rm pr}_{1}^{*}j - {\rm pr}_{2}^{*}j)^{2}
\right|^{2}
=
-\delta_{\varDelta}+2 \delta_{B}.
\end{equation}

On the other hand, since the Dedekind $\eta$-function $\eta(\tau)$ is a modular function of half weight on ${\mathfrak H}$, we have 
\[
 -dd^c \log \Vert \eta(\tau)^{24} \Vert^2 = 
 12 \widetilde{\omega_{Y(1)}} - \delta_{\infty}
\]
as a current on $X(1)$, where $\widetilde{\omega_{Y(1)}}$ is the extension of the K\"ahler form 
of $Y(1)$ induced from the Poincar\'e metric on ${\mathfrak H}$. It follows from Theorem~\ref{prop:c:-3} that 
as a current on $\Sigma^{2}X(1)$ we have 
\[
-dd^{c}\log
(
\prod_{\gamma\,{\rm even}}\|\Phi_{\gamma}\|^{2}
)
=
-dd^{c}\log
(
\prod_{\langle J\rangle \neq \binom{1 2 3}{4 5 6}}\overline{\varpi}_{\langle J\rangle}^{*}\|\Phi\|^{2}
)
= 36  \widetilde{\omega_{\Sigma^2 Y(1)}} - 3 \delta_B.  
\]
Together with Proposition~\ref{prop:a:b}, 
we get the following equation of currents on $\Sigma^{2}X(1)$:
\begin{equation}
\label{eqn:5:1}
-dd^{c}\log
\left[
\frac{(\prod_{\gamma\,{\rm odd}}\|\Phi_{\gamma}\|^{2})^{3}}{(\prod_{\gamma\,{\rm even}}\|\Phi_{\gamma}\|^{2})^{2}}
\right]
=
-dd^{c}\log
\left[
\frac{(\prod_{\varrho\in{\mathfrak S}_{3}}\overline{\varpi}_{\varrho}^{*}\|\Phi\|^{2})^{3}}{(\prod_{\langle J\rangle}\overline{\varpi}_{\langle J\rangle}^{*}\|\Phi\|^{2})^{2}}
\right]
=
-3\,\delta_{\varDelta}+ 6\,\delta_{B}.
\end{equation}
Combining \eqref{eqn:5:2} and \eqref{eqn:5:1}, we get the equation of currents on $\Sigma^{2}X(1)$
\begin{equation}
\label{eqn:current:3}
-dd^{c}\log
\left[
\left|
\frac{\prod_{\gamma\,{\rm odd}}\Phi_{\gamma}^{3}}{\prod_{\gamma\,{\rm even}}\Phi_{\gamma}^{2}}
\right|^{2}
\cdot 
\left|
({\rm pr}_{1}^{*}j - {\rm pr}_{2}^{*}j)^{2}
\right|^{-6}
\right]
=
0.
\end{equation}
Since $\Sigma^{2}X(1)$ is compact, there exists by \eqref{eqn:current:3} a nonzero constant $C$ such that
\begin{equation}
\label{eqn:final:c}
{\prod_{\gamma\,{\rm odd}}\Phi_{\gamma}^{6}}/{\prod_{\gamma\,{\rm even}}\Phi_{\gamma}^{4}}
=
C\,\left({\rm pr}_{1}^{*}j - {\rm pr}_{2}^{*}j\right)^{12}. 
\end{equation}

\subsection{Determination of the constant}
Let $a(n)$ be the $n$-th Fourier coefficient of $j(\tau)-744$. 
As before, we put $p = \exp(2 \pi i \tau)$ and $q = \exp(2 \pi \tau')$. 
The denominator formula for the Monster Lie algebra \cite[Lemma~7.1]{Borcherds92} states that 
\begin{equation}
\label{eqn:denominator:formula:monster:Lie:algebra}
 j(\tau) - j(\tau') 
   =
  \left(p^{-1}-q^{-1}\right)
  \prod_{m,n> 0}(1 - p^m q^n)^{a(mn)}
\end{equation}
for all $\tau,\tau'\in{\mathfrak H}$. 
Since $\eta(\tau) = p^{1/24} \prod_{n > 0}(1 - p^n)$, 
we get by Theorem~\ref{prop:c:-3} 
\[
\prod_{\gamma\,{\rm even}}\Phi_{\gamma}(\tau,\tau')^{2} 
= 
\prod_{\langle J\rangle\not=\binom{123}{456}}\Phi_{\ell(J)}\left( j_{\ell(J)}^{-1}\varpi_{\langle J\rangle}(\tau,\tau') \right)^2
\in 
2^{96} (p q)^6 (1 + {\mathfrak m}).  
\]
Together with Proposition~\ref{prop:odd:leading:level:1:and:2}, we get 
\begin{align}
\label{eqn:leading:term:ratio:Phi:1}
& \left(p^{-1}-q^{-1}\right)^{-12}
\frac{\prod_{\gamma\,{\rm odd}}\Phi_{\gamma}(\tau,\tau')^{6}}{\prod_{\gamma\,{\rm even}}\Phi_{\gamma}(\tau,\tau')^{4}}
\\
\notag
& \quad \equiv
2^{-96} 
\left(p^{-1}-q^{-1}\right)^{-12}
\frac{(p-q)^{12}}{(pq)^{12}} (1+{\mathfrak m})
= 
2^{-96}
\mod
{\mathfrak m}.
\end{align}
Comparing \eqref{eqn:denominator:formula:monster:Lie:algebra} and
\eqref{eqn:leading:term:ratio:Phi:1}, we get
\begin{equation}
\label{eqn:leading:term:ratio:Phi:2}
\left(j(\tau)-j(\tau')\right)^{-12}
\frac{\prod_{\gamma\,{\rm odd}}\Phi_{\gamma}(\tau,\tau')^{6}}{\prod_{\gamma\,{\rm even}}\Phi_{\gamma}(\tau,\tau')^{4}}
\equiv
2^{-96}
\mod{\mathfrak m}.
\end{equation}
By \eqref{eqn:final:c}, \eqref{eqn:leading:term:ratio:Phi:2}, we get $C=2^{-96}$.
This completes the proof of Theorem~\ref{thm:MainTheorem}.
\qed

\section{Automorphic forms arising from non-Enriques involutions}
\label{sec:comparison:2}
In the preceding sections, we have discussed the pullbacks of the Borcherds $\Phi$-functions attached to the fixed-point-free 
involutions on the Kummer surfaces ${\rm Km}(E\times E')$. In this section, we study the relation of these automorphic forms
and another automorphic form canonically attached to another anti-symplectic involution closely related to the $K3$ surfaces
with period lattice ${\mathbb K}$. We also discuss the same problem for the family of Enriques surfaces with period lattice
${\mathbb U}\oplus{\mathbb U}(2)$.

\subsection{The case of ${\mathbb U}(2)\oplus{\mathbb U}(2)$}
\par
Recall that ${\rm Km}(E\times E')$ carries an anti-symplectic involution $\mu\colon {\rm Km}(E\times E') \to {\rm Km}(E\times E')$,
whose fixed locus consists of $8$ smooth rational curves and whose anti-invariant lattice is isometric to
${\mathbb K}={\mathbb U}(2)\oplus{\mathbb U}(2)$ (see Section~\ref{sect:period:disc:locus}). 
By \cite[Sect.\,10]{MaYoshikawa20}, one can canonically attach to the anti-symplectic involution $\mu$
a Borcherds product $\Psi_{\mathbb K}$ on the modular variety $\Omega_{\mathbb K}/O({\mathbb K})$, whose weight is $24$ and which
is nowhere vanishing on $\Omega_{\mathbb K}$ by \cite[Th.\,8.1]{Yoshikawa13}. 
More precisely, $\Psi_{\mathbb K}$ is the Borcherds lift of the vector valued modular form for ${\rm Mp}_{2}({\mathbf Z})$ induced from
$\eta(\tau)^{-8}\eta(2\tau)^{8}\eta(4\tau)^{-8}\theta_{{\mathbb A}_{1}^{+}}(\tau)^{8}$ with values in the group ring ${\mathbf C}[A_{\mathbb K}]$.
(See \cite[Def.\,7.6]{Yoshikawa13} for an explicit formula for the ${\mathbf C}[A_{\mathbb K}]$-valued elliptic modular form 
$F_{\mathbb K}(\tau)$ induced from $\eta(\tau)^{-8}\eta(2\tau)^{8}\eta(4\tau)^{-8}\theta_{{\mathbb A}_{1}^{+}}(\tau)^{8}$.)
We normalize $\Psi_{\mathbb K}$ as follows. The Baily-Borel compactification of $\Omega_{\mathbb K}/O({\mathbb K})$ has a unique 
zero-dimensional cusp. By Borcherds \cite[Th.\,13.3]{Borcherds98}, $\Psi_{\mathbb K}$ is expressed as follows near 
the zero-dimensional cusp:
$$
\Psi_{\mathbb K}(z) = q_{1}^{2}q_{2}^{2}(1+{\mathfrak m}),
$$
where $q_{k}=\exp(2\pi i z_{k})$ $(k=1,2)$ and ${\mathfrak m}=q_{1}{\mathbf C}\{q_{1},q_{2}\}+q_{2}{\mathbf C}\{q_{1},q_{2}\}$.
By this expression and the fact that $\Psi_{\mathbb K}$ is nowhere vanishing on ${\Omega}_{\mathbb K}$, we conclude that 
$$
\Psi_{\mathbb K}(z) = \eta(z_{1})^{48}\eta(z_{2})^{48}.
$$
For any $\gamma\in A_{\mathbb K}$, the ratio $\Phi_{\gamma}^{6}/\Psi_{\mathbb K}$ is a rational function on 
$\Omega_{\mathbb K}/O({\mathbb K})$. When $\gamma$ is even, $\Phi_{\gamma}^{6}/\Psi_{\mathbb K}$ is given as follows
by \cite[Cor.\,7.6]{KawaguchiMukaiYoshikawa18}:
$$
\frac{\theta_{\epsilon_{1}}(\tau_{1})^{16}}{\theta_{\epsilon_{2}}(\tau_{1})^{16}}
\cdot
\frac{\theta_{\epsilon_{1}}(\tau_{1})^{16}}{\theta_{\epsilon_{3}}(\tau_{1})^{16}}
\cdot
\frac{\theta_{\delta_{1}}(\tau_{2})^{16}}{\theta_{\delta_{2}}(\tau_{2})^{16}}
\cdot
\frac{\theta_{\delta_{1}}(\tau_{2})^{16}}{\theta_{\delta_{3}}(\tau_{2})^{16}},
$$ 
where $\{ \epsilon_{1}, \epsilon_{2}, \epsilon_{3} \} = \{ \delta_{1}, \delta_{2}, \delta_{3} \} = \{0,2,3\}$. 
Since $\theta_{\epsilon_1}(\tau)^{8}/ \theta_{\epsilon_2}(\tau)^{4}\theta_{\epsilon_3}(\tau)^{4}$ is one of the functions 
$$
\lambda(\tau)(\lambda(\tau)-1), 
\qquad
\lambda(\tau)/(\lambda(\tau)-1)^{2},
\qquad
(\lambda(\tau)-1)/\lambda(\tau)^{2},
$$
the product of any two of these functions on ${\mathfrak H}\times{\mathfrak H}$ is a Borcherds product. This gives a slightly
more geometric explanation for some formulas of Yang--Yin--Yu \cite[Prop.\,3.3]{YangYinYu18}.

\subsection{The case of ${\mathbb U}\oplus{\mathbb U}(2)$}
\par
Set 
$$
{\mathbb K}':={\mathbb U}\oplus{\mathbb U}(2).
$$
The Enriques surfaces with period lattice ${\mathbb K}'$ were studied by Horikawa \cite{Horikawa78}, Barth--Peters \cite{BP83}, 
Mukai--Namikawa \cite{MukaiNamikawa84} and Mukai \cite[Appendix A]{Mukai10}. 
In this subsection, we study the relation between the two Borcherds products naturally attached to this two-parameter family of 
Enriques surfaces and the Weber function. 
The universal covering $K3$ surface of such an Enriques surface carries a non-symplectic involution,
say $\theta$, whose fixed locus consists of a smooth elliptic curve and $8$ smooth rational curves. Then ${\mathbb K}'$ is isometric 
to the anti-invariant sublattice of the $\theta$-action on the second integral cohomology group of the $K3$ surface. 
Considering the two distinct involutions, one Enriques and the other
non-Enriques, we obtain two Borcherds products on the modular variety $\Omega_{{\mathbb K}'}/O({\mathbb K}')$
by \cite[Th.\,10.2]{MaYoshikawa20}. 
One is the restriction of the Borcherds $\Phi$-function to $\Omega_{{\mathbb K}'}$ induced by the natural inclusion 
of lattices ${\mathbb K}'\hookrightarrow\LAM={\mathbb K}'\oplus{\mathbb E}_{8}(2)$, which we denote $\Phi_{{\mathbb K}'}$. 
Since ${\mathbb E}_{8}(2)$ has no roots,
as an automorphic form on $\Omega_{{\mathbb K}'}$, the divisor of $\Phi_{{\mathbb K}'}$ is given by the discriminant divisor
of $\Omega_{{\mathbb K}'}$. 
Comparing the infinite product expansion near the two zero-dimensional cusps of $\Omega_{{\mathbb K}'}/O({\mathbb K}')$, 
one level $1$ and the other level $2$,
$\Phi_{{\mathbb K}'}$ is the Borcherds lift of the ${\mathbf C}[A_{{\mathbb K}'}]$-valued modular form for ${\rm Mp}_{2}({\mathbf Z})$ 
induced from $\eta(\tau)^{-8}\eta(2\tau)^{8}\eta(4\tau)^{-8}\theta_{{\mathbb E}_{8}^{+}(2)}(\tau)^{8}$. 
(See e.g. \cite[Prop.\,7.1]{Yoshikawa13} for the construction of such a ${\mathbf C}[A_{{\mathbb K}'}]$-valued elliptic modular from.)
By \cite[Th.\,13.3]{Borcherds98}, the infinite product expansions of $\Phi_{{\mathbb K}'}$ at the zero-dimensional cusps are given as follows. 
As before, $\{{\mathbf e}_{1}, {\mathbf f}_{1}\}$ (resp. $\{{\mathbf e}_{2}, {\mathbf f}_{2}\}$) is the standard basis of
${\mathbb U}$ (resp. ${\mathbb U}(2)$). Once the cusp is chosen, $\Omega_{{\mathbb K}'}$ is identified with 
${\mathfrak H}\times{\mathfrak H}$ via the map analogous to \eqref{eqn:tube:domain:level2}. At the level $2$ cusp, 
\begin{equation}
\label{eqn:Phi:K':level2}
\begin{aligned}
\Phi_{{\mathbb K}'}(z) 
&= 
2^{8}e^{2\pi i\langle {\mathbf e}_{1}, z_{1}{\mathbf e}_{1}+z_{2}{\mathbf f}_{1}\rangle}
\left(1-e^{2\pi i \langle -{\mathbf e}_{1}+{\mathbf f}_{1},z_{1}{\mathbf e}_{1}+z_{2}{\mathbf f}_{1}\rangle} \right) (1+{\mathfrak m})
\\
&=
2^{8}(q_{2}-q_{1})(1+{\mathfrak m}),
\end{aligned}
\end{equation}
where $q_{k}=\exp(2\pi i z_{k})$ and ${\mathfrak m} = q_{1}{\mathbf C}\{ q_{1}, q_{2} \} + q_{2}{\mathbf C}\{ q_{1}, q_{2} \}$.
At the level $1$ cusp,
\begin{equation}
\label{eqn:Phi:K':level1}
\Phi_{{\mathbb K}'}(w) = 1+ {\mathfrak m}',
\end{equation}
where $q'_{k}=\exp(2\pi i w_{k})$ and ${\mathfrak m}' = q'_{1}{\mathbf C}\{ q'_{1}, q'_{2} \} + q'_{2}{\mathbf C}\{ q'_{1}, q'_{2} \}$.
\par
The automorphic form attached to $\theta$ is given by \cite[Th.\,8.1]{Yoshikawa13}. 
This automorphic form is denoted by $\Psi_{{\mathbb K}'}$, which is the Borcherds lift of the ${\mathbf C}[A_{{\mathbb K}'}]$-valued 
modular form for ${\rm Mp}_{2}({\mathbf Z})$ induced from
$\eta(\tau)^{-8}\eta(2\tau)^{8}\eta(4\tau)^{-8}\theta_{{\mathbb A}_{1}^{+}}(\tau)^{8}$.
See \cite[Def.\,7.6]{Yoshikawa13} for an explicit formula for the ${\mathbf C}[A_{\mathbb K}]$-valued elliptic modular form 
$F_{{\mathbb K}'}(\tau)$ induced from $\eta(\tau)^{-8}\eta(2\tau)^{8}\eta(4\tau)^{-8}\theta_{{\mathbb A}_{1}^{+}}(\tau)^{8}$.
By \cite[Th.\,8.1]{Yoshikawa13}, $\Psi_{{\mathbb K}'}$ has weight $36$ and vanishes of order $1$ on the discriminant divisor 
of $\Omega_{{\mathbb K}'}$. 
Since the Weyl vector of $F_{{\mathbb U}(2)}$ is given by $2{\mathbf e}_{2}+2{\mathbf f}_{2}$ by 
\cite[Th.\,10.4]{Borcherds98}, \cite[p.321 Correction]{Borcherds00} 
and since ${\mathbb U}(2)$ has no roots, at the level $1$ cusp, we get by \cite[Th.\,13.3]{Borcherds98}
\begin{equation}
\label{eqn:Psi:K':level1}
\Psi_{{\mathbb K}'}(w)
=
e^{2\pi i\langle 2{\mathbf e}_{2}+2{\mathbf f}_{2}, w_{1}{\mathbf e}_{2}+w_{2}{\mathbf f}_{1}\rangle} (1+{\mathfrak m}')
=
(q'_{1}q'_{2})^{4} (1+{\mathfrak m}').
\end{equation}
We normalize the ambiguity of the constants of $\Psi_{{\mathbb K}'}$ by the equation \eqref{eqn:Psi:K':level1}.
Since the Weyl vector of $F_{{\mathbb U}}$ is given by $5{\mathbf e}_{1}+4{\mathbf f}_{1}$ by 
\cite[Th.\,10.4]{Borcherds98}, \cite[p.321 Correction]{Borcherds00}, at the level $2$ cusp, we get by \cite[Th.\,13.3]{Borcherds98}
\begin{equation}
\label{eqn:Psi:K':level2}
\begin{aligned}
\Psi_{{\mathbb K}'}(z)
&=
C e^{2\pi i\langle 5{\mathbf e}_{1}+4{\mathbf f}_{1}, z_{1}{\mathbf e}+z_{2}{\mathbf f}_{1}\rangle}
\left(1-e^{2\pi i \langle -{\mathbf e}_{1}+{\mathbf f}_{1},z_{1}{\mathbf e}+z_{2}{\mathbf f}_{1}\rangle} \right) (1+{\mathfrak m})
\\
&=
C (q_{1}q_{2})^{4}(q_{2}-q_{1}) (1+{\mathfrak m}),
\end{aligned}
\end{equation}
where $C$ is an explicit constant determined later.
Set
$$
F:= \Psi_{{\mathbb K}'}/\Phi_{{\mathbb K}'}^{9}.
$$
Then $F$ is a meromorphic function on the Baily--Borel compactification of $\Omega_{{\mathbb K}'}/O({\mathbb K}')$.
We identify $F$ with the corresponding function on the second symmetric product of the compactified modular curve $X_{0}(2)^{*}$.
Let $\Delta$ be the diagonal locus of $\Sigma^{2}X_{0}(2)^{*}$. Let $B_{\infty,\infty}$, $B_{0,0}$, $B_{0,\infty}$ be the boundary 
divisors corresponding to the loci $\{i\infty\}\times X_{0}(2)^{*}+X_{0}(2)^{*}\times\{i\infty\}$, 
$\{0\}\times X_{0}(2)^{*}+X_{0}(2)^{*}\times\{0\}$, $\{0\}\times X_{0}(2)^{*}+X_{0}(2)^{*}\times\{i\infty\}$, respectively. 
We remark that the divisor $B_{\infty,\infty}$ coincides with $B_{0,0}$ in the Baily--Borel compactification of 
$\Omega_{{\mathbb K}'}/O({\mathbb K}')$.
By \eqref{eqn:Phi:K':level2}, \eqref{eqn:Phi:K':level1}, \eqref{eqn:Psi:K':level2}, \eqref{eqn:Psi:K':level1} and the equality of divisors
${\rm div}\,\Phi_{{\mathbb K}'}={\rm div}\,\Psi_{{\mathbb K}'}=\Delta$ on $\Omega_{{\mathbb K}'}/O({\mathbb K}')$, we get the following
equality of divisors on $\Sigma^{2}X_{0}(2)^{*}$ (note that $(q_{1}-q_{2})^{2}$ is the local equation of $\Delta$):
\begin{equation}
\label{eqn:div:F:K'}
{\rm div}\,F = -4\Delta +4B_{\infty,\infty}+4B_{0,\infty}+4B_{0,0}.
\end{equation}
\par
Recall that the Weber function is the modular function for $\Gamma_{0}(2)$ defined as 
$$
W(z):=2^{12}\eta(2z)^{24}/\eta(z)^{24}.
$$
We define
$$
G(z) := \frac{( W(z_{1}) - W(z_{2}) )^{2}}{W(z_{1})W(z_{2})} = \left( \sqrt{\frac{W(z_{1})}{W(z_{2})}} - \sqrt{\frac{W(z_{2})}{W(z_{1})}} \right)^{2}.
$$
Since $[ (1/2, -z_{1}z_{2}, -z_{1},-z_{2}) ] = [ (-1/(2z_{1}), z_{2}, 1, z_{2}/z_{1} ) ]$, 
the coordinates $(z_{1},z_{2})$ and $(w_{1},w_{2})$ are related by the equations $w_{1}=-1/(2z_{1})$, $w_{2}=z_{2}$.
Since $W(\tau)=2^{12}e^{2\pi i \tau}\prod_{n>0}(1+e^{2\pi i n\tau})^{24}$, we see that at the level $2$ cusp, $G$ is expressed as
\begin{equation}
\label{eqn:G:K':level2}
G(z)=\frac{(q_{1}-q_{2})^{2}}{q_{1}q_{2}}(1+{\mathfrak m}).
\end{equation}
Since $W(z_{1})=W(-1/(2w_{1}))=\eta(w_{1})^{24}/\eta(2w_{1})^{24}=(q_{1}')^{-1}(1+{\mathfrak m}')$, 
we see that at the level $1$ cusp, $G$ is expressed as
\begin{equation}
\label{eqn:G:K':level1}
G(z) 
= 
\frac{( W(-1/(2w_{1})) - W(w_{2}) )^{2}}{W(-1/(2w_{1}))W(w_{2})} 
= 
\frac{(q_{1}')^{-2}(1+{\mathfrak m})}{2^{12}(q_{1}')^{-1}q_{2}'(1+{\mathfrak m})}
=
\frac{2^{-12}}{q_{1}'q_{2}'}(1+{\mathfrak m}).
\end{equation}
Since $G$ vanishes of order $1$ on the discriminant divisor, i.e., the diagonal locus $\Delta$, it follows from
\eqref{eqn:G:K':level2}, \eqref{eqn:G:K':level1} that
\begin{equation}
\label{eqn:div:G:K'}
{\rm div}\, G = \Delta - B_{\infty,\infty} - B_{0,\infty} - B_{0,0}.
\end{equation}
By \eqref{eqn:div:F:K'}, \eqref{eqn:div:G:K'}, we see that $FG^{4}$ is a constant function on $\Sigma^{2}X_{0}(2)^{*}$.
Comparing \eqref{eqn:Phi:K':level2}, \eqref{eqn:Psi:K':level2}, \eqref{eqn:G:K':level2}, at the level $2$ cusp, we get
\begin{equation}
\label{eqn:FG4:K':level2}
F(z)G(z)^{4} = 2^{-72}C(1+{\mathfrak m}).
\end{equation}
Similarly, comparing \eqref{eqn:Phi:K':level1}, \eqref{eqn:Psi:K':level1}, \eqref{eqn:G:K':level1}, at the level $1$ cusp, we get
\begin{equation}
\label{eqn:FG4:K':level1}
F(z)G(z)^{4} = 2^{-48}(1+{\mathfrak m}).
\end{equation}
By \eqref{eqn:FG4:K':level2} and \eqref{eqn:FG4:K':level1}, we conclude $FG^{4}=2^{-48}$ and $C=2^{24}$.
As a result, we obtain the following formula:
$$
\Psi_{{\mathbb K}'}(z, \eta_{1^{-8}2^{8}4^{-8}}(\theta_{{\mathbb E}_{8}(2)^{+}}-9\theta_{{\mathbb A}_{1}^{+}}))
=
\frac{\Psi_{{\mathbb K}'}(z)}{\Phi_{{\mathbb K}'}(z)^{4}}
=
2^{-48}\left( \sqrt{\frac{W(z_{1})}{W(z_{2})}} - \sqrt{\frac{W(z_{2})}{W(z_{1})}} \right)^{-8}.
$$
We remark that when $s$ is a divisor of $24$, $W(z_{1})^{24/s}-W(z_{2})^{24/s}$ is expressed as a product of certain Borcherds lifts by
Li--Yang \cite[Th.\,1.8]{LiYang20}.

\medskip
\subsection{Open problems}

We pose some problems, which may merit further study.

\medskip
\noindent
\textit{Problem $1$.}
We study the lattice embeddings ${\mathbb K} \hookrightarrow {\LAM}$ to relate the difference of the $j$-invariants and the Borcherds $\Phi$-function. With other lattice embeddings, does our method  produce other relations between seemingly unrelated modular functions? 
Two particularly interesting cases seem

1) the Enriques structures on very general Jacobian Kummer surfaces, whose classification was given by Ohashi \cite{Ohashi09}, and

2) the Enriques surfaces with cohomologically trivial involutions, studied by Horikawa \cite{Horikawa78}, Barth--Peters \cite{BP83}, Mukai--Namikawa \cite{MukaiNamikawa84} and Mukai \cite[Appendix A]{Mukai10}.

The former is parametrized by Siegel 3-folds and the latter by $X^1(2) \times X^1(2)$ minus the diagonal. 
In \cite{YangYin19}, Yang--Yin showed that the difference of the modular $\lambda$-function is a Borcherds product.
Is it possible to understand this formula in terms of the odd Enriques structures of the Kummer surfaces of product type?
Is it possible to give an explicit formula for the ratio $\Phi_{\gamma}^{2}/\Phi_{\delta}^{2}$ in terms of the $\lambda$-function
for odd $\gamma$ and even $\delta$?

\medskip
\noindent
\textit{Problem $2$.}
Our computation of the restrictions of $\Phi$ for various embeddings $\KK \hookrightarrow \LAM$ is achieved by 
geometric considerations of the period mapping for Enriques surfaces. Is there an alternate (non-geometric) proof 
using techniques of Borcherds products such as Schofer's formula \cite{Schofer09} or Ma's formula \cite{Ma18}? 

\medskip
\noindent
\textit{Problem $3$.}
By Freitag--Salvati-Manni \cite{FreitagSalvatiManni17}, $\Phi$ can be viewed as a theta series, and 
our calculation in Sect.~\ref{sec:comparison} may be viewed as the determination of the leading terms of the theta series of the orthogonal complement of ${\mathbb K}$ in $\LAM$.
It is also natural to expect that $\Phi_{\gamma}$ may be expressed as some type of theta series. 
For even $\gamma$, this is indeed the case by \cite{KawaguchiMukaiYoshikawa18}. 
For odd $\gamma$, we do not know any explicit formula for $\Phi_{\gamma}$.  

\medskip
\noindent
\textit{Problem $4$.}
By \cite{Yo04}, $\Phi$ is the equivariant analytic torsion of $K3$ surfaces with fixed-point-free involution, and 
analytic torsions are firmly tied with Gillet--Soul\'e's arithmetic Riemann--Roch formula \cite{GS}. Is it possible to deduce the formula \eqref{eqn:Main:Theorem} from the arithmetic Riemann--Roch formula or its equivariant extension~\cite{KR}?

\section*{Appendix: Some properties of the Enriques lattice}%
\label{sec:appendix:A}

\renewcommand{\thetheorem}{A.\arabic{theorem}}
\renewcommand{\theclaim}{A.\arabic{theorem}.\arabic{claim}}
\renewcommand{\theequation}{A.\arabic{equation}}
\setcounter{theorem}{0}
\setcounter{equation}{0}

We prove some technical results used in the proof of Lemma~\ref{lemma:boundary:component}.
To simplify the notation, we write $d^{\perp}$ for $d^{\perp_\LAM}$. 
By Lemma~\ref{lemma:isotoropic:sublattices}, we have $d^{\perp}\cong{\mathbb I}_{2,9}(2)$.  

\begin{lemma}
\label{lemma:appendix:surjectivity}
Let $O^{+}(\LAM)_{d}$ be the stabilizer of $d$ in $O^+(\LAM)$. Then 
the restriction map $O^{+}(\LAM)_{d}\ni g \mapsto g|_{d^{\perp}}\in O^{+}(d^{\perp})$ is surjective.
\end{lemma}

\begin{pf}
Recall that for an even $2$-elementary lattice $L$, a vector $\lambda\in L^{\lor}$ is said to be {\em characteristic} 
if $\langle\lambda,x\rangle_{L}\equiv x^{2}\mod{\bf Z}$ for all $x\in L^{\lor}$. By \cite[p.~150]{Nikulin80}, 
a characteristic vector always exists. 
Let $\lambda_{1}\in (\ZZ d)^{\lor}$ and $\lambda_{2}\in(d^{\perp})^{\lor}$ be characteristic vectors of the $2$-elementary lattices ${\bf Z}d$ and $d^{\perp}$, respectively. 
Note that $(\ZZ d)^{\lor} = {\bf Z}(d/2)$, so 
$\lambda_1 \in {\bf Z}(d/2)$. 
The choice of $\lambda_{1}$ (resp. $\lambda_{2}$) is unique up to a vector of ${\bf Z}d$ (resp. $d^{\perp}$).
We define $L \colonequals {\bf Z}(\lambda_{1}+\lambda_{2})+{\bf Z}d\oplus d^{\perp}$.

\medskip
{\it Claim.}\quad 
The lattice $L$ is isometric to $\LAM$. 
\medskip

Since $d/2$ is a characteristic vector of $\ZZ d$, we take $\lambda_{1}=d/2$. 
We see that $(3,-1,\ldots,-1)/2\in{\mathbb I}_{2,9}(2)^{\lor}$ is a 
characteristic vector of ${\mathbb I}_{2,9}(2)$, which we take as $\lambda_{2}$ via the identification of 
$d^\perp \cong {\mathbb I}_{2,9}(2)$.
Then $\lambda_{1}^{2}=-\frac{1}{2}$ and $\lambda_{2}^{2}=\frac{1}{2}$.
Any element of $L$ can be expressed as $x_{1}+\lambda_{1}+x_{2}+\lambda_{2}$, where  $x_{1}\in{\bf Z}d$ and $x_{2}\in d^{\perp}$. 
Since 
\[
(x_{1}+\lambda_{1}+x_{2}+\lambda_{2})^2 = 
(x_{1}+\lambda_{1})^{2}+(x_{2}+\lambda_{2})^{2}\equiv \lambda_{1}^{2}+\lambda_{2}^{2} = 0\mod 2{\bf Z},
\]
$L$ is an even lattice.
Considering the inclusions of lattices ${\bf Z}d\oplus(d^{\perp})\subset L\subset L^{\lor}\subset{\bf Z}(d/2)\oplus(d^{\perp})^{\lor}$,
$A_L$ is a quotient of a subspace of $A_{{\bf Z}d}\oplus A_{d^{\perp}}\cong{\bf F}_{2}^{\oplus12}$. Hence $L$ is $2$-elementary. 
Since ${\rm sign}(d^{\perp})=(2,9)$ and ${\rm sign}({\mathbf Z}d)=(0,1)$, we get ${\rm sign}(L)=(2,10)$. 
Since $L/\{{\bf Z}d\oplus d^{\perp}\}={\bf F}_{2}(\overline{\lambda}_{1}+\overline{\lambda}_{2})\cong{\bf F}_{2}$
and since $b_{{\bf Z}d\oplus(d^{\perp})}(\cdot,\cdot)$ is non-degenerate, 
\[
L^{\lor}/\{{\bf Z}d\oplus d^{\perp} \}
=\{x\in A_{{\bf Z}d\oplus d^{\perp}}\mid b_{{\bf Z}d\oplus d^{\perp}}(x,\overline{\lambda}_{1}+\overline{\lambda}_{2})\equiv0\}
\cong{\bf F}_{2}^{\oplus 11}. 
\]
Thus we get $\rank A_L \colonequals \dim_{{\bf F}_{2}}L^{\lor}/L=10$.
Finally, let $y=y_{1}+y_{2}\in L^{\lor}$ be an arbitrary vector, where $y_{1}\in{\bf Z}(d/2)$ and $y_{2}\in(d^{\perp})^{\lor}$.
Since $y\in L^{\lor}$ and since $\lambda_{1}$ and $\lambda_{2}$ are characteristic vectors, we get
\[
0
\equiv 
b_{{\bf Z}d\oplus(d^{\perp})}(y,\overline{\lambda}_{1}+\overline{\lambda}_{2})
\equiv 
\langle y_{1},\lambda_{1}\rangle+\langle y_{2},\lambda_{2}\rangle
\equiv
y_{1}^{2}+y_{2}^{2} = y^{2}\mod{\bf Z}.
\]
Thus $\delta(L)=0$.
All together, we obtain the claim $L\cong\LAM$ by \cite[Th.~3.6.2]{Nikulin80}.

Let $\gamma \in O(d^{\perp})$ be an arbitrary element. 
We set $g \colonequals 1_{{\bf Z}d}\oplus \gamma\in O({\bf Z}d\oplus d^{\perp})$.
Since, by \cite[Lemma 3.9.1]{Nikulin80}, $1_{{\bf Z}d}$ and $\gamma$ respectively preserve the set of characteristic vectors of ${\bf Z}(d/2)$ and $(d^{\perp})^{\lor}$, 
$g$ preserves $L$. By the above claim, we have $g \in O(\LAM)_{d}$ and $g|_{d^{\perp}}=\gamma$. 
If $\gamma \in O^{+}(d^{\perp})$, then $\gamma$ preserves the connected components of $\Omega_{d^{\perp}}$.
Since $g$ is an extension of $\gamma$, $g$ preserves the connected components of $\Omega_{\LAM}$.
That is, if $\gamma\in O^{+}(d^{\perp})$, then $g\in O^{+}(\LAM)_d$.
\end{pf}

\begin{lemma}
\label{lemma:appendix:extension}
Let $d \in \Delta_\LAM$. 
If $F\subset d^{\perp}$ is a totally isotropic primitive sublattice of rank $2$, then any element of ${\rm SL}(F)$ lifts to an element of $O^{+}(d^{\perp})$.
\end{lemma}

\begin{pf}
Since $d^{\perp_{\LAM}}\cong{\mathbb I}_{2,9}(2)$ by Lemma~\ref{lemma:isotoropic:sublattices}, 
it suffices to prove that, if $F$ is a totally isotropic primitive sublattice of ${\mathbb I}_{2,9}$
of rank $2$, then any element of ${\rm SL}(F)$ lifts to an element of 
$O^{+}({\mathbb I}_{2,9})  = O^{+}\left({\mathbb I}_{2,9}(2)\right)$.
We take an identification ${\mathbb I}_{2,9}\cong{\mathbb U}\oplus{\mathbb U}\oplus{\mathbb I}_{7}(-1)$, 
where ${\mathbb I}_{7}$ is the positive-definite unimodular lattice of rank $7$ given by the Gram matrix ${\bf 1}_7$. 
Let ${\bf e}, {\bf e'}$ (resp. ${\bf f}, {\bf f'}$) be the standard free basis of the left (resp. middle) lattice $\UU$ of ${\mathbb U}\oplus{\mathbb U}\oplus{\mathbb I}_{7}(-1)$. 
By \cite[Prop.~1.17.1]{Nikulin80}, ${\mathbb I}_{2,9}$ has a unique primitive totally isotropic sublattice of rank $2$ up to $O({\mathbb I}_{2,9})$. 
Thus we may assume that 
$F = \ZZ {\bf e} + \ZZ {\bf f}$. Let $\gamma = \binom{a \, b}{c \, d} \in \SL(F)$ be any element. 
Since $O^{+}({\mathbb U}\oplus{\mathbb U})\subset O^{+}({\mathbb U}\oplus{\mathbb U}\oplus{\mathbb I}_{7}(-1))$ in the obvious way, 
it suffices to show that $\gamma$ lifts to an element of $O^+(\UU \oplus \UU)$. We define the lattice homomorphism 
$
g\colon \UU \oplus \UU \to \UU \oplus \UU 
$
as 
\[
  g({\bf e}) = a\, {\bf e} + c\, {\bf f}, \quad  g({\bf f}) = b\, {\bf e} + d\, {\bf f}, \quad 
  g({\bf e'}) = d\, {\bf e'} - b\, {\bf f'}, \quad  g({\bf f'}) = -c\, {\bf e'} + a\, {\bf f'}. 
\]
Then $g \in O(\UU \oplus \UU)$ and $\rest{g}{F} = \gamma$. Further, with the identification 
of $x {\bf e} + y{\bf e'} + z {\bf f} + w {\bf f'} \in \PP((\UU \oplus \UU) \otimes \CC)$ 
with $(x:y:z:w) \in \PP(\CC^{4})$,  
we regard  
\[
\Omega_{\UU \oplus \UU}
= \{ (x:y:z:w) \in \PP(\CC^{4}) \mid x y + zw = 0, \, x \overline{y} + \overline{x} y + z \overline{w} + \overline{z} w > 0\}
\]
and $\Omega_{\UU \oplus \UU}^+  = \Omega_{\UU \oplus \UU}
\cap \{\Ima (z/x) >0, \, \Ima (w/x) >0\}$. 
Then we see that $g$ maps $(1:1:i:i) \in\Omega_{\UU \oplus \UU}^+$ to an element of $\Omega_{\UU \oplus \UU}^+$. 
We conclude that $\gamma$ lifts to $g \in O^+(\UU \oplus \UU)$. 
\end{pf}



\begin{thebibliography}{99}


\bibitem{BP83}
Barth, W., Peters, C.
\newblock 
{\em Automorphisms of Enriques  Surfaces},
\newblock 
Invent. Math.
\newblock 
{\bf 73}
\newblock 
(1983),
\newblock 
383--411.


\bibitem{BPV84}
Barth, W., Hulek, K., Peters, C., Van de Ven, A.
\newblock 
{\em Compact Complex Surfaces},
\newblock
Second edition,
\newblock 
Springer,
\newblock 
Berlin
\newblock 
(2004). 


\bibitem{Borcherds92}
Borcherds, R.E.
\newblock 
{\em Monstrous moonshine and monstrous Lie superalgebras},
\newblock 
Invent. Math.
\newblock 
{\bf 109}
\newblock 
(1992),
\newblock 
405--444.


\bibitem{Borcherds95}
Borcherds, R.E.
\newblock 
{\em Automorphic forms on $O_{s+2, s}(\RR)$ and infinite products},
\newblock 
Invent. Math. 
\newblock 
{\bf 120}
\newblock 
(1995),
\newblock 
161--213.


\bibitem{Borcherds96}
Borcherds, R.E.
\newblock 
{\em The moduli space of Enriques surfaces 
and the fake monster Lie superalgebra},
\newblock 
Topology
\newblock 
{\bf 35}
\newblock 
(1996),
\newblock 
699--710.

\bibitem{Borcherds98}
Borcherds, R.E.
\newblock 
{\em Automorphic forms with singularities on Grassmannians},
\newblock 
Invent. Math. 
\newblock 
{\bf 132}
\newblock 
(1998),
\newblock 
491--562.


\bibitem{Borcherds00}
Borcherds, R.E.
\newblock 
{\em Reflection groups of Lorentzian lattices},
\newblock 
Duke Math. J.
\newblock 
{\bf 104}
\newblock 
(2000),
\newblock 
319--366.


\bibitem{Borel72}
Borel, A.
\newblock 
{\em Some metric properties of arithmetic quotients of symmetric spaces and an extension theorem},
\newblock 
J. Differential Geom.
\newblock 
{\bf 6}
\newblock 
(1972),
\newblock 
543--560.


\bibitem{FreitagSalvatiManni17}
Freitag, E., Salvati Manni, R.
\newblock 
{\em Octavic theta series},
\newblock 
Asian J. Math.
\newblock 
{\bf 21}
\newblock 
(2017),
\newblock
483--498.


\bibitem{GS}
Gillet, H., Soul\'e, C.
\newblock
{\em An arithmetic Riemann--Roch theorem},
\newblock
Invent. Math.
\newblock
{\bf 110}
\newblock
(1992),
\newblock 
473--543.

\bibitem{GHS}
Gritsenko, V., Hulek, K., Sankaran, G. K.
\newblock
{\em Abelianisation of orthogonal groups and the fundamental group of modular varieties},
\newblock
J. Algebra
\newblock
{\bf 322}
\newblock
(2009),
\newblock 
463--478.


\bibitem{GZ}
Gross, B.H., Zagier, D.B. 
\newblock
{\em On singular moduli},
\newblock 
J. reine angew. Math. 
\newblock 
{\bf 355}
\newblock 
(1985),
\newblock 
191--220.

\bibitem{Horikawa78}
Horikawa, E.
\newblock
{\em On the periods of Enriques surfaces. I},
\newblock 
Math. Ann.
\newblock
{\bf 234}
\newblock
(1978)
\newblock
73--88.


\bibitem{KawaguchiMukaiYoshikawa18}
Kawaguchi, S., Mukai, S., Yoshikawa, K.-I.
\newblock
{\em Resultants and the Borcherds $\Phi$-function},
\newblock
Amer. J. Math. 
\newblock
{\bf 140}
\newblock
(2018),
\newblock 
1471--1519. 


\bibitem{KR}
K\"ohler, K., Roessler, D. 
\newblock
{\em A fixed point formula of Lefschetz type in Arakelov geometry. I. Statement and proof.}
\newblock
Invent. Math. 
\newblock
{\bf 145}
\newblock
(2001),
\newblock 
333--396.


\bibitem{Kondo86}
Kondo, S.
\newblock
{\em Enriques surfaces with finite automorphism groups},
\newblock 
Japan. J. Math.
\newblock
{\bf 12}
\newblock
(1986),
\newblock 
191--282.

\bibitem{LiYang20}
\newblock
Li, Y., Yang, T.
\newblock
{\em On a conjecture of Yui and Zagier},
\newblock
Algebra and Number Theory
\newblock
{\bf 14}
\newblock
(2020),
\newblock
2197--2238.



\bibitem{Ma18}
Ma, S.
\newblock
{\em Quasi-pullback of Borcherds products},
Bull. London Math. Soc.
\newblock
{\bf 51}
\newblock
(2019),
\newblock
1061--1078.


\bibitem{MaYoshikawa20}
Ma, S., Yoshikawa, K.-I.
\newblock
{\em K3 surfaces with involution, equivariant analytic torsion, and automorphic forms on the moduli space IV}, 
\newblock
Compositio Math.
\newblock
{\bf 156}
\newblock
(2020),
\newblock 
1965--2019. 


\bibitem{Mukai10}
Mukai, S.
\newblock
{\em Numerically trivial involutions of Kummer type of an Enriques surface},
\newblock 
Kyoto J. Math.
\newblock
{\bf 50}
\newblock
(2010)
\newblock
889--904.


\bibitem{MukaiNamikawa84}
Mukai, S., Namikawa, Y.
\newblock
{\em Automorphisms of Enriques surfaces which act trivially on the cohomology groups},
\newblock 
Invent. Math.
\newblock
{\bf 77}
\newblock
(1984),
\newblock 
383--397.

\bibitem{Namikawa85}
Namikawa, Y.
\newblock 
{\em Periods of Enriques surfaces},
\newblock 
Math. Ann.
\newblock 
{\bf 270}
\newblock 
(1985),
\newblock 
201--222.


\bibitem{Nikulin80}
Nikulin, V.V.
\newblock 
{\em Integral Symmetric bilinear forms and some of their applications},
\newblock 
Math. USSR Izv.
\newblock 
{\bf 14}
\newblock 
(1980)
\newblock 
103--167.


\bibitem{Ohashi07}
Ohashi, H.
\newblock
{\em On the number of Enriques quotients of 
a $K3$ surface},
\newblock 
Publ. RIMS, Kyoto Univ.
\newblock
{\bf 43}
\newblock
(2007),
\newblock 
181--200.


\bibitem{Ohashi09}
Ohashi, H.
\newblock
{\em Enriques surfaces covered by Jacobian Kummer surfaces},
\newblock 
Nagoya Math. J.
\newblock
{\bf 195}
\newblock
(2009),
\newblock 
165--186.


\bibitem{Scattone87}
Scattone, F.
\newblock 
{\em On the compactification of moduli spaces for algebraic
$K3$ surfaces},
\newblock 
Memoirs of the AMS
\newblock 
{\bf 70}
\newblock 
374
\newblock
(1987). 


\bibitem{Schofer09}
Schofer, J.
\newblock
{\em Borcherds forms and generalizations of singular moduli},
\newblock
J. reine angew. Math.
\newblock
{\bf 629}
\newblock
(2009),
\newblock
1--36.


\bibitem{Sterk}
Sterk, H.
\newblock
{\em Compactifications of the period space of Enriques surfaces Part I},
\newblock 
Math. Z. 
\newblock
{\bf 207}
\newblock
(1991),
\newblock 
1--36.


\bibitem{YangYinYu18}
Yang, T., Yin, H., Yu, P.
\newblock
{\em The lambda invariants at CM points},
\newblock
IMRN
\newblock
doi:10.1093/imrn/rnz230


\bibitem{YangYin19}
\newblock
Yang, T., Yin, H.
\newblock
{\em Difference of modular functions and their CM value factorization},
\newblock
Trans. Amer. Math. Soc.
\newblock
{\bf 371}
\newblock
(2019),
\newblock
3451--3482.




\bibitem{Yo04}
Yoshikawa, K.-I.
\newblock
{\em K3 surfaces with involution, equivariant analytic torsion, and automorphic forms on the moduli space}, 
\newblock
Invent. Math.
\newblock
{\bf 156}
\newblock
(2004),
\newblock 
53--117. 


\bibitem{Yoshikawa09}
Yoshikawa, K.-I.
\newblock
{\em Calabi--Yau threefolds of Borcea--Voisin, analytic torsion, and Borcherds products}, 
\newblock
Ast\'erisque
\newblock
{\bf 328}
\newblock
(2009),
\newblock 
355--393. 


\bibitem{Yoshikawa13}
Yoshikawa, K.-I.
\newblock
{\em $K3$ surfaces with involution, equivariant analytic torsion, and automorphic forms on the moduli space II},
\newblock
J. reine angew. Math. 
\newblock
{\bf 677}
\newblock
(2013),
\newblock
15--70.



\end{thebibliography}
\end{document}